\documentclass[a4paper, reqno]{amsart}

\usepackage[utf8]{inputenc}	
\usepackage[T1]{fontenc}	% einige Zeichensätze und Sonderzeichen - verpixelt die Schrift, wenn das Paket "cm-super" nicht installiert ist!!!
\usepackage{a4wide}
\usepackage{amsmath}		% ganz viele Mathe-Sachen, z. B. Umgebungen, weitere Zeichen, ...
\usepackage{amssymb}		% einige weitere mathematischen Sonderzeichen
\usepackage{enumerate}	% mehr Möglichkeiten, enumerations (einfach) umzudefinieren
\usepackage{nicefrac}		% 45°-Bruchstriche, z. b. für Brüche im Text oder Faktorgruppen
\usepackage{array}			% mehr Möglichkeiten bei Tabellen und Arrays, z.B. > und <, außerdem optische Verbesserung
\usepackage{calc}				% einfachere Rechnungen, z. B. \setlength{...}{2*(2mm-0.9pt)+1em}
\usepackage{flafter}		% Floats werden nicht vor ihrer Deklarationsstelle gesetzt
\usepackage{esdiff}
\usepackage{pgffor}   %for loops in latex
\usepackage{pgf}
\usepackage{mathabx}
\usepackage{multicol}

\usepackage
[
	pdfauthor   = {Lisa Maria Kreusser},
  pdftitle    = {Title},
  pdfsubject  = {},
  pdfkeywords = {},
  bookmarks=true,
  bookmarksopen=true,
  pagebackref=false,
  hyperindex=true,
  colorlinks=true,
  linkcolor=blue,
  citecolor=blue,
  filecolor=blue,
  urlcolor=blue,
  ]
  {hyperref}						% verlinkt alles, was linkable ist, z.B. Fußnoten, das Verzeichnisse, Referenzen, sowie aktive Links
\usepackage{textcomp}		% Sonderzeichen etc.
\usepackage{color}			% Farben (Treiber?)
\usepackage{graphicx}		% Grafikelemente, z.B. Bilder einfügen, Boxen drehen, ...
\usepackage{amscd}			% kommutative Diagramme

\usepackage{mathrsfs}		% geschwungene Buchstaben, Aufruf über \mathscr{}
\usepackage{graphicx}
\usepackage[caption=false,subrefformat=parens,labelformat=parens]{subfig}
\newcommand\N{\mathbb{N}}
\newcommand\R{\mathbb{R}}

\newcommand\C{\mathbb{C}}

\newcommand\bl{\left(}
\newcommand\br{\right)}

\newcommand*\di{\mathop{}\!\mathrm{d}}

\renewcommand\epsilon{\varepsilon}
\renewcommand\theta{\vartheta}

\newtheoremstyle{mytheoremstyle} % name
    {6pt}                    % Space above
    {6pt}                    % Space below
    {\itshape}                   % Body font
    {}						% Indent amount
    {\bf}                   % Theorem head font
    {}                          % Punctuation after theorem head
    {.5em}                       % Space after thm head: " " = normal interword space; \newline = linebreak
    {}  % Theorem head spec (can be left empty, meaning ânormalâ)
\theoremstyle{mytheoremstyle}
\newtheorem{satz}{Satz}[section]

\newtheorem{proposition}[satz]{Proposition}

\numberwithin{equation}{section}
\usepackage{cleveref}
\title{An Anisotropic Interaction Model for Simulating Fingerprints}
\author[B. D\"{u}ring, C. Gottschlich, S. Huckemann, L. M. Kreusser, C.-B. Sch\"{o}nlieb]{Bertram D\"{u}ring, Carsten Gottschlich, Stephan Huckemann, Lisa Maria Kreusser and Carola-Bibiane Sch\"{o}nlieb}

\begin{document}

\begin{abstract}
Evidence suggests that both the interaction of so-called Merkel cells and the epidermal stress distribution play an important role in the formation of fingerprint patterns during pregnancy. To model the formation of fingerprint patterns in a biologically meaningful way these patterns have to become stationary. For the creation of synthetic fingerprints it is also very desirable that rescaling the model parameters leads to rescaled distances between the stationary fingerprint ridges. Based on these observations, as well as the model introduced by  K\"ucken and Champod we propose a new model for  the formation of fingerprint patterns during pregnancy. In this anisotropic interaction model the interaction forces not only depend on the distance vector between the cells and the model parameters, but additionally on an underlying tensor field, representing a stress field. This dependence on the tensor field leads to  complex, anisotropic  patterns. We study the resulting stationary patterns both analytically and numerically. In particular, we show that fingerprint patterns can be modeled as stationary solutions by choosing the underlying tensor field appropriately. 
\end{abstract}

\maketitle

\section{Introduction}\label{sec:introduction}

Large databases are required for developing, validating and comparing 
the performance of fingerprint indexing and identification algorithms. The goal of these algorithms is to search and find a fingerprint in a  database
(or providing the search result that the query fingerprint is not stored in that database). 
The database sizes for fingerprint identification can vary between several thousand fingerprints e.g.\ watchlists in border crossing scenarios
or hundreds of millions of fingerprints in case of the national biometric ID programme of India.

Clearly, fingerprint identification is of great importance in forensic science and is increasingly used in biometric applications.
Unfortunately, collecting databases of real fingerprints for research purposes is usually very cost-intensive,
requires time and effort, and in many countries, it is constrained 
by laws addressing important aspects such as data protection and privacy.
Therefore, it is very desirable to avoid all these disadvantages 
by simulating large fingerprint databases on a computer. 

The creation of synthetic fingerprint images is of great interest 
to the community of biometric and forensic researchers, as well as practitioners. 
The SFinGe method \cite{SFINGE} has been proposed to this end by Cappelli et al. in 2000.
This method can produce fingerprint images which look realistic enough to deceive 
attendees of a pattern recognition conference,
however, systematic differences between real fingerprints and synthetic images by SFinGE
regarding the minutiae pattern have been found  which allow to distinguish between the two \cite{SeparatingRealFromSynthetic}.
Recently, the realistic fingerprint creator (RFC) \cite{RealisticFingerprintCreator} has been suggested 
to overcome the issue of unrealistic minutiae distributions.
SFinGe and RFC are both based on Gabor filters \cite{GaborFilter} for image creation.
A different approach to fingerprint creation has been introduced by K\"{u}cken and Champod in \cite{Merkel}.
They strive to directly model the process of fingerprint pattern formation as it occurs in nature 
and their approach is inspired by existing knowledge from biology, anatomy and dermatology.
Two commonalities of Gabor filters based and biology-inspired approaches are that both start 
with random initial conditions and both perform changes in an iterative fashion. Kücken and Champod suggest a model describing the formation of fingerprint patterns over time based on the interaction of certain cells and mechanical stress in the epidermis  \cite{fingerprintbiology}.
% In the model by K\"ucken and Champod,  
% Merkel cells and their interactions are interpreted as agents of these changes
% and iterations model the flow of time.

In principle, a nature-inspired model nourishes the hope of producing more realistic fingerprints
and potentially also to gain insights into the process of natural fingerprint pattern formation. 
%In the biological community it has been suggested extensively that the interaction of mechanical stress, trophic factors from incoming nerves and interactions between  certain cells (so-called Merkel cells) lead to the formation of fingerprints \cite{Merkel,champod2016fingerprints}. 
%In \cite{Merkel}  K\"{u}cken and Champod introduced a model describing the formation of fingerprint patterns based on the interaction of Merkel cells and mechanical stress in the epidermis \cite{fingerprintbiology}.
An extensive literature \cite{champod2016fingerprints,Dell1986,fingerprintbiology,Merkel,Moore1986,Morohunfola,Wertheim2011} in the biological community suggests that fingerprint patterns are formed due to the interaction of mechanical stress, trophic factors from incoming nerves and interactions between so-called Merkel cells. Merkel cells are epidermal cells that appear in the volar skin at about the 7th week of pregnancy. From that time onward they start to multiply and organise themselves in lines exactly where the primary ridges arise \cite{Merkel}.  K\"{u}cken and Champod  \cite{Merkel} model this pattern formation process on a  domain $\Omega\subset \R^2$ as the rearrangement of Merkel cells from a random initial configuration into roughly parallel ridges along the lines of smallest compressive stress. However, the resulting patterns do not seem to be stationary which is desirable for describing the formation of fingerprints accurately.

%K\"ucken and Champod considered $N$ particles, representing the cells,
%Merkel cells, also referred to as particles,
% at positions $x_j=x_j(t)\in\R^2,~j=1,\ldots,N,$ at time $t$  whose interactions satisfy an evolution equation of the form \eqref{eq:particlemodel}, equipped with initial data $x_j(0)=x_j^{in},~j=1,\ldots,N$. 
%Here, $F(x_j-x_k,T(x_j))$ denotes the total force that particle $k$ exerts on particle $j$ subject to an underlying stress tensor field $T(x_j)$ at $x_j$, describing the local stress field.
%The dependence on $T(x_j)$ is based on the experimental results in \cite{KIM1995411} where an alignment of the particles along the local stress lines is observed. 
%Hence, the evolution of particle $j$ at location $x_j$ depends on the the local stress tensor field $T(x_j)$. Note that \eqref{eq:particlemodel} results from Newton's second law by neglecting the inertia term. 
%Further note that the value of $\mu$ only depends on the choice of time scaling.

Motivated by their approach we propose a general class of  evolutionary particle models with an anisotropic interaction force in two space dimensions where the form of the interaction forces is inspired by biology.  This class of models can be regarded as an biology-inspired adaptation of the K\"{u}cken-Champod model \cite{Merkel}.
%which describes the formation of fingerprint patterns based on the interaction of Merkel cells and mechanical stress in the epidermis \cite{fingerprintbiology}. 
In \cite{patternformationanisotropicmodel} a generalization  of the K\"{u}cken-Champod model  is proposed and its stationary states are studied both analytically and numerically in the spatially homogeneous case.

In this work, we consider  $N$ interacting particles whose positions $x_j=x_j(t)\in\Omega\subset\R^2$, ${j=1,\ldots,N},$ at time $t$   satisfy
\begin{align}\label{eq:particlemodel}
\frac{\di x_j}{\di t}=\frac{1}{N}\sum_{\substack{k=1\\k\neq j}}^N F(x_j-x_k,T(x_j)),
\end{align}
equipped with initial data $x_j(0)=x_j^{in},~j=1,\ldots,N$. The term  $F(x_j-x_k,T(x_j))$ in \eqref{eq:particlemodel} denotes the force which a particle at position $x_k$ exerts on a particle at position $x_j$. This force depends on an underlying stress tensor field $T(x_j)$ at location $x_j$. The existence of such a tensor field $T(x_j)$  is based on the experimental results in \cite{KIM1995411} where an alignment of the particles along the local stress lines is observed.  We define the tensor field $T(x_j)$ by the  directions of smallest stress  at location $x_j$  by a unit vector field  $s=s(x)\in\R^2$ and introduce a corresponding orthonormal vector field $l=l(x)\in\R^2$, representing the directions of largest stress. Then the force is given by
\begin{align}\label{eq:totalforcenew}
F(d=d(x_j,x_k),T(x_j))=f_s(|d|)(s(x_j)\cdot d)s(x_j)+f_l(|d|)(l(x_j)\cdot d)l(x_j)
\end{align}
for coefficient functions $f_s$ and $f_l$. 

In previous works on the Kücken-Champod model \cite{Merkel} and its generalization \cite{patternformationanisotropicmodel} a dynamical system of ordinary differential equations of the form \eqref{eq:particlemodel} was considered where the force that particle $k$ exerts on particle $j$ is given by
\begin{align}\label{eq:totalforce}
F(d=d(x_j,x_k),T(x_j))=F_A(d,T(x_j))+F_R(d),
\end{align} 
i.e.\ the sum of repulsion and attraction forces, $F_R$ and $F_A$, respectively. Here, the attraction force  depends on the underlying tensor field $T(x_j)$ at $x_j$, modeling  the local stress field.  The matrix $T(x_j)$ encodes the direction of the fingerprint lines at $x_j$, defined by 
\begin{align}\label{eq:tensorfield}
T(x):=\chi s(x)\otimes s(x) +l(x)\otimes l(x), \quad\chi\in[0,1],  
\end{align}
and orthonormal  vector fields $s=s(x)$, $l=l(x)\in\R^2$. 
%The vector fields $s, l$ represent the directions of smallest and largest stress, respectively.
%Denoting the directions of smallest and largest compression of the stress tensor field $T(x_j)$ at $x_j$ by the  $s=s(x_j)\in \R^2$ and $l=l(x_j)\in\R^2$, respectively, the tensor field $T(x_j)$ at $x_j$
For studying the pattern formation with an underlying spatially homogeneous tensor field $T$ producing straight parallel ridges, e.g. $$T=\begin{pmatrix}
1& 0\\ 0 &\chi
\end{pmatrix}, $$   is considered \cite{patternformationanisotropicmodel}.  The repulsion and attraction forces in the K\"ucken-Champod model \cite{Merkel} and its generalization in \cite{patternformationanisotropicmodel} are of the form
\begin{align}\label{eq:repulsionforce}
F_R(d)=f_R(|d|)d
\end{align}
and
\begin{align}\label{eq:attractionforce}
F_A(d=d(x_j,x_k),T(x_j))=f_A(|d|)T(x_j)d,
\end{align}
respectively. Note that the direction of the attraction force $F_A$ and hence also the direction of the total force $F$ are regulated by the parameter $\chi$ in the definition of the tensor field $T$. The parameter $\chi$ introduces an anisotropy to the equation leading to  complex, anisotropic patterns. 

For $\chi=1$ the model \eqref{eq:particlemodel} with interaction forces of the form \eqref{eq:totalforce} for repulsion and attraction force \eqref{eq:repulsionforce} and \eqref{eq:attractionforce} reduces to a gradient flow
\begin{align}\label{eq:standardmodel}
\frac{\di x_j}{\di t}%=-N\nabla_{x_j} E
=\frac{1}{N}\sum_{\substack{k=1\\k\neq j}}^N F(x_j-x_k)
\end{align}
and $F(d)=-\nabla W(d)$ for a radially symmetric interaction potential $W$. The continuum equation associated with the isotropic particle model \eqref{eq:standardmodel} is given by
\begin{align*}%\label{eq:standardmodelmacroscopic}
\rho_t+\nabla\cdot\bl \rho u\br=0,\qquad u=-\nabla W\ast \rho
\end{align*}
where  $u=u(t,x)$ is the macroscopic velocity field. This continuum model, referred to as the aggregation equation has been studied extensively recently, mainly in terms of its gradient flow structure, the blow-up dynamics for fully attractive potentials and the rich variety of steady states, see   \cite{gradientflows,balague_preprint,nonlocalinteraction,Confinement,swarmequilibria,bertozzi2009,bertozzi2012,Canizo2015,Carrillo2016,carrillo2011,Carrillo2012306,Carrillo2012550,Carrillo2016304,Carrillo:2003,Carrillo2006,Fellner2010,Fellner20111436,Li2004,nonlocalInteractionEquation,opac-b1122739,vonBrecht2012,PredictingPatternFormation} and the references therein. There has been a trend recently to connect the microscopic and the macroscopic descriptions via  kinetic modeling, see for instance \cite{Bellomo2012,Carrillo2010,Ha2008} for different kinetic models in swarming, \cite{Fornasier2011,Ha2009} for the particle to hydrodynamics passage and \cite{Karper} for the hydrodynamic limit of a kinetic model.    It seems that not many results are currently available in the field of anisotropies. In \cite{Evers2015,Evers2017} anisotropy is modeled by adding weights to the interaction terms. One can show that the model in \cite{Evers2015,Evers2017} is related to our model if a tensor field $T$ is introduced as the velocity direction.

Fingerprint simulation results are shown for certain model parameters in \cite{Merkel} where the underlying tensor field is constructed based on fingerprint images using the NBIS package from the National Institute of Standards and Technology. However, \cite{Merkel} is purely descriptive, the choice of parameters is not discussed and the model \eqref{eq:particlemodel} was not studied mathematically. The model \eqref{eq:particlemodel} was studied analytically and numerically for the first time in \cite{patternformationanisotropicmodel}. Here, the authors justify why the particles align along the vector field lines $s$ provided the parameter $\chi$ is chosen sufficiently small so that the total force is purely repulsive along $s$. Besides,  the authors investigate the stationary states to the particle model \eqref{eq:particlemodel} for a spatially homogeneous underlying tensor field where the chosen model parameters are consistent with the work of Kücken and Champod in \cite{Merkel}. For the simulation of fingerprints, however, non-homogeneous tensor fields have to be considered, making the analysis of the model significantly more difficult. No analytical results of the long-time behavior of \eqref{eq:particlemodel} for non-homogeneous tensor fields are currently available. Besides, numerical results for the given model parameters and different non-homogeneous tensor fields are shown over time in \cite{patternformationanisotropicmodel} and one can clearly see that the resulting patterns are not stationary. The simulation results for realistic tensor fields for the simulation of fingerprints in \cite{Merkel} seem to be far away from being stationary too. This is illustrated in Figure 9 in \cite{Merkel} where snapshots of the solution are shown for a spatially homogeneous tensor field which should have been parallel lines for steady states. In the biological community, however, it is well-known that fingerprint patterns with their ridge lines and minutiae configuration are determined during pregnancy and remain the same during lifetime (as long as no fingerprint alterations occur). Hence, we are particularly interested in stationary solutions of the system \eqref{eq:particlemodel}.

\begin{figure}[htbp]
	\centering
	\subfloat[Original]{\includegraphics[width=0.24\textwidth]{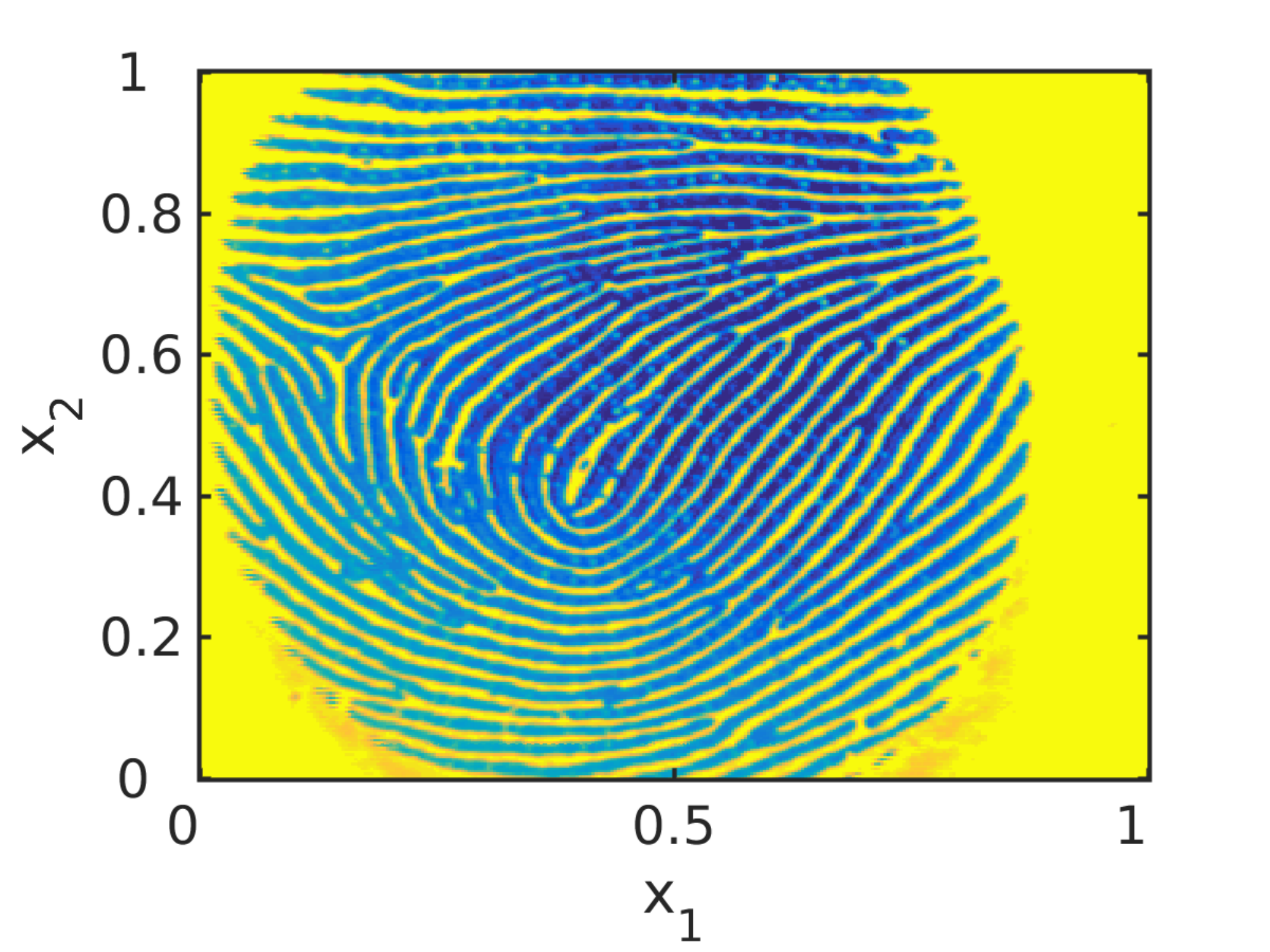}\label{fig:goalorigfingerprint}}\hfill
	%      \subfloat[Arguments $\theta$]{\includegraphics[width=0.24\textwidth]{NumericalResultsFingerprints/FullFingerprintAnglesCut}\label{fig:goalarguments}}\hfill
	\subfloat[$s$ with original]{\includegraphics[width=0.24\textwidth]{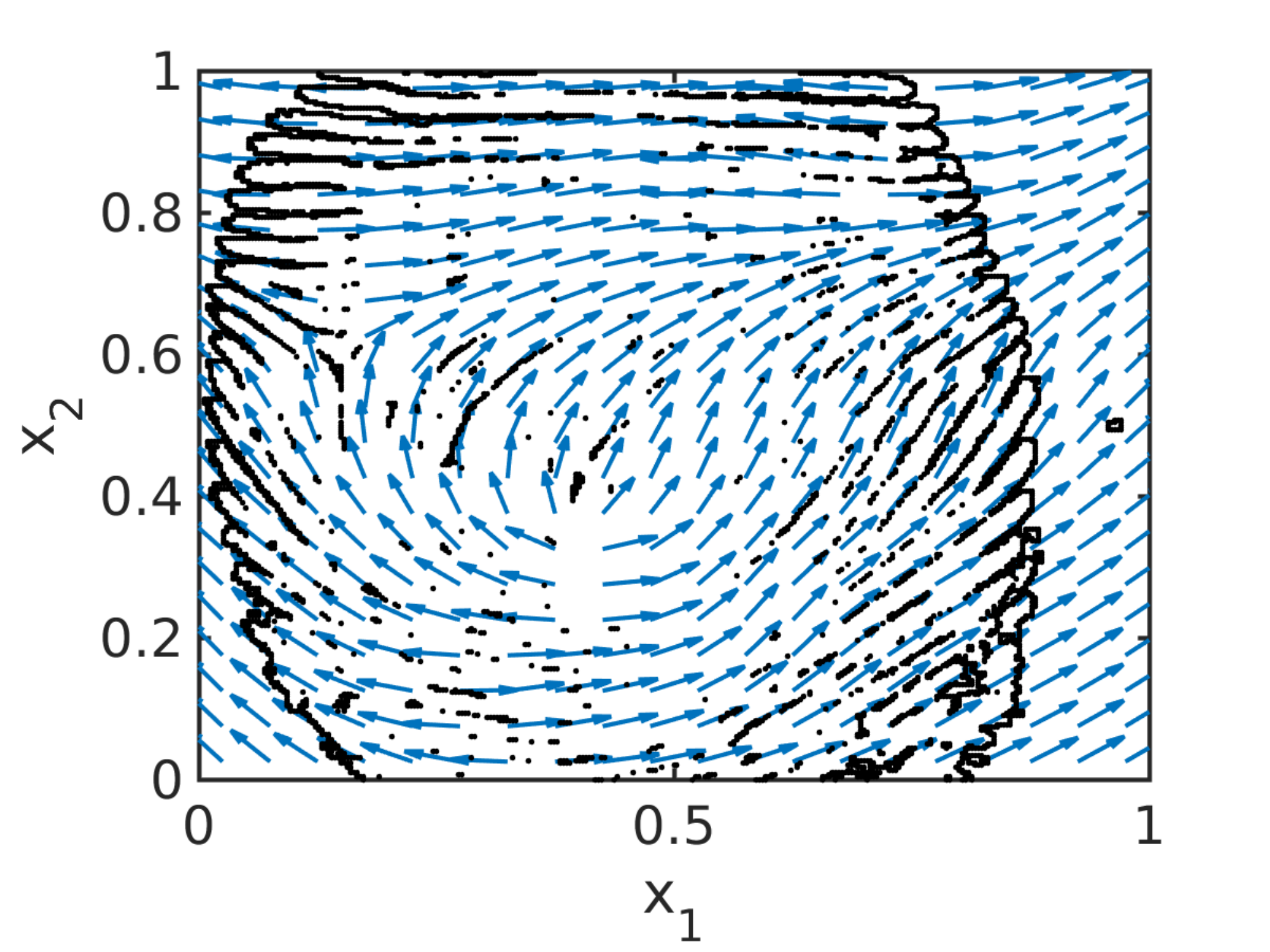}\label{fig:goaltensorfieldmask}}\hfill
	\subfloat[$s$]{\includegraphics[width=0.24\textwidth]
		{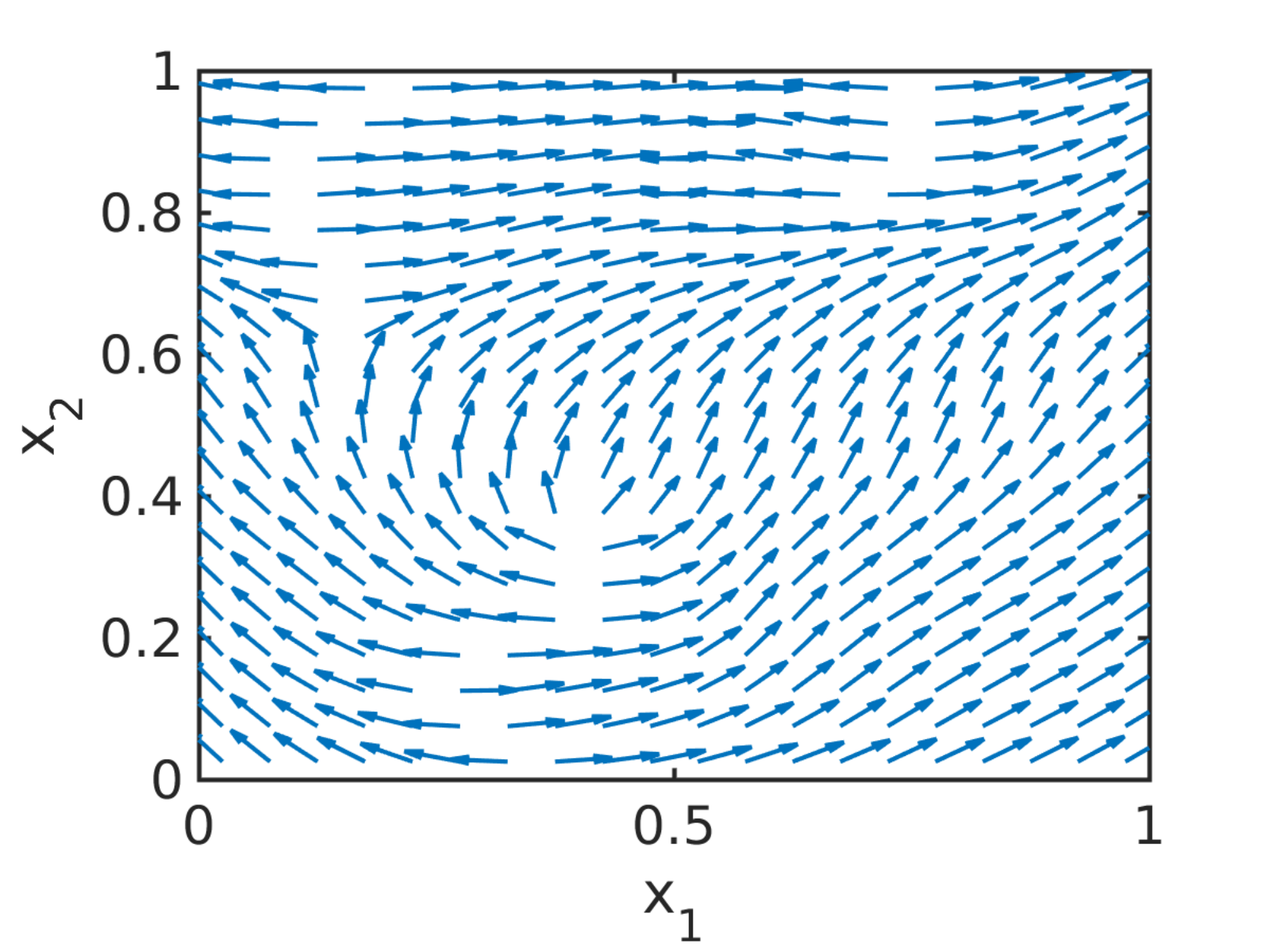}\label{fig:goaltensorfield}}\hfill
	\subfloat[Stationary]{\includegraphics[width=0.24\textwidth]{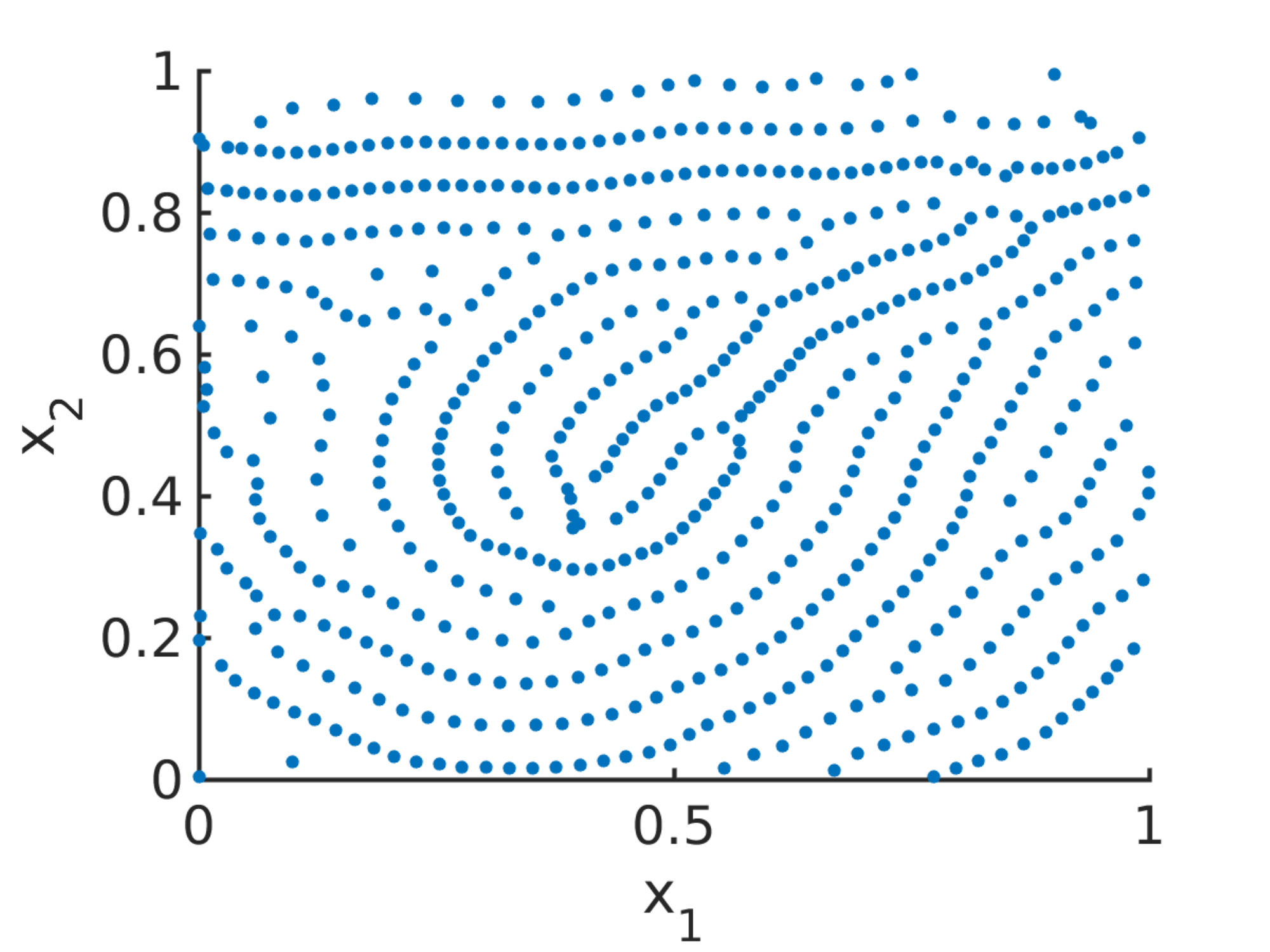}\label{fig:goalstationary}}\hfill
	\caption{Original fingerprint image and lines of smallest stress $s=s(x)$ for the reconstructed tensor field $T=T(x)$ with an overlying mask of the original fingerprint image in black, as well as stationary solution to the interaction model \eqref{eq:particlemodel} for interaction forces of the form \eqref{eq:totalforcenew} and randomly uniformly distributed initial data}\label{fig:stationaryharmonicrealfingercut}
\end{figure}

The goal of this work is to develop an efficient algorithm for creating synthetic fingerprint patterns as  stationary solutions of an evolutionary dynamical system of the form \eqref{eq:particlemodel} as illustrated in Figure \subref*{fig:goalstationary} for the underlying tensor field in Figure \subref*{fig:goaltensorfield}. 
As a first step we study the existence of stationary solutions to \eqref{eq:particlemodel} for spatially homogeneous underlying tensor fields and extend these results to certain spatially inhomogeneous tensor fields. Based on these analytical results as well as the stability analysis of line patterns in \cite{stabilityanalysisanisotropicmodel} we can expect stable stationary patterns along the vector field $s$. Since the solutions to the particle model \eqref{eq:particlemodel} with the parameters suggested by Kücken and Champod do not seem to be stationary, we investigate the impact of the interaction forces on the resulting pattern formation numerically. In particular the size of the total attraction force plays a crucial role in the pattern formation. We adjust the model parameters  accordingly and simulate fingerprints which seem to be close to being stationary, resulting in an extension of the numerical results in \cite{patternformationanisotropicmodel} for inhomogeneous tensor fields. Based on real fingerprint images as in Figure \subref*{fig:goalorigfingerprint} we determine the underlying tensor field $T$ with lines of smallest stress $s$ by extrapolating the direction field outside of the fingerprint image based on \cite{LineSensor}. In Figure \subref*{fig:goaltensorfieldmask} we overlay a mask of the original fingerprint image on the estimated tensor field with direction field $s$ and in Figure \subref*{fig:goaltensorfield} only the direction field $s$ is shown. Besides, we include a novel method for the generation of the underlying tensor fields in our numerical simulations which is based on quadratic differentials as a global model for orientation fields of fingerprints \cite{tensorfieldfingerprints}. 

In the fingerprint community major features of a fingerprint, called minutiae, are of great interest. Examples include  ridge bifurcation, i.e.\ a single ridge dividing into two ridges. We study how they evolve over time, both heuristically and numerically. Finally, we propose a new bio-inspired model for the creation of synthetic fingerprint patterns which  not only allows us to simulate fingerprint patterns as stationary solution of the particle model \eqref{eq:particlemodel} but also adjust the distances between the fingerprint lines by rescaling the model parameters. This is the first step towards modeling fingerprint patterns with specific features in the future.

Studying the model \eqref{eq:particlemodel} and in particular its pattern formation  result in a better understanding of the fingerprint pattern formation process. Due to the generality of the formulation of the anisotropic interaction model \eqref{eq:particlemodel} this can be regarded as an important step towards understanding the formation of fingerprints and may be applicable to other anisotropic interactions in nature.

This work is organized as follows. In Section \ref{sec:descriptionmodel} the  K\"{u}cken-Champod model \cite{Merkel} is introduced and we propose a new bio-inspired modeling approach. Section \ref{sec:existencesteadystates} deals with the existence of steady states to \eqref{eq:particlemodel} for spatially homogeneous tensor fields and its extension to locally spatially homogeneous tensor fields.    In Section  \ref{sec:fingerprintsimulations} we  adapt the  parameters in the force coefficients \eqref{eq:repulsionforcemodel} and \eqref{eq:attractionforcemodel} of the K\"{u}cken-Champod model in such a way that   fingerprint patterns can be obtained as stationary solutions to the particle model \eqref{eq:particlemodel}. Based on these results, we propose the  bio-inspired model, described in Section \ref{sec:descriptionmodel}, to simulate fingerprints with variable  distances between the fingerprint lines. For the creation of realistic fingerprints we consider a novel methods for obtaining the underlying tensor field based on quadratic differentials as well as images of real fingerprint data.

\section{Description of the model}\label{sec:descriptionmodel}

In the sequel, we consider  particle models of the form \eqref{eq:particlemodel} where the force $F$ is of the form \eqref{eq:totalforcenew} or \eqref{eq:totalforce} where the repulsion and attraction forces are given by \eqref{eq:repulsionforce} and \eqref{eq:attractionforce}, respectively. Note that a model of the form \eqref{eq:particlemodel} can be rewritten as
\begin{align}\label{eq:particlemodelrewritten}
\begin{split}
\frac{\di x_j}{\di t}&=v_j\\
v_j&=\frac{1}{N}\sum_{\substack{k=1\\k\neq j}}^N F(x_j-x_k,T(x_j))
\end{split}
\end{align} 
and can be derived from Newton's second law
\begin{align*}
\begin{split}
\frac{\di \overline{x}_j}{\di \tau}=\overline{v}_j\\
m \frac{\di \overline{v}_j}{\di \tau}+\lambda \overline{v}_j=\overline{F}_j
\end{split}
\end{align*}
for particles of identical mass $m$. Here, $\lambda$ denotes the coefficient of friction and $\overline{F}_j$ is the total force acting on particle $j$. Rescaling in time  $\tau=\frac{m}{\epsilon\lambda}t$ for small $\epsilon>0$, setting  $x_j:=\frac{\epsilon\lambda}{m}\overline{x}_j$, $v_j:=\overline{v}_j$ and 
$$\overline{F}_j=\frac{\lambda}{N}\sum_{\substack{k=1\\k\neq j}}^N F(x_j-x_k,T(x_j))$$
results in the rescaled second order model
\begin{align}\label{eq:secondorderparticles}
\begin{split}
\frac{\di x_j}{\di t}&=v_j\\
\epsilon \frac{\di v_j}{\di t}&=-v_j+\frac{1}{N}\sum_{\substack{k=1\\k\neq j}}^N F(x_j-x_k,T(x_j))
\end{split}
\end{align}
for small $\epsilon>0$. 
Starting from \eqref{eq:secondorderparticles} the first order model \eqref{eq:particlemodelrewritten} was justified and formally derived in \cite{Bodnar} and similar to the rigorous limit from the isotropic second order model \eqref{eq:secondorderparticles} to the isotropic first order model \eqref{eq:particlemodelrewritten} as $\epsilon\to 0$  in \cite{Fetecau2015} one can proof the rigorous limit of the anisotropic model \eqref{eq:secondorderparticles}. Note that setting $\epsilon=0$ in  \eqref{eq:secondorderparticles} leads to \eqref{eq:particlemodelrewritten}, corresponding to instantaneous changes in velocities.

\subsection{Kücken-Champod model}
In the papers \cite{patternformationanisotropicmodel,Merkel} systems of evolutionary differential equations of the form \eqref{eq:particlemodel} are considered where the total force, the attraction and the repulsion forces are of the forms \eqref{eq:totalforce}, \eqref{eq:repulsionforce} and \eqref{eq:attractionforce}, respectively, and the underlying tensor field $T$ is defined as \eqref{eq:tensorfield}. The coefficient functions $f_R$ and $f_A$ of the repulsion force $F_R$  \eqref{eq:repulsionforce} and the attraction force \eqref{eq:attractionforce} in the Kücken-Champod model are given by
\begin{align}\label{eq:repulsionforcemodel}
f_R(d)=(\alpha |d|^2+\beta)\exp(-e_R |d|)
\end{align}
and 
\begin{align}\label{eq:attractionforcemodel}
f_A(d)=-\gamma|d|\exp(-e_A|d|)
\end{align}
for nonnegative constants $\alpha$, $\beta$, $\gamma$, $e_A$ and $e_R$, and, again, $d=d(x_j,x_k)=x_j-x_k\in\R^2$.  To be consistent with the work of Kücken and Champod \cite{Merkel} we assume that the total force \eqref{eq:totalforce}  exhibits short-range repulsion and long-range attraction along $l$ and we choose the parameters in an initial study as:
\begin{align}\label{eq:parametervaluesRepulsionAttraction}
\begin{split}
\alpha&=270, \quad \beta=0.1, \quad \gamma=35, \quad
e_A=95, \quad e_R=100, \quad \chi=0.2.
\end{split}
\end{align} 
These parameters are chosen in such a way that the resulting plots of the force coefficients are as close as possible to the ones shown by K\"ucken and Champod in \cite{Merkel}. Here,  the parameter $\chi\in[0,1]$ determines the direction of the interaction. For $\chi=1$ the attraction force between two particles is aligned along their distance vector, while for $\chi=0$ the attraction between two particles is oriented exactly along the lines of largest stress \cite{patternformationanisotropicmodel}. 

\begin{figure}[htbp]
	\centering
	\subfloat[Force coefficients  $f_R$ and  $f_A$]{\includegraphics[width=0.45\textwidth]{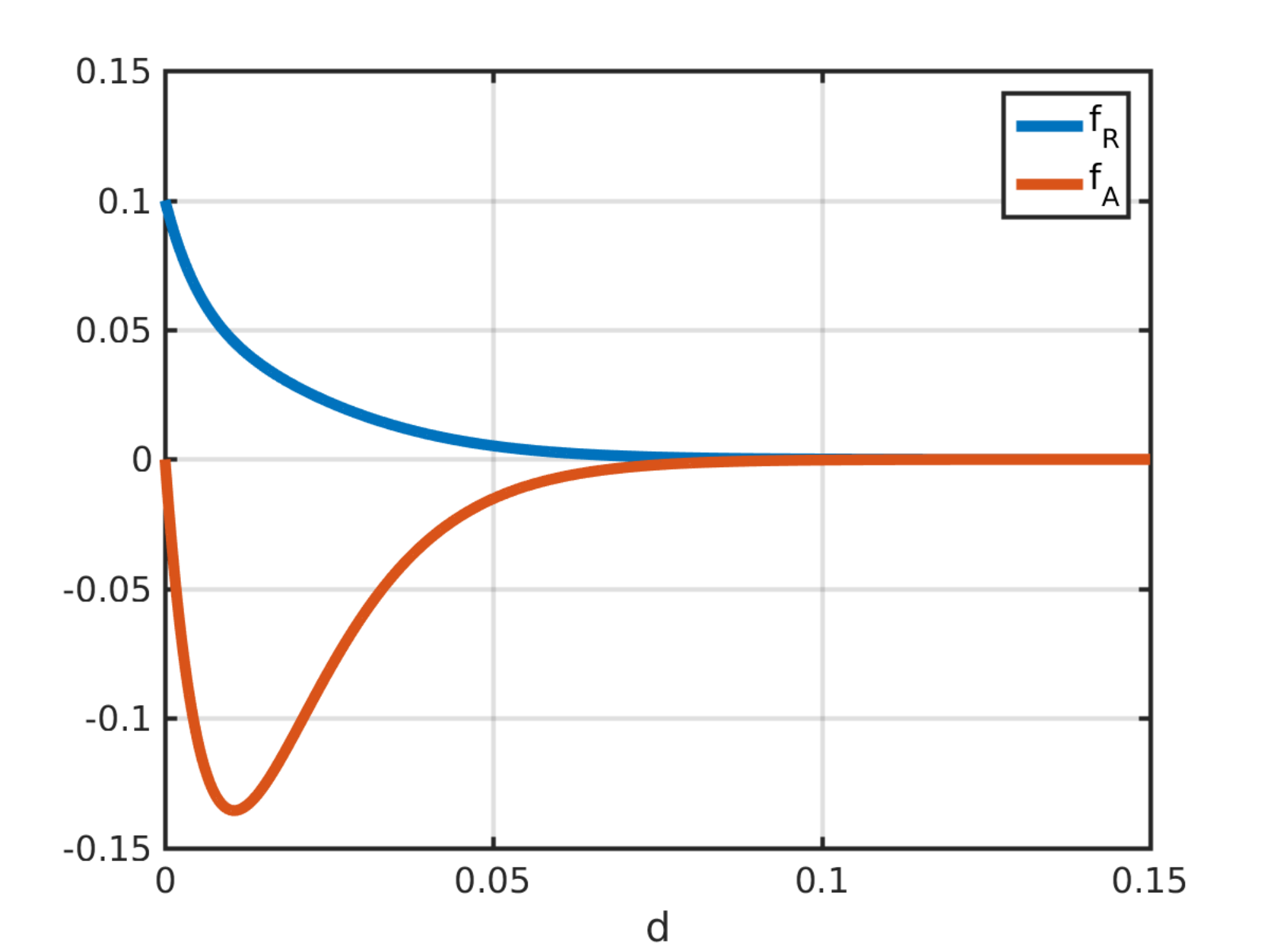}\label{fig:forces}}
	\subfloat[Total force coefficients]{\includegraphics[width=0.45\textwidth]{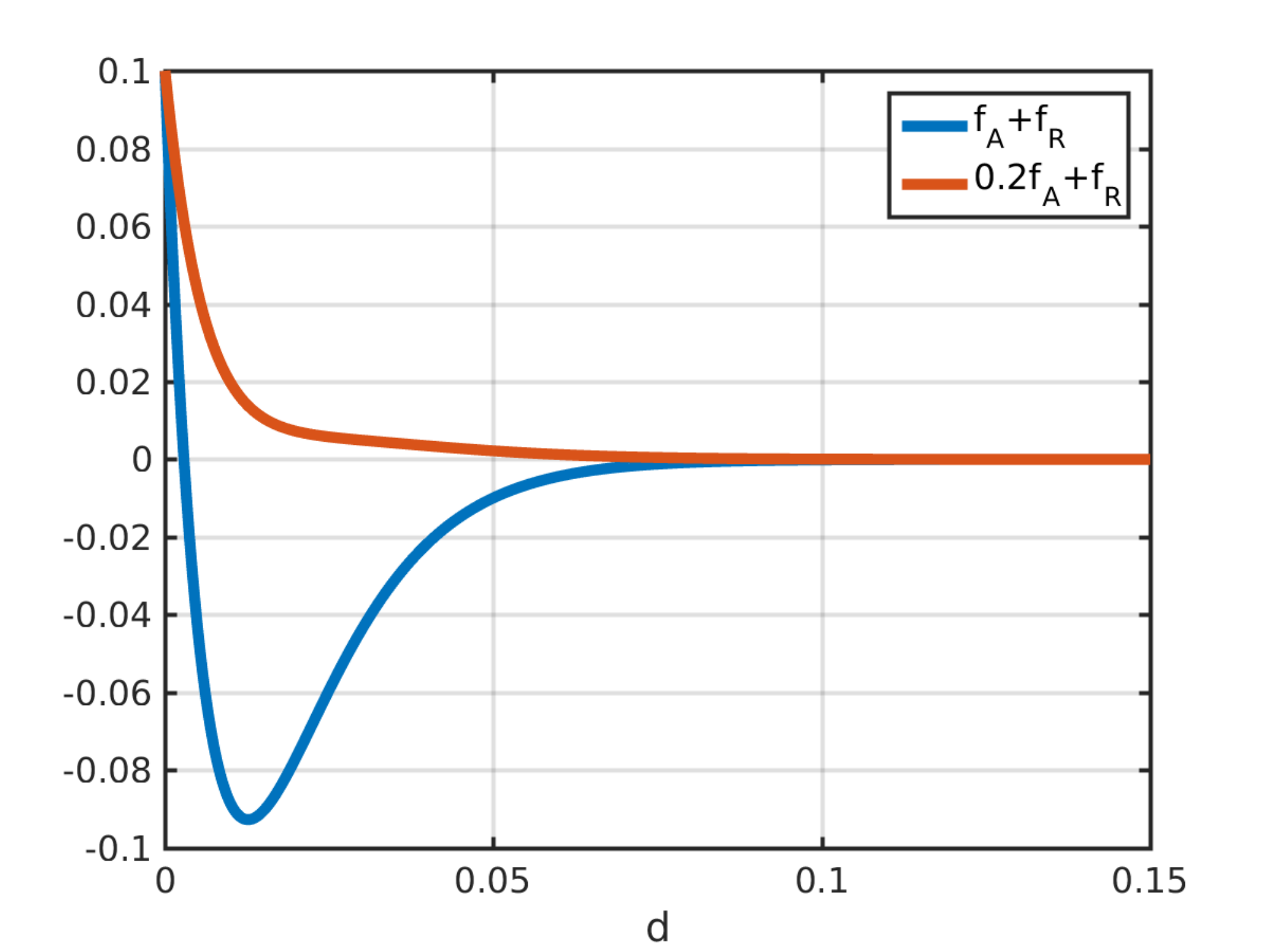}\label{fig:sumofforces}}
	\caption{Coefficients $f_R$ in \eqref{eq:repulsionforcemodel} and $f_A$ in \eqref{eq:attractionforcemodel} of repulsion force \eqref{eq:repulsionforce} and attraction force \eqref{eq:attractionforce}, respectively, as well as total force coefficients  along the lines of largest and smallest stress for $\chi=0.2$ (i.e.\ $f_A+f_R$ and $0.2f_A+f_R$, respectively) for parameter values in  \eqref{eq:parametervaluesRepulsionAttraction}}
\end{figure}

In Figure \subref*{fig:forces} the coefficient functions  \eqref{eq:repulsionforcemodel} and \eqref{eq:attractionforcemodel} for the repulsion and attraction forces \eqref{eq:repulsionforce} and \eqref{eq:attractionforce} in the K\"{u}cken-Champod model \eqref{eq:particlemodel}  are plotted for the parameters in \eqref{eq:parametervaluesRepulsionAttraction}. 

The sums of the coefficients of the forces $f_R+f_A$ and $f_R+\chi f_A$ for $\chi=0.2$ are illustrated in Figure \subref*{fig:sumofforces}. Note that $f_R+f_A$ and $f_R+\chi f_A$ are the force coefficients along $l$ and $s$, respectively. For the choice of parameters in \eqref{eq:parametervaluesRepulsionAttraction}   repulsion dominates for short distances along the lines of largest stress to prevent the collision of particles and the force is long-range attractive along the lines of largest stress  leading to accumulations of the particles. The absolute value of the attractive force decreases with the distance between particles. Along the lines of smallest stress the particles are always repulsive for $\chi=0.2$, independent of the distance, though the  repulsion force gets weaker for longer distances. 

\subsection{Bio-inspired model}
We propose a system of ordinary differential equations of the form \eqref{eq:particlemodel} where the forces are of the form \eqref{eq:totalforcenew}. Note that plugging  the repulsion and attraction forces \eqref{eq:repulsionforce} and \eqref{eq:attractionforce} as well as the definition \eqref{eq:tensorfield} of the tensor field $T$  into the force term \eqref{eq:totalforce} results in forces of the form \eqref{eq:totalforcenew}. Hence, to generalise the Kücken-Champod model we require for the coefficient functions $f_s$ and $f_l$:
\begin{align*}
f_s\approx \chi f_A+f_R, \qquad f_l\approx f_A+f_R.
\end{align*}

\begin{figure}[htbp]
	\centering
	\includegraphics[width=0.45\textwidth]{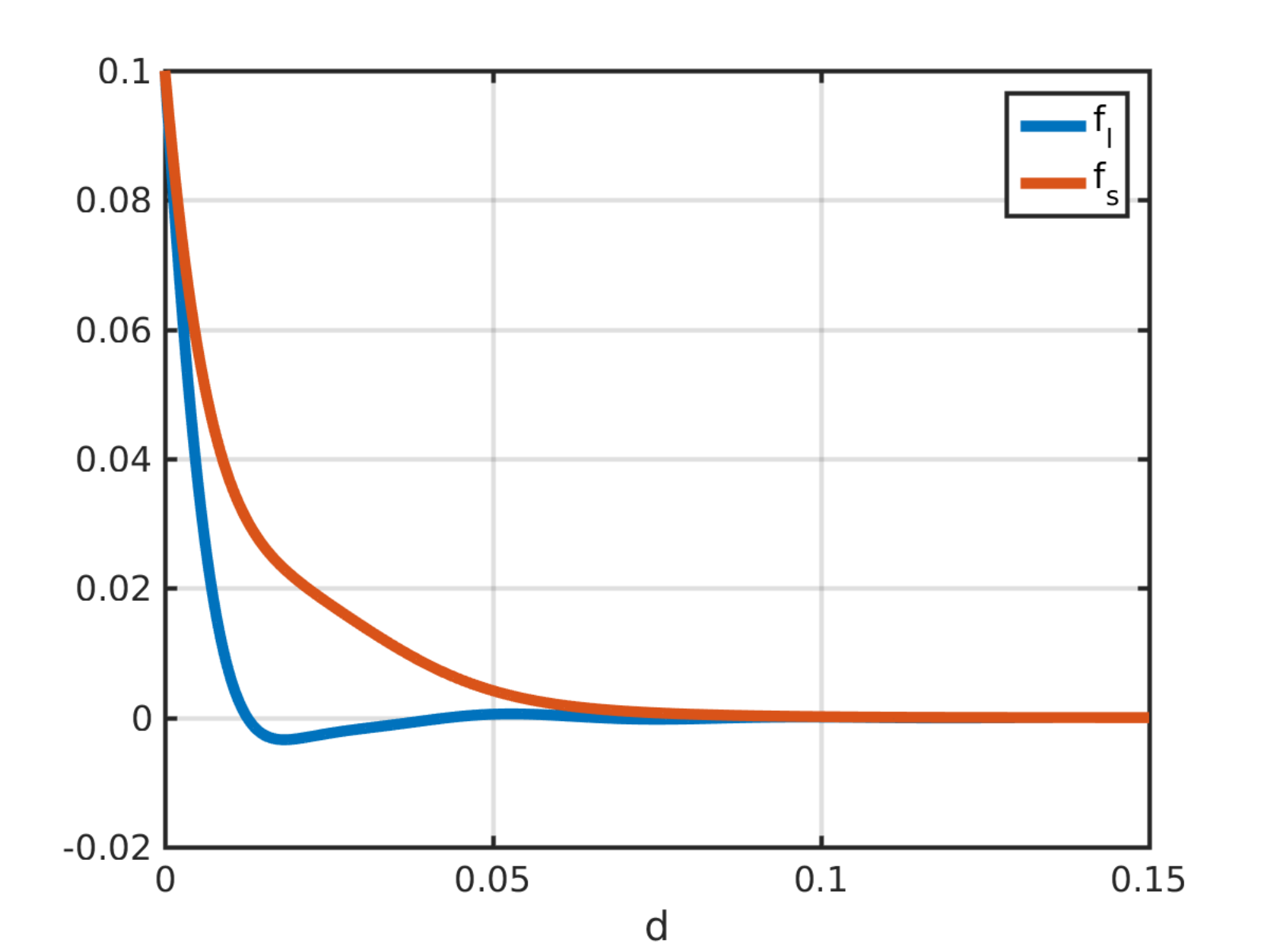}
	\caption{Coefficients $f_s$ and $f_l$ in \eqref{eq:forcecoeffharmonic}  for parameter values in  \eqref{eq:parameterharmonic}}\label{fig:forcesharmonic}
\end{figure}

We model the force coefficients $f_s$ and $f_l$ in \eqref{eq:totalforcenew} as solutions to a damped harmonic oscillator. Like for the coefficient functions \eqref{eq:repulsionforcemodel}, \eqref{eq:attractionforcemodel} in the Kücken-Champod model  we consider  exponentially decaying forces  describing that short-range interactions between the particles  are much stronger than long-range interactions. Besides, the repulsion and attraction forces suggested in the Kücken-Champod model dominate on different regimes. For a more unified modeling approach one may regard this interplay of repulsion and attraction forces as oscillations. This motivates to model the force coefficients $f_s$ and $f_l$ in \eqref{eq:totalforcenew} as solutions to a damped harmonic oscillator which is also a common modeling approach in cell biology. Hence, we consider the following ansatz functions for the force coefficients $f_s$ and $f_l$:
\begin{align}\label{eq:forcecoeffharmonic}
\begin{split}
f_s(d)&=  c\exp(e_{s_1}|d|)+c_{s}\sin\bl\frac{\pi|d|}{a_s}\br\exp(e_{s_2}|d|) ,\\ f_l(d)&=c\cos\bl\frac{\pi|d|}{a_{l}}\br\exp(e_{l_1}|d|)+c_{l}\sin\bl\frac{\pi|d|}{a_{l}}\br\exp(e_{l_2}|d|)
\end{split}
\end{align}
for real parameters $c,c_s,c_l,e_{s_1},e_{s_2},e_{l_1},e_{l_2},a_s,a_l$.
In the sequel, we will see that for the parameter choice
\begin{align}\label{eq:parameterharmonic}
\begin{split}
c=0.1,\quad c_{s}=-0.05,\quad e_{s_1}=-65.0,\quad e_{s_2}=-100.0,\quad a_s=0.03 \\
c_{l}=0.005,\quad e_{l_1}=-160.0,\quad e_{l_2}=-40.0,\quad a_{l}=0.022
\end{split}
\end{align}
it is possible to obtain stationary fingerprint patterns and that rescaling of the coefficient functions $f_s$ and $f_l$ leads to stationary patterns with scaled line distances. The force coefficients $f_s$ and $f_l$ for the parameters in \eqref{eq:parameterharmonic} are shown in Figure \ref{fig:forcesharmonic}. In comparison with the force coefficients $F_A+f_R$ and $0.2 f_A+f_R$ along $l$ and $s$, respectively,  the force $f_s$ along $s$ is also purely repulsive, while the force $f_l$ is less attractive which is necessary for obtaining stationary patterns as discussed in Section \ref{sec:numerics}.

\subsection{General setting}\label{sec:generalsetting}
In the sequel, we consider the particle model \eqref{eq:particlemodel} with force terms of the form $F(x_j-x_k,T(x_j))$, such as \eqref{eq:totalforcenew} and \eqref{eq:totalforce}. As in \cite{patternformationanisotropicmodel} we consider the domain $\Omega=\mathbb{T}^2$ where $\mathbb{T}^2$ is  the $2$-dimensional unit torus that can be identified with the unit square $[0,1)\times [0,1)\subset \R^2$ with periodic boundary conditions. Note that the particles on the domain $\Omega$ are separated by a distance of at most 0.5. Motivated by this we require for $j\in\{1,2\}$ and all $x\in \Omega$
%Because of the periodic boundary conditions we assume additionally that particles cannot interact with copies of themselves as well as with more than one copy of other particles. To guarantee this we assume that interactions can only occur if the particles are separated by a distance  of less than 0.5 in each spatial direction, i.e.\ for $j\in\{1,2\}$ and all $x\in \Omega$ we require that 
\begin{align}\label{eq:minimumimagecriterion}
F(x-x',T(x))\cdot e_j=0\quad \text{for}\quad |x-x'|\geq 0.5
\end{align}
where $e_j$ denotes the standard basis for the Euclidean plane. The forces satisfy this assumption if a spherical cutoff radius of length $0.5$ is introduced for the forces in \eqref{eq:totalforcenew} or \eqref{eq:totalforce}, respectively. This assumptions guarantees that the size of the domain is large enough compared to the range of the total force. In particular, non-physical artifacts due to periodic boundary conditions are prevented. A cutoff radius is also very useful to make numerical simulations more efficient.

\section{Mathematical analysis of steady states}\label{sec:existencesteadystates}

To use the particle system \eqref{eq:particlemodel} for the simulation of fingerprints it is of great interest to have a better understanding about the form of the steady states. The steady states are formed by a number of lines which are referred to as ridges. As discussed in  Section \ref{sec:descriptionmodel} we consider purely repulsive forces along $s$.
In this section, we study the existence of  steady states for the particle model \eqref{eq:particlemodel} for spatially homogeneous and certain inhomogeneous tensor fields $T$ analytically. The stability of these line patterns is further investigated in \cite{stabilityanalysisanisotropicmodel}. In particular, the authors show that line patterns for purely repulsive forces along $s$ can only be stable if the patterns are aligned in direction of the vector field  $s$.

\subsection{Spatially homogeneous tensor field}
For spatially homogeneous tensor fields $T$ it is sufficient to restrict ourselves to the tensor field given  by $s=(0,1)$ and $l=(1,0)$ since stationary solutions to the Kücken-Champod model for any other tensor field can be obtained by coordinate transform \cite{patternformationanisotropicmodel}. Further note that steady states are translation invariant, i.e.\ if $x_1,\ldots,x_N$ is a steady state, so is $x_1+z,\ldots,x_N+z$ for any $z\in\R^2$. Hence it is sufficient to consider one specific constellation of particles for analysing the steady states of \eqref{eq:particlemodel}. Because of the stability analysis in \cite{stabilityanalysisanisotropicmodel} we restrict ourselves to line patterns along $s=(0,1)$, i.e.\ we consider patterns of vertical lines.

%Model \eqref{eq:particlemodel} for an underlying spatially homogeneous tensor field was analyzed mathematically and numerically in \cite{patternformationanisotropicmodel}. 

\begin{proposition}
	Given  $d\in(0,1]$ such that $n:=\frac{1}{d}\in\N$ and let $N\in\N$ be given such that $\frac{N}{n}\in \N$. Then $n$ parallel equidistant vertical lines of distance $d$ of uniformly distributed particles along each line are a steady state to the particle model \eqref{eq:particlemodel} for forces of the form \eqref{eq:totalforcenew} or \eqref{eq:totalforce} where the  repulsion and attraction forces are of the form \eqref{eq:repulsionforce} and \eqref{eq:attractionforce}, respectively.
\end{proposition}
Note that the choice of the distance $d$ of the parallel vertical lines is consistent with the periodic boundary conditions.
\begin{proof}
	Because of the translational invariance of steady states it is sufficient to consider any $n$ equidistant parallel vertical lines of $\frac{N}{n}$ particles distributed uniformly along each line. Without loss of generality we assume that the positions of the particles are given by
	\begin{align*}
	\bar{x}_j=\bl \frac{\bl j-j\bmod \frac{N}{n}\br\frac{n}{N}}{n},\frac{j\bmod \frac{N}{n}}{\frac{N}{n}}\br=\frac{1}{N}\bl  j-j\bmod \frac{N}{n},n\bl j\bmod \frac{N}{n}\br\br\in\R^2.
	\end{align*}
	
	Because of the periodic boundary conditions of the domain as well as the fact that the particles are uniformly distributed along parallel lines, it is sufficient to require that 
	\begin{align}\label{eq:straightlinecondition}
	\sum_{k=1}^{N-1} F(\bar{x}_N-\bar{x}_k,T(\bar{x}_N))=0
	\end{align}
	for steady states of the particle model \eqref{eq:particlemodel}. Note that for forces of the form  \eqref{eq:totalforcenew} or \eqref{eq:totalforce} where the  repulsion and attraction forces are of the form \eqref{eq:repulsionforce} and \eqref{eq:attractionforce}, respectively, we have 
	\begin{align}\label{eq:newtonthree}
	F(d,T(\bar{x}_N))=-F(-d,T(\bar{x}_N))\quad\text{for all}\quad d\in\R^2.
	\end{align} 
	As a first step we show that 
	\begin{align}\label{eq:forcebalancehelp}
	\sum_{k=1}^{\frac{N}{n}-1} F(\bar{x}_N-\bar{x}_k,T(\bar{x}_N))=0.
	\end{align}
	Note that $\bar{x}_k\in\{0\}\times [0,1]$ for $k=1,\ldots,\frac{N}{n}$ and $\bar{x}_N=(0,0)$ by the periodic boundary conditions, i.e.\ we consider all the particles of the vertical line with $x_1$-coordinate $x_1=0$. If $\frac{N}{n}$ is odd, then \eqref{eq:forcebalancehelp} is satisfied by the balance of forces \eqref{eq:newtonthree}. For even $\frac{N}{n}$ we have 
	\begin{align*}
	F(\bar{x}_N-\bar{x}_k,T(\bar{x}_N))=-F(\bar{x}_N-\bar{x}_{\frac{N}{n}-k},T(\bar{x}_N))
	\end{align*}
	for $k=1,\ldots,\frac{N}{2n}-1$. Besides, $$F(\bar{x}_N-\bar{x}_{\frac{N}{2n}},T(\bar{x}_N))=0$$ since  $|\bar{x}_N-\bar{x}_{\frac{N}{2n}}|=0.5$ and the assumption of the finite range of the forces in \eqref{eq:minimumimagecriterion}, implying that \eqref{eq:forcebalancehelp} is satisfied. If there is an odd number $n$ of parallel equidistant  vertical lines, then the condition for steady states \eqref{eq:straightlinecondition} is satisfied by \eqref{eq:newtonthree}. For $n$ even, the forces due to particles on the vertical lines at $x_1=kd$ balances the interaction forces due to particles on the  vertical lines at $x_1=(n-k)d$ for $k=1,\ldots,\frac{n}{2}-1$ by \eqref{eq:newtonthree}, so it suffices to consider the particles on the vertical line at $x_1=\frac{n}{2}d$, i.e.\ the particles at positions $\bar{x}_{k}$ for $k=\frac{N}{2},\ldots,\frac{N}{2}+\frac{N}{n}-1$. Note that $$\sum_{k=\frac{N}{2}}^{\frac{N}{2}+\frac{N}{n}-1} F(\bar{x}_N-\bar{x}_{k},T(\bar{x}_N))=0$$ since  $|\bar{x}_N-\bar{x}_{k}|\geq 0.5$ for $k=\frac{N}{2},\ldots,\frac{N}{2}+\frac{N}{n}-1$ and the assumption of the finite range of the forces in \eqref{eq:minimumimagecriterion}. This implies that the condition for steady states \eqref{eq:straightlinecondition} is satisfied. Hence, $\bar{x}_1,\ldots,\bar{x}_N$ form a steady state of the microscopic model \eqref{eq:particlemodel}.
\end{proof}

\subsection{Non-constant tensor fields}
Many non-constant tensor fields can locally be regarded as spatially homogeneous tensor fields. Note that by the assumptions in Section \ref{sec:generalsetting} we consider forces of finite range. In particular, we have local forces for the forces \eqref{eq:totalforcenew} with coefficients \eqref{eq:forcecoeffharmonic} and parameters \eqref{eq:parameterharmonic} as well as for forces of the form \eqref{eq:totalforce} with force coefficients \eqref{eq:repulsionforcemodel}, \eqref{eq:attractionforcemodel} and parameters \eqref{eq:parametervaluesRepulsionAttraction}. Hence, roughly parallel lines form a steady states of the particle model \eqref{eq:particlemodel} for locally spatially homogeneous tensor fields. However, it is not clear whether these possible steady states are also stable.

\section{Simulation of fingerprint patterns} \label{sec:fingerprintsimulations}

In this section we investigate how to simulate fingerprint patterns by extending the theoretical and the numerical results in \cite{patternformationanisotropicmodel}. In particular, we consider more realistic tensor fields for the formation of fingerprint patterns and study  the dependence of the parameter values in the K\"{u}cken-Champod model on the resulting fingerprints.

\subsection{Local fields in a fingerprint image} \label{sec:fingerprinttensor}
\begin{figure}[htbp]
	\subfloat[Delta]{\includegraphics[width=0.45\textwidth]{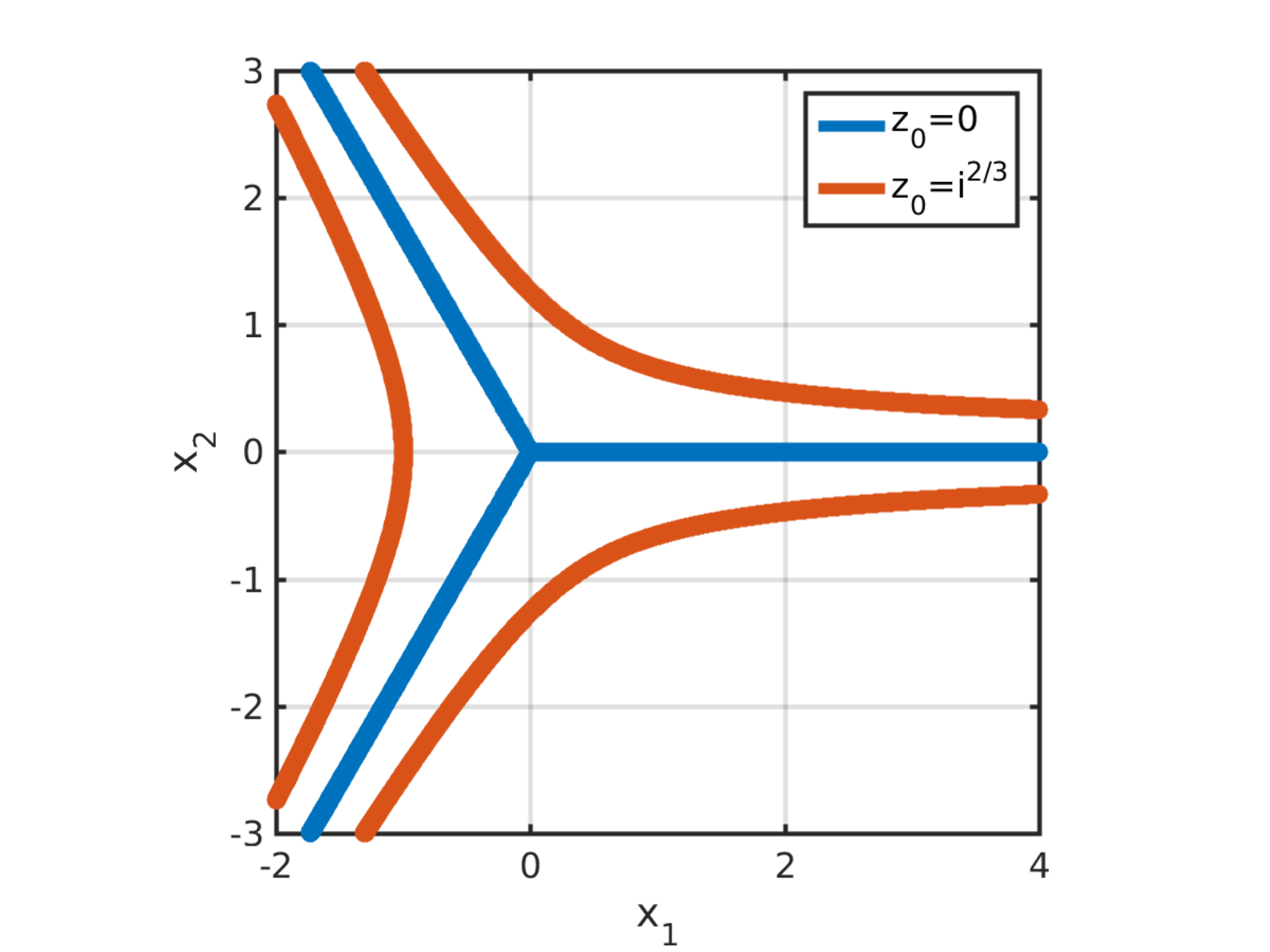}\label{fig:deltaode}}
	\subfloat[Core]{\includegraphics[width=0.45\textwidth]{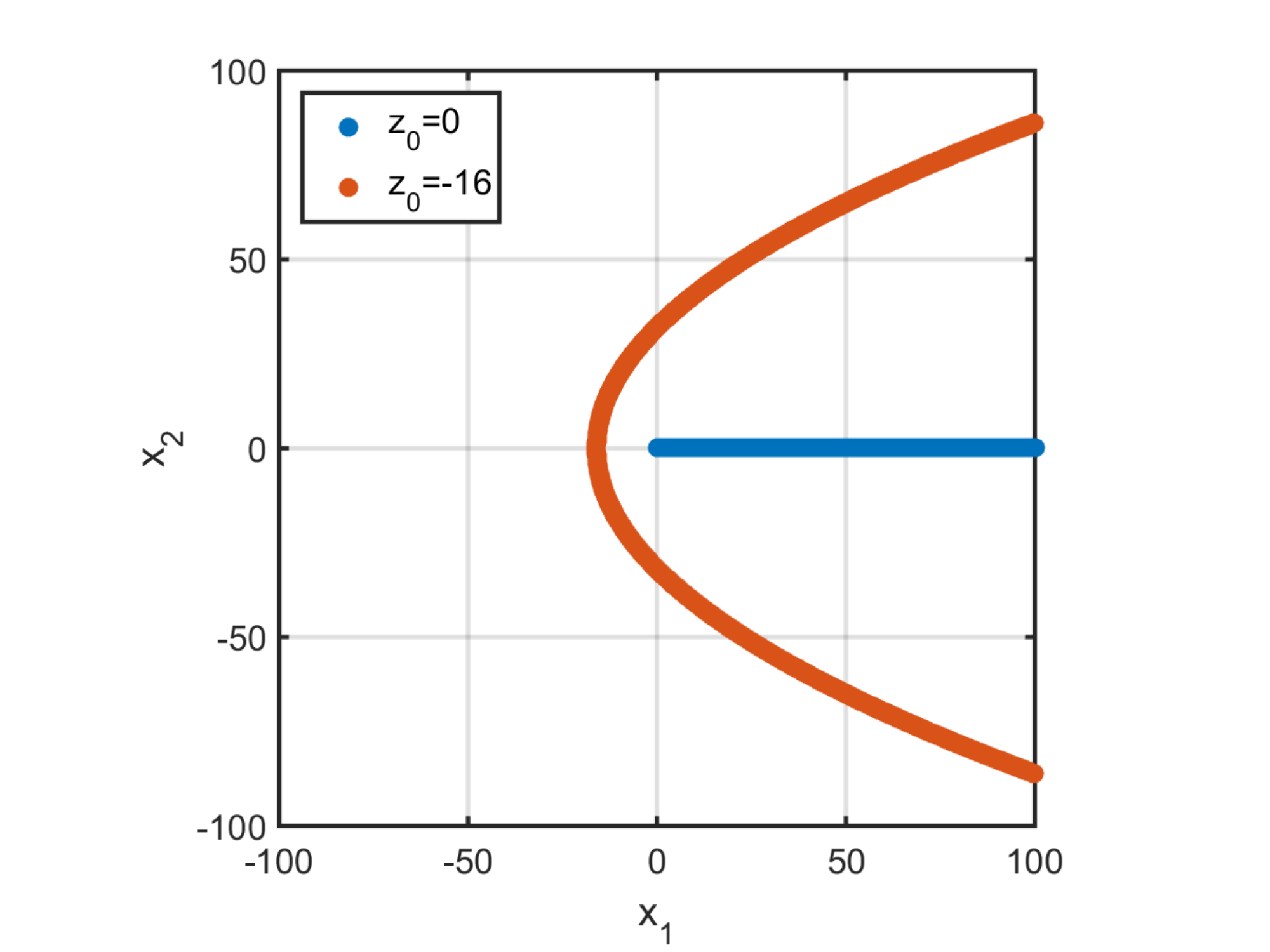}\label{fig:coreode}}
	\caption{Solution curves \eqref{eq:solutiondelta} and \eqref{eq:solutioncore} to the initial value problems \eqref{eq:odedelta} and \eqref{eq:odecore}, respectively, generating fields of quadratic differentials for a delta and a core}
\end{figure}

In order to use the particle model \eqref{eq:particlemodel} to simulate fingerprint patterns a realistic tensor field is needed. It is well known that fingerprints are composed of two key directional features known as cores and deltas. Hence, we consider the construction of the tensor fields for these two features first. In \cite{tensorfieldfingerprints} Huckemann et al. propose to use quadratic differentials for generating the global fields in a fingerprint image. The local field is then generated by the singular points of the field. A core is the endpoint of a single line (cp. Figure \subref*{fig:coreode}) and a delta occurs at the junction of three lines (cp. Figure \subref*{fig:deltaode}).

\begin{figure}[htbp]
	\subfloat[$s$ for delta]{\includegraphics[width=0.45\textwidth]{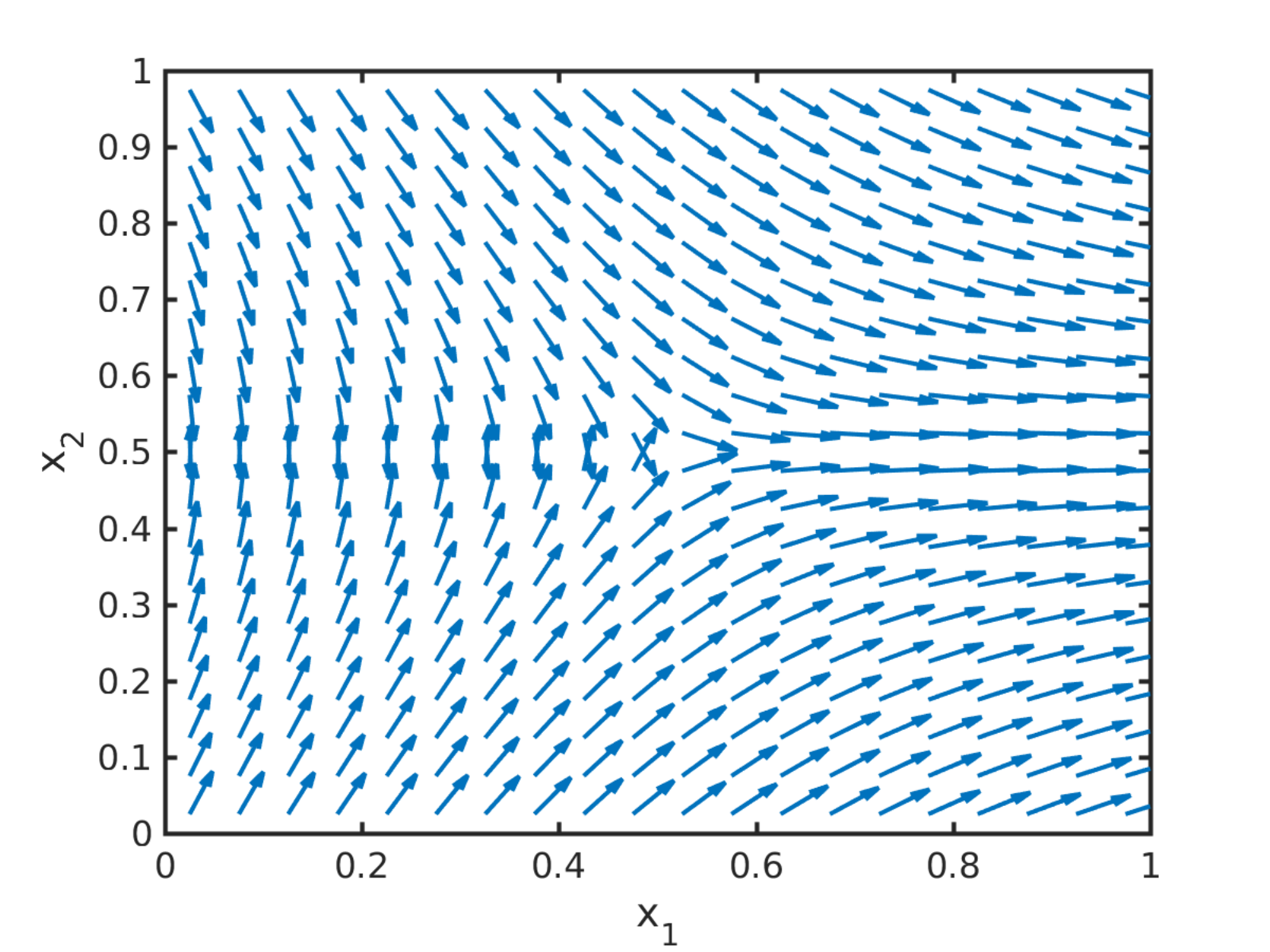}\label{fig:deltatensor}}
	\subfloat[$s$ for core]{\includegraphics[width=0.45\textwidth]{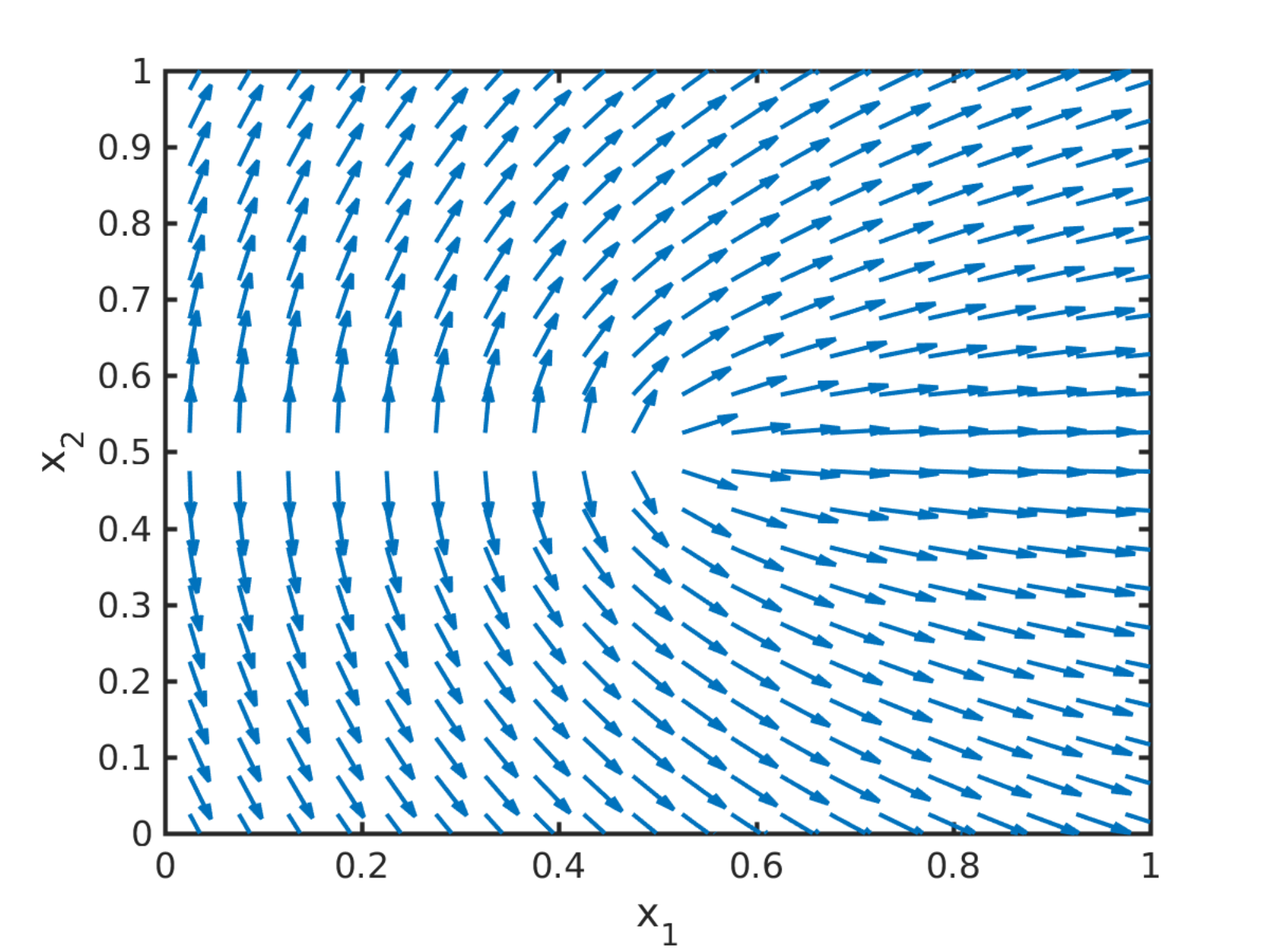}}
	\caption{Lines of smallest stress $s=s(x)$ of tensor fields $T$ for a delta and a core}\label{fig:deltacoretensor}
\end{figure}

For simplicity we consider the origin ($\zeta=0$) as the only singular point, but the idea can be extended to arbitrary singular points $\zeta\in\C$. As outlined in \cite{tensorfieldfingerprints} one can model the field near the origin $\zeta=0$ by considering the initial value problem
\begin{align}\label{eq:odedelta}
z(r)\dot{z}(r)^2=\phi(r),\quad z(r_0)=z_0
\end{align}
for a smooth, positive function $\phi=\phi(r), ~r\in\R,$ and $z_0\in\C$. For $\phi=\frac{2}{3}$ the solution to the differential equation is given by 
\begin{align}\label{eq:solutiondelta}
z(r)=\bl r+z_0^{3/2}\br^{2/3}.
\end{align}
Note that the shape of the solution curves does not change for reparametrizations, provided $\phi>0$. By varying $z_0\in\C$ the associated solution curves form a delta at the origin ($\zeta=0$) as illustrated in Figure \subref*{fig:deltaode}. Hence, we require $$z\di z^2>0$$ for a delta at the origin. Note that $z=|z|\exp( i \text{arg}(z))$ where $\text{arg}(z)$ denotes the principal argument of the complex number $z\in\C$. Further note that $\di z$ can be regarded as the direction of the smallest stress at $z\in\C$ if $\R^2$ is identified with $\C$. As outlined in Section \ref{sec:descriptionmodel} the direction of smallest stress is denoted by the unit vector $s=s(z)$ for $z\in\R^2$ implying  $\di z=\pm \exp( -i \text{arg}(z) /2)$. Thus, the lines of smallest stress on a domain $\Omega\subset\C$ can be obtained  by evaluating  $\di z$ for all $z\in\Omega$. Note that $\exp( -i \text{arg}(z) /2)$ and $-\exp( -i \text{arg}(z) /2)$ result in the same lines of the stress field. 
Similarly, the  initial value problem
\begin{align}\label{eq:odecore}
\frac{1}{z(r)}\dot{z}(r)^2=\phi(r),\quad z(r_0)=z_0
\end{align}
generates a field with a core at the origin. Up to reparameterization the solution is given by 
\begin{align}\label{eq:solutioncore}
z(r)=\bl r+z_0^{1/2}\br^2
\end{align}
and the solution curves are illustrated for different initial conditions $z_0\in\C$ in Figure \subref*{fig:coreode}. The condition $$\frac{\di z^2}{z}>0$$ has to be satisfied for a core at the origin, implying $\di z=\pm \exp( i \text{arg}(z) /2)$, and as before $\pm \exp( i \text{arg}(z) /2)$ result in the same lines. Further note that a delta or a core at any $\zeta\in\C$ can be obtained by linear transformation. In Figure \ref{fig:deltacoretensor} the tensor field for a delta and a core at the singular point $(0.5,0.5)$ are plotted on the unit square $[0,1]^2$.

\subsection{Numerical methods}\label{sec:numerics}
In this section, we describe the general setting for investigating the long-time behavior of solutions to the particle model \eqref{eq:particlemodel}, motivated by \cite{patternformationanisotropicmodel}. 

We consider the particle model \eqref{eq:particlemodel}  where the forces are of the form \eqref{eq:totalforcenew} or \eqref{eq:totalforce} and  investigate the patterns of the corresponding stationary solutions. As in \cite{patternformationanisotropicmodel} and outlined in Section \ref{sec:generalsetting} we consider the domain $\Omega=\mathbb{T}^2$, i.e.\ the unit square $[0,1)\times [0,1)\subset \R^2$ with periodic boundary conditions, and we consider a cutoff of the forces as in \eqref{eq:minimumimagecriterion} to make the simulations more efficient. 

To solve the $N$ particle ODE system \eqref{eq:particlemodel} we apply either the simple explicit Euler scheme or higher order methods such as the Runge-Kutta-Dormand-Prince method, all resulting in very similar simulation results. For the numerical simulations we consider $\Delta t=0.2$ for the size of the time step.

\subsection{Numerical study of the Kücken-Champod model}\label{sec:numerickcmodel}
\begin{figure}[htbp]
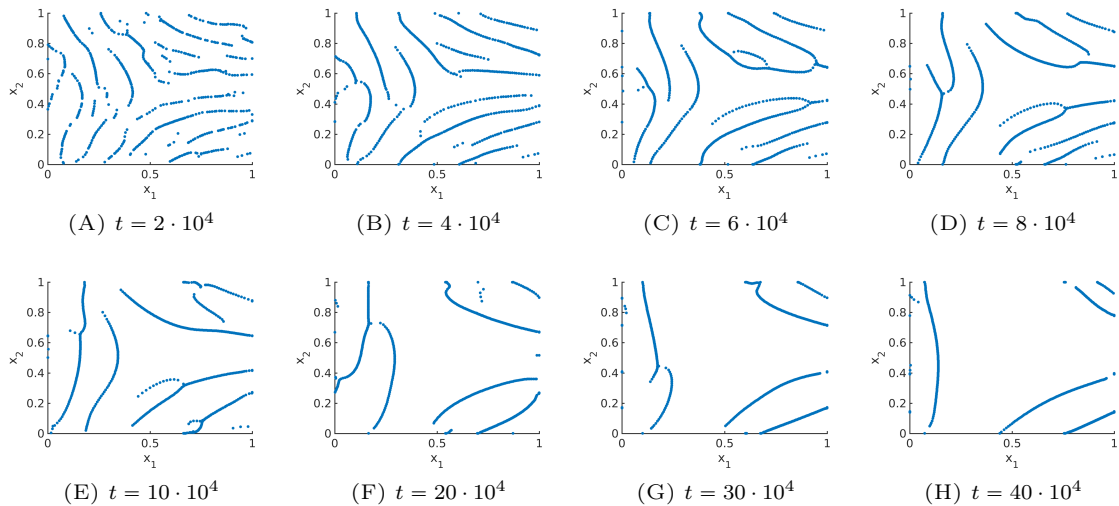

	\foreach \x in {2,4,...,10}{%
		\subfloat[$t=\x\cdot 10^4$]{\includegraphics[width=0.24\textwidth]{DeltaN600FirstAttemptTime\x.pdf}}\hfill
	}%
	\foreach \x in {20,30,...,40}{%
		\subfloat[$t=\x\cdot 10^4$]{\includegraphics[width=0.24\textwidth]{DeltaN600FirstAttemptTime\x.pdf}}\hfill
	}%
	\caption{Numerical solution to the K\"{u}cken-Champod model \eqref{eq:particlemodel} for $N=600$ and $\chi=0.2$ at different times $t$ where the stress field represents a delta and the cutoff radius is $0.5$}\label{fig:deltafirsttry}
\end{figure}

%\subsection{Complex stationary patterns for non-homogeneous tensor fields} \label{sec:extensionfingerprints}
Using the tensor fields introduced in Section \ref{sec:fingerprinttensor} we consider the interaction model \eqref{eq:particlemodel} with forces of the form \eqref{eq:totalforce} to simulate fingerprint patterns. Here, the repulsion and attraction forces are of the forms \eqref{eq:repulsionforce} and \eqref{eq:attractionforce} with force coefficients \eqref{eq:repulsionforcemodel} and \eqref{eq:attractionforcemodel}, respectively, and we consider the parameters in \eqref{eq:parametervaluesRepulsionAttraction} to make the simulations as close as possible to the model suggested by K\"ucken and Champod in \cite{Merkel}.  It is well known that fingerprints develop during pregnancy and stay  the same afterwards provided no fingerprint alterations occur. In order to simulate biologically meaningful fingerprints we  aim to model fingerprint patterns as stationary solution to  the particle model \eqref{eq:particlemodel}. Based on the analysis of steady states in Section \ref{sec:existencesteadystates} it is possible to obtain stationary patterns consisting of multiple roughly parallel ridges along the lines of smallest stress. However, the force coefficients need to be chosen appropriately so that the resulting patterns are also stable. For the simulations in Figure \ref{fig:deltafirsttry} we consider  the tensor field for the delta constructed in Section \ref{sec:fingerprinttensor} and depicted in Figure \subref*{fig:deltatensor}. 
One can clearly see in Figure \ref{fig:deltafirsttry} that the particles are aligned along the lines of smallest stress $s=s(x)$ initially, but the patterns dissolve over time and  the simulation results have little similarity with fingerprint patterns over large time intervals. Besides, the patterns are clearly no stable steady states  in Figure \ref{fig:deltafirsttry}. Hence, the question arises  why the patterns simplify so much over time for non-homogeneous tensor fields in contrast to the stationary patterns arising for spatially homogeneous tensor fields, cf. \cite{patternformationanisotropicmodel}, and how this can be prohibited.  

\begin{figure}[htbp]
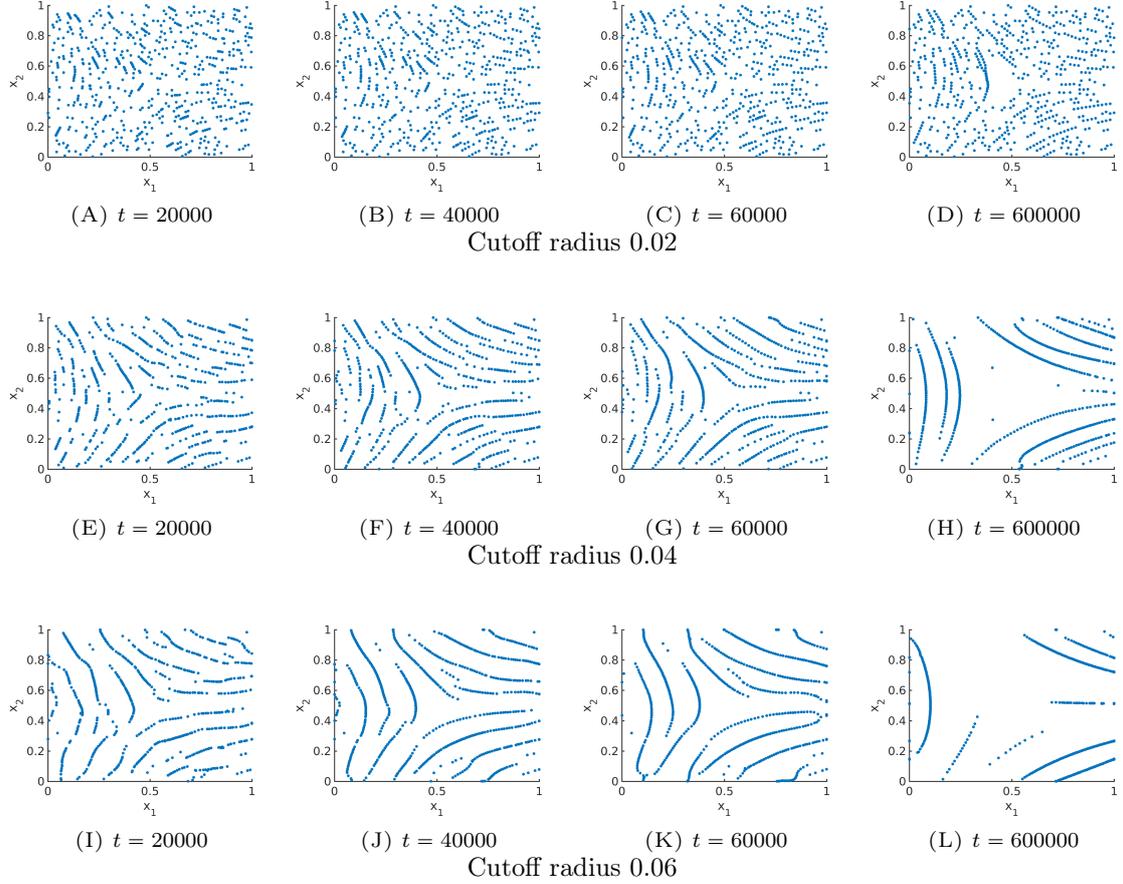

	\centering
	\begin{multicols}{4}
		\begin{minipage}{\textwidth}
			\centering
			\foreach \x in {2,4,...,6,60}{%
				\subfloat[$t=\x0000$]{\includegraphics[width=0.24\textwidth]{DeltaN600Cutoff002FirstAttemptCutOffTime\x.pdf}}\hfill
			}
			\vspace{3mm}
			Cutoff radius 0.02
		\end{minipage}
		\begin{minipage}{\textwidth}
			\centering
			\foreach \x in {2,4,...,6,60}{%
				\subfloat[$t=\x0000$]{\includegraphics[width=0.24\textwidth]{DeltaN600Cutoff004FirstAttemptCutOffTime\x.pdf}}\hfill
			}
			\vspace{3mm}
			Cutoff radius 0.04
		\end{minipage}
		\begin{minipage}{\textwidth}
			\centering
			\foreach \x in {2,4,...,6,60}{%
				\subfloat[$t=\x0000$]{\includegraphics[width=0.24\textwidth]{DeltaN600Cutoff006FirstAttemptCutOffTime\x.pdf}}\hfill
			}
			\vspace{3mm}
			Cutoff radius 0.06
		\end{minipage}
	\end{multicols}
	\caption{Numerical solution to the K\"{u}cken-Champod model \eqref{eq:particlemodel} for different cutoff radii for $N=600$ and $\chi=0.2$  at different times $t$ where the stress field represents a delta}\label{fig:deltafirsttrycutoff}
\end{figure}

%\subsubsection{Cutoff radius for efficient simulations}
To study the long-time behaviour of the numerical solution, it is desirable to have efficient numerical simulations and of course efficient simulations are also necessary to 
%apply the methods described in this section for fingerprint simulations 
to simulate fingerprints based on cell interactions in practice. 
%By the minimum image criterion a cutoff radius has to be introduced for the repulsion and attraction forces 
In Section \ref{sec:generalsetting} we introduced a cutoff radius for the forces, given by \eqref{eq:minimumimagecriterion}, 
in order to deal with the periodic boundary conditions. Since the forces in the K\"{u}cken-Champod model \eqref{eq:particlemodel} decrease exponentially, the interaction force between two particles is very small if their distance is sufficiently large.  This is also illustrated in Figure \subref*{fig:forces} for the parameters in \eqref{eq:parametervaluesRepulsionAttraction}. Hence, defining the cutoff radius as $0.1$ changes the values of the forces only slightly, but it allows us to compute the numerical solution to the  K\"{u}cken-Champod model \eqref{eq:particlemodel} by using cell lists \cite{celllists}. The idea of cell lists is to subdivide the simulation domain into cells with edge lengths greater than or equal to the cutoff radius of the interaction forces. All particles are sorted into these cells and only particles in the same or neighbouring cells have to be considered for interactions. This results in significantly faster simulations since we only have to consider those particle pairs with relevant sizes of the interaction forces. Note that the cutoff radius has an impact on the number of lines that occur in the solution as shown in Figure \ref{fig:deltafirsttrycutoff} in comparison to a cutoff radius of $0.5$ in Figure \ref{fig:deltafirsttry}. In particular the cutoff radius should not be chosen to small because this prevents the accumulation of particles.

%\subsubsection{Motivation for adaptation of the K\"{u}cken-Champod model}
The simulation results for the K\"{u}cken-Champod model \eqref{eq:particlemodel} in Figures \ref{fig:deltafirsttry} and \ref{fig:deltafirsttrycutoff} illustrate that the particles align in roughly parallel lines along the lines of smallest stress initially, but the number of roughly parallel lines decreases as time goes on. In particular, the   complex  patterns that occur initially  are not stationary. We can expect a similar behavior (i.e.\ initial alignment along the lines of smallest stress of the stress tensor field and subsequent accumulation) of the numerical solution  if the parameters in the coefficient functions of the repulsion and attraction force in  \eqref{eq:repulsionforcemodel} and \eqref{eq:attractionforcemodel} are slightly changed provided  they are repulsive along the lines of smallest stress, as well as short-range repulsive and long-range attractive along the lines of largest stress. Denoting the directions of smallest and largest stress by $s$ and $l$, respectively, the transition of the initial pattern of multiple lines to fewer and fewer lines along $s$ suggests that the attraction forces are very strong resulting in an accumulation of the particles. Note that  this transition is also observed for the long-time behavior of the numerical solution to the K\"{u}cken-Champod model \eqref{eq:particlemodel} for spatially homogeneous tensor fields in \cite{patternformationanisotropicmodel} where lines merge over time until finally a steady state of equidistant parallel lines is reached.

\begin{figure}[htbp]
	\centering
	\subfloat[$s$]{        \includegraphics[width=0.24\textwidth]{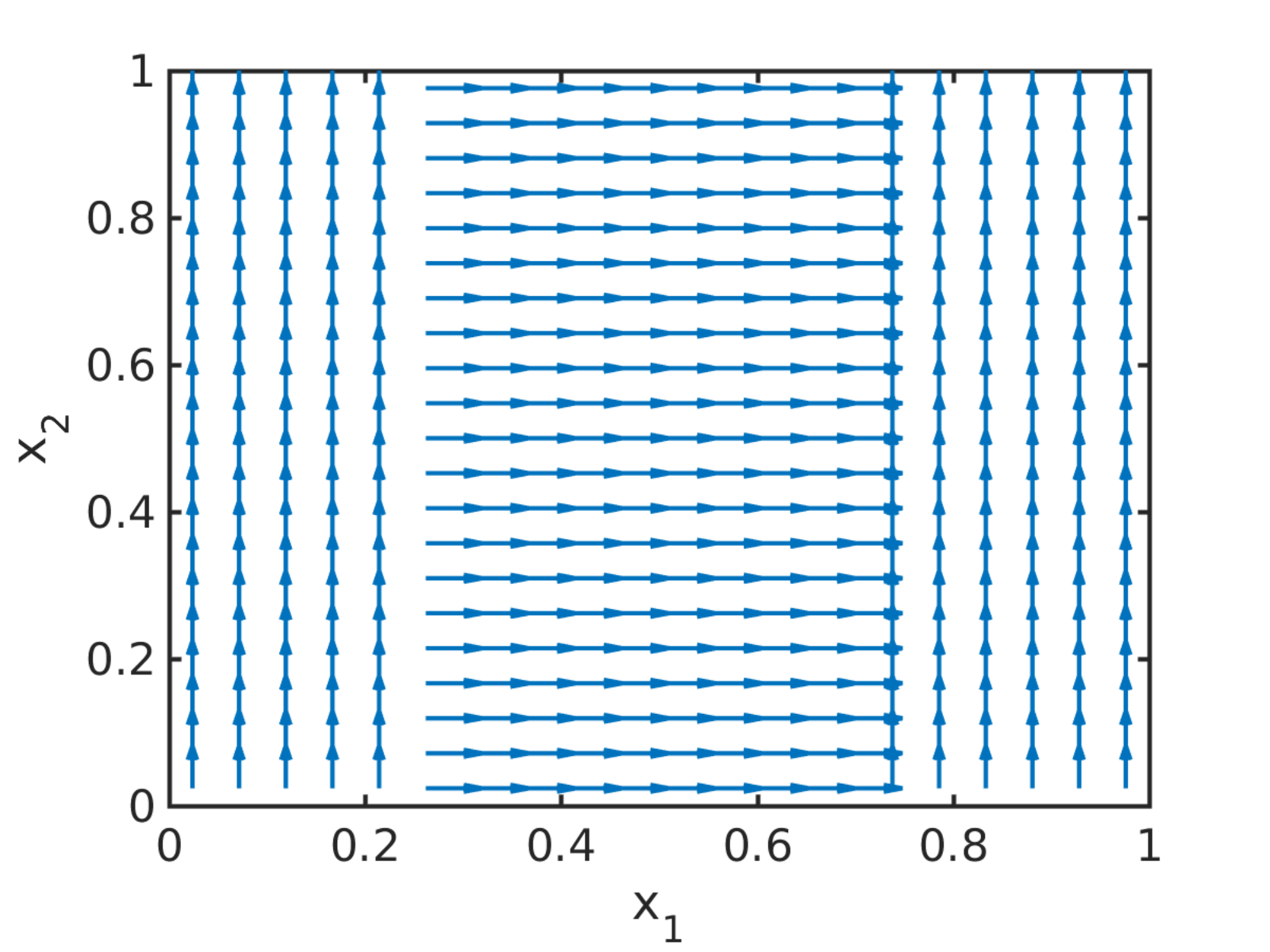}}
	\subfloat[$t=40000$]{       \includegraphics[width=0.24\textwidth]{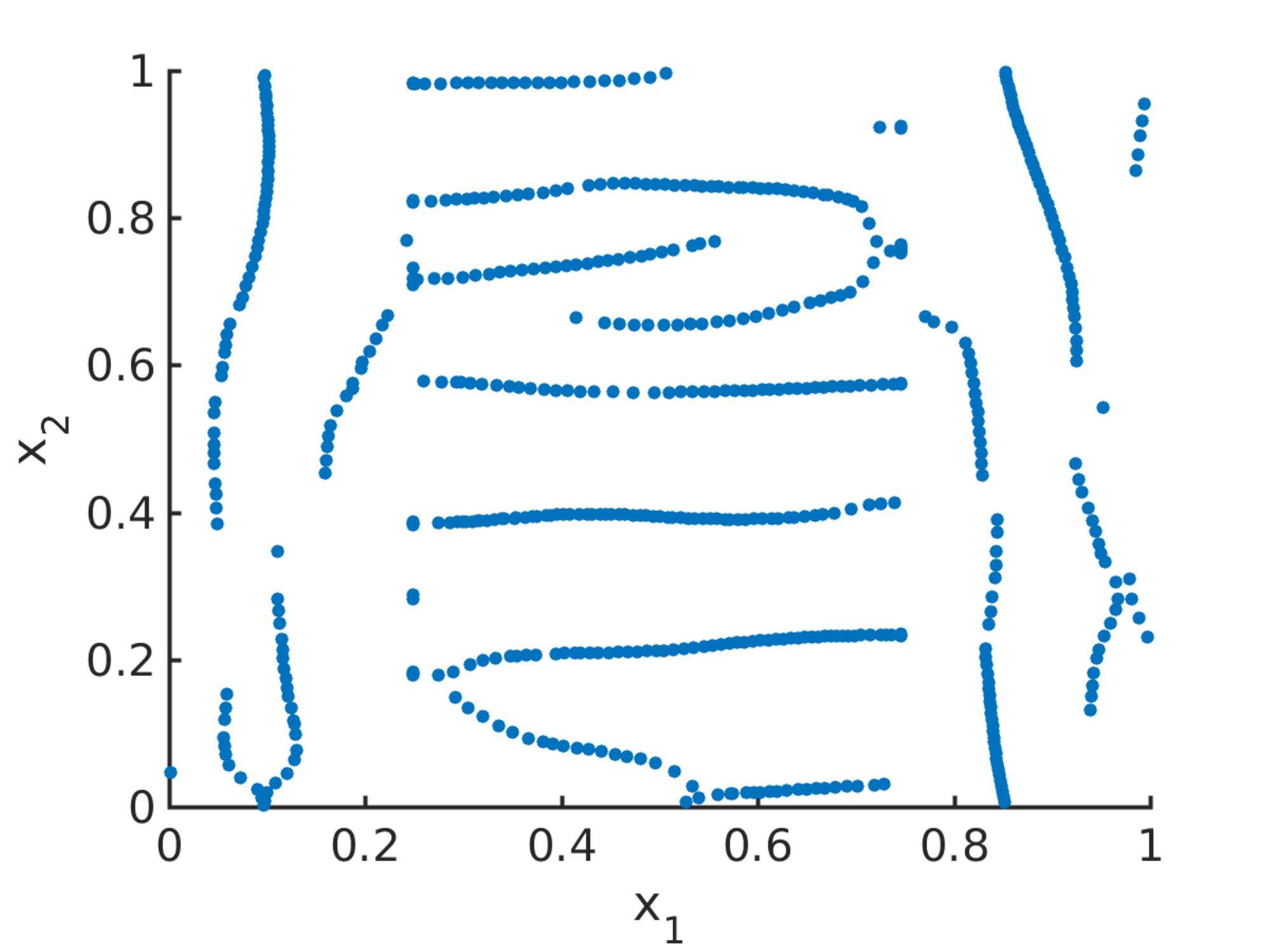}}
	\subfloat[$t=200000$]{        \includegraphics[width=0.24\textwidth]{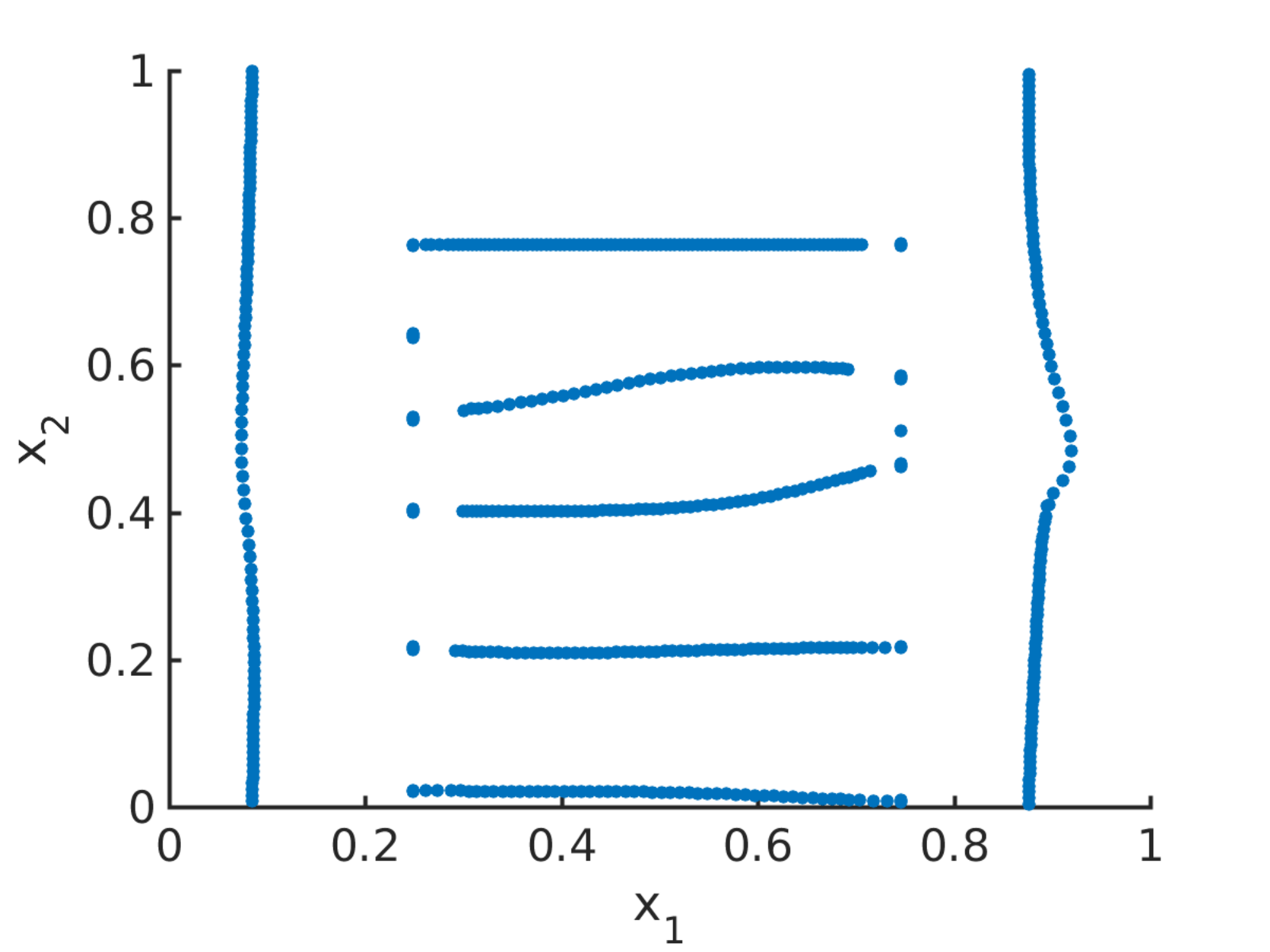}}
	\subfloat[$t=400000$]{     \includegraphics[width=0.24\textwidth]{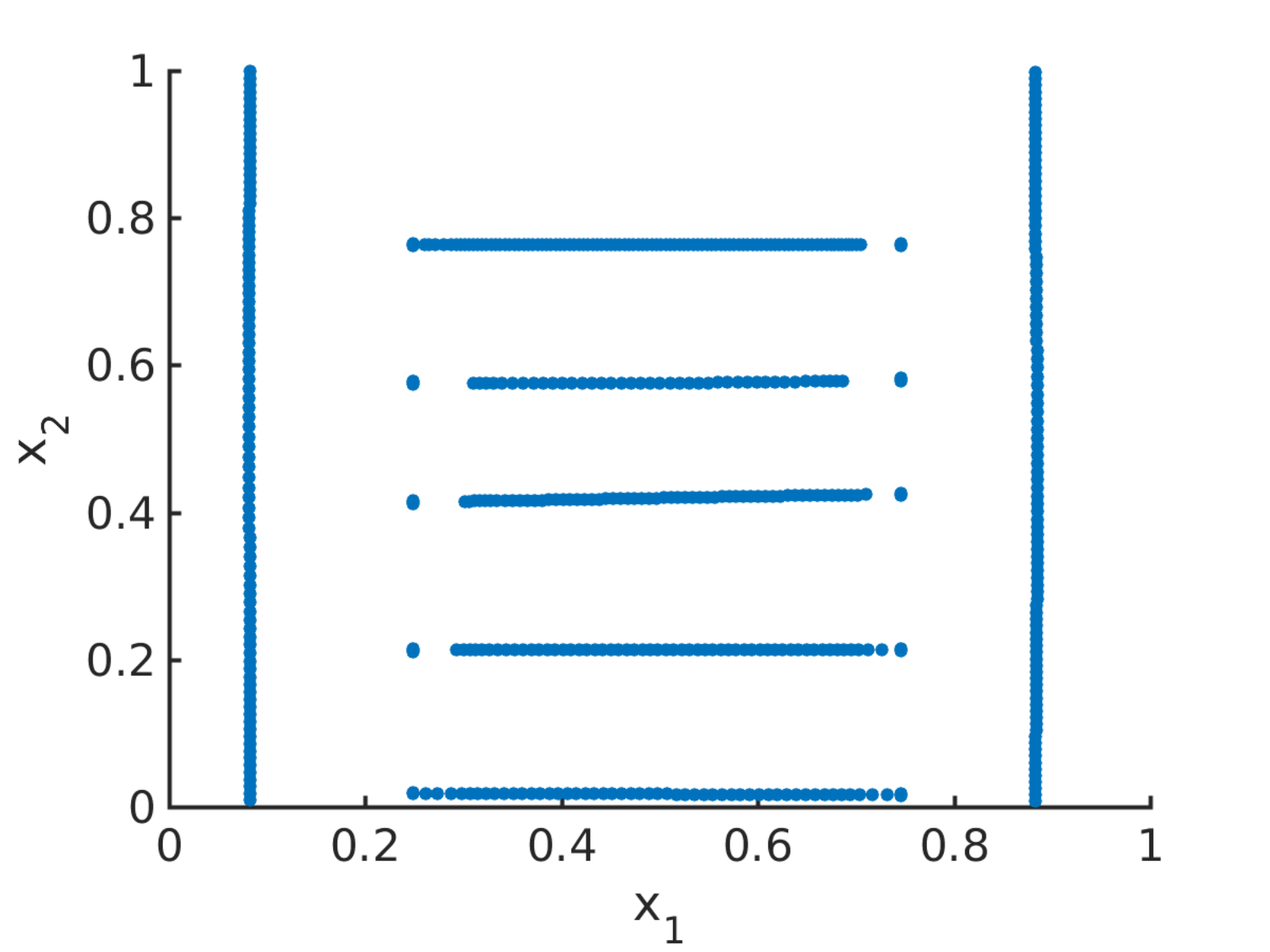}}
	\caption{Non-homogeneous tensor fields $T=T(x)$ given by $s=s(x)$ (left) and the numerical solution to the K\"{u}cken-Champod model \eqref{eq:particlemodel}  at different times $t$ for $\chi=0.2$,    $T=T(x)$, $N=600$ and randomly uniformly distributed initial data (right)}\label{fig:numericalsol_nonhomtensor_longtime}
\end{figure}

In Figure \ref{fig:numericalsol_nonhomtensor_longtime} we show the numerical solution to \eqref{eq:particlemodel}  for a piecewise spatially homogeneous tensor field, randomly uniformly distributed initial data and $N=600$, resulting in stationary line patterns along the lines of smallest stress $s=s(x)$. In particular, this tensor field is not smooth. This suggests that smoothness and periodicity are not necessary to obtain stationary solutions aligned along the lines of smallest stress. 

\begin{figure}[htbp]
	\centering
	\subfloat[$\delta=0.1$]{\includegraphics[width=0.2\textwidth]
		{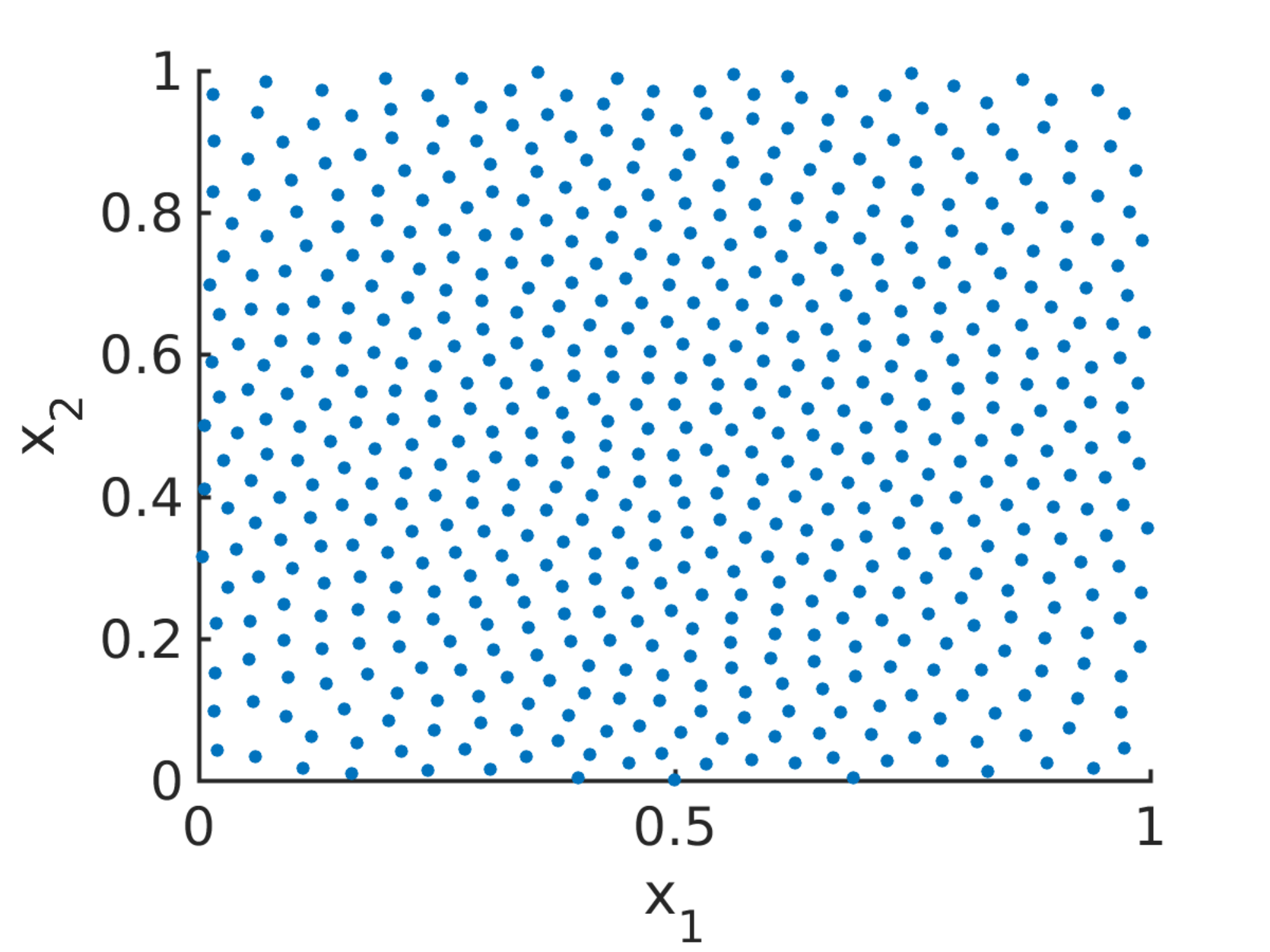}}
	\foreach \x in {3,5,...,9}{%
		\subfloat[$\delta=0.\x$]{\includegraphics[width=0.2\textwidth]
			{Particle_Euler_periodic_Repulsion0\x_Attraction_N600_circleinit_dt02.pdf}}
	}%
	\caption{Stationary solution to the K\"{u}cken-Champod model \eqref{eq:particlemodel} for  force $F(d,T)= \delta F_A(d,T)+ F_R(d)$  for different values of $\delta$ (i.e.\ different sizes of the attraction force $F_A$) and  different axis scalings where  $\chi=1$,  $N=600$ and  radially symmetric initial data (equiangularly distributed on a circle with center $(0.5, 0.5)$) so that the corresponding stationary solutions are also radially symmetric 
		and radius $0.005$}\label{fig:numericalsol_dependence_attraction}
\end{figure}

The big impact of the choice of the attraction force along the lines of largest stress can be seen by considering Figure \ref{fig:numericalsol_dependence_attraction}. Here, we  assume  that the total force is given by $F(d,T)=\delta F_A(d,T)+F_R(d)$ for $\delta\in [0,1]$ for the  spatially homogeneous tensor field $T=\chi s \otimes s+l\otimes l$ with $l=(1,0)$, $s=(0,1)$ and $\chi=1$ instead of the definition of $F$ as the sum of $F_A$ and $F_R$ in \eqref{eq:totalforce}, i.e.\ we vary the size of the attraction force and consider a radially symmetric force $F$. In Figure  \ref{fig:numericalsol_dependence_attraction} the steady states to the interaction mode \eqref{eq:particlemodel} are shown for different factors $\delta$ of the attraction force $F_A$,  where $N=600$ and  initial data distributed equiangularly on a circle with center $(0.5, 0.5)$
and radius $0.005$ is considered. One can see in Figure \ref{fig:numericalsol_dependence_attraction} that $\delta=0.1$ results in a stationary solution spread over the entire domain, while ring patterns arise as $\delta$ increases. 
The intermediate state, occurring for $\delta=0.3$, is of interest in the sequel, as it is an example of a %two-dimensional pattern in contrast to the one-dimensional ring patterns. 
more complex pattern and in particular not all the particles accumulate on one single ring as for $\delta=0.5,\delta=0.7$ and $\delta=0.9$ due to too attractive forces. 

\begin{figure}[htbp]
	\centering
	\subfloat[Normal scaling]{\includegraphics[width=0.45\textwidth]{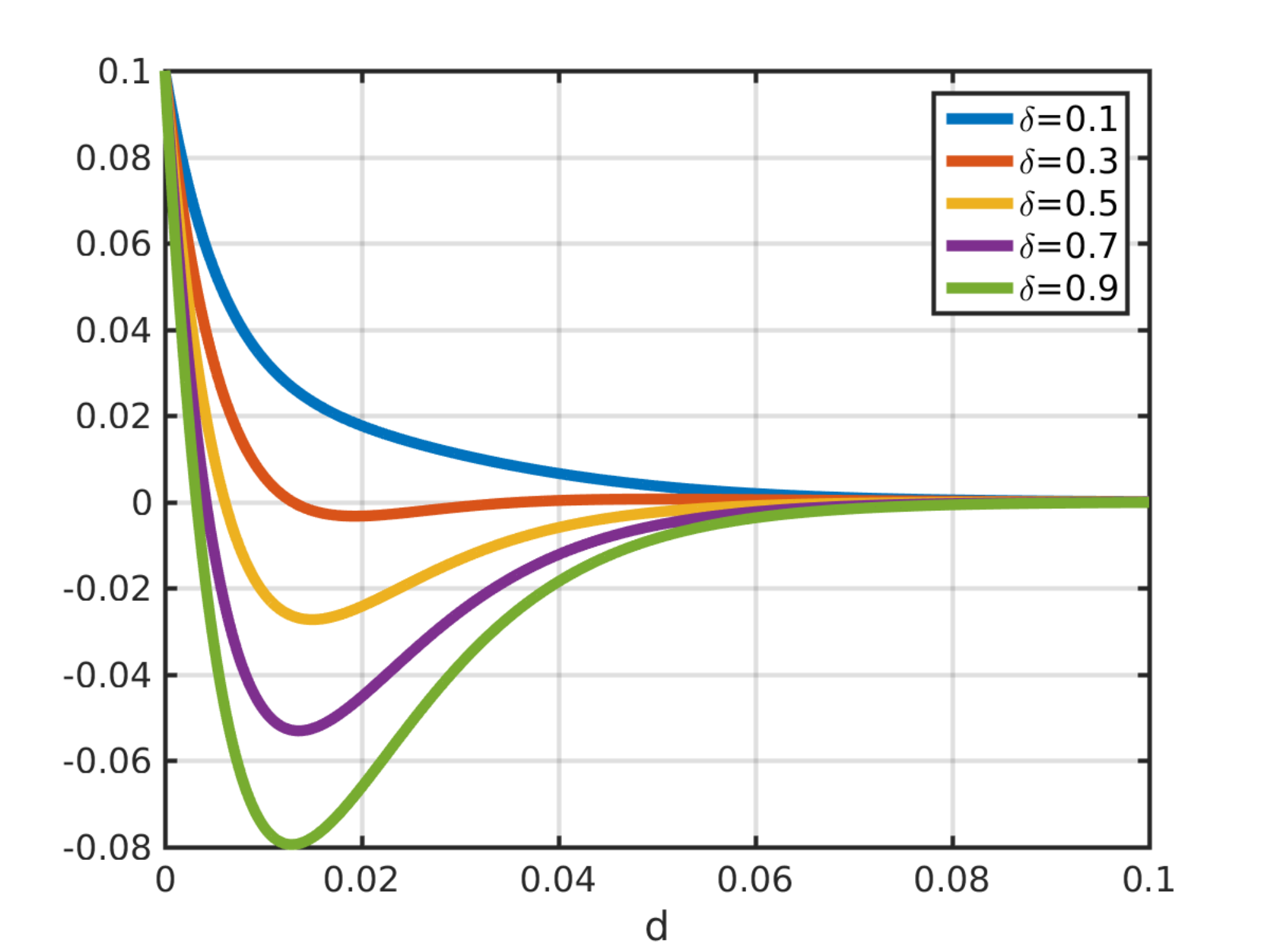}}
	\subfloat[Zoom]{\includegraphics[width=0.45\textwidth]{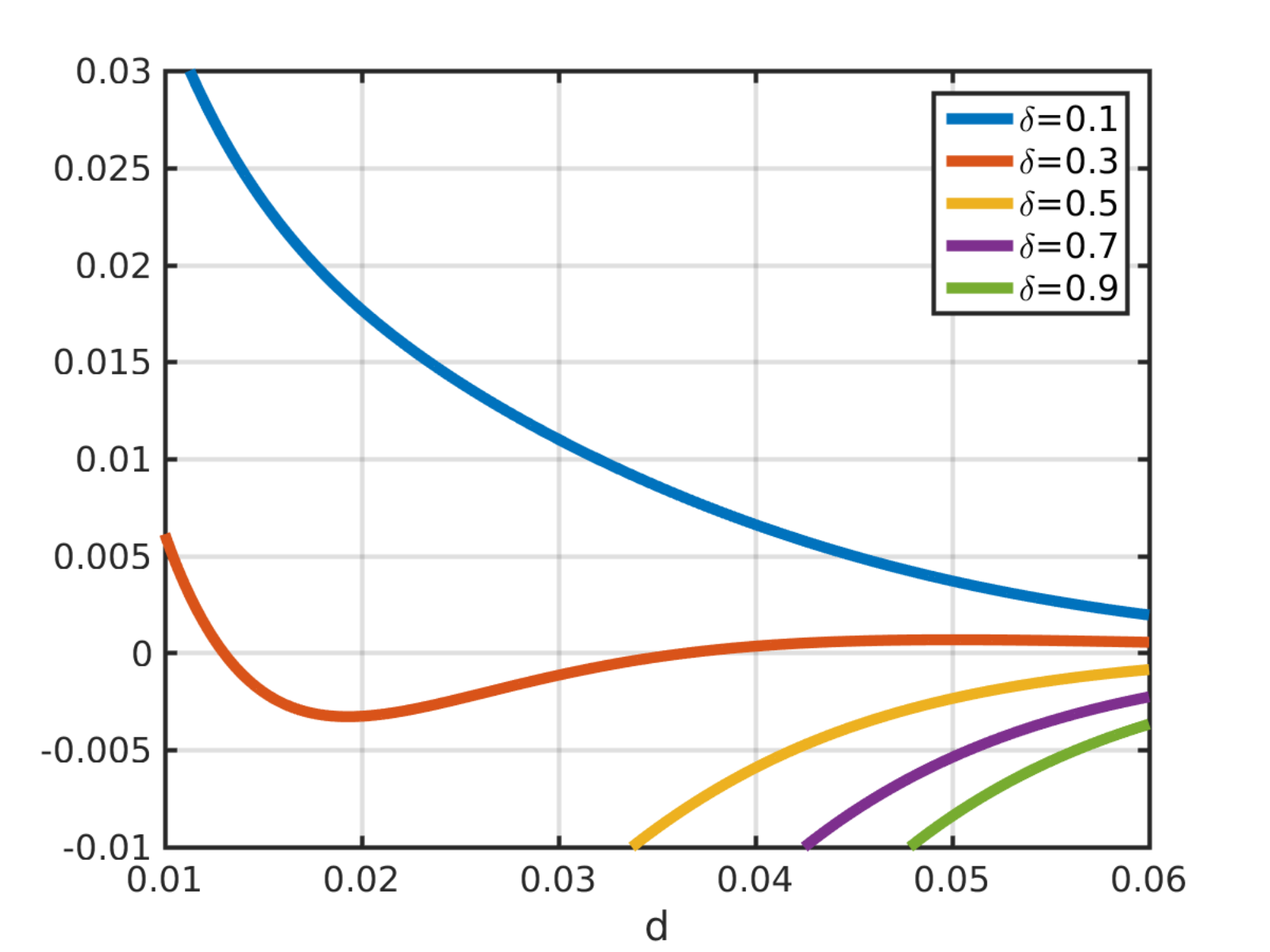}}
	\caption{Total force coefficients $\delta f_A+f_R$  along the lines of largest  stress for different values of $\delta$ and different scaling\label{fig:forcesmodelvaryattraction}}
\end{figure}

The forces considered in Figure \ref{fig:numericalsol_dependence_attraction} and given by $\delta f_A+f_R$ along the lines of largest stress are plotted in Figure \ref{fig:forcesmodelvaryattraction} for different values of $\delta$. As observed in the stationary states in Figure \ref{fig:forcesmodelvaryattraction},  the force along the lines of largest stress is purely repulsive for $\delta=0.1$, medium- and long-range attractive for $\delta\geq 0.5$, as well as medium-range attractive and long-range repulsive for $\delta=0.3$. In particular, the medium-range attractive forces for $\delta=0.3$ are significantly smaller than for larger values of $\delta$.

\section{A new model for simulating fingerprints}

Based on the analysis of stationary states in Section \ref{sec:existencesteadystates} as well as the numerical investigation of the Kücken-Champod model in Section \ref{sec:numerickcmodel} we propose a new model for simulating the formation of fingerprints based on cell interactions where  fingerprints are obtained as stationary states to our model. As a next step we propose a bio-inspired model for the creation of synthetic fingerprint patterns which can not only be used to model the formation of fingerprints as stationary solutions but also allows to adjust the ridge distances of the fingerprint lines.

\subsection{Stationary patterns}

In this section we investigate how fingerprints can be obtained as stationary solutions to the K\"{u}cken-Champod model \eqref{eq:particlemodel} where the coefficients of the repulsive and attractive forces are given by \eqref{eq:repulsionforcemodel} and \eqref{eq:attractionforcemodel}, respectively.

\subsubsection{Adaptation of the forces in the K\"{u}cken-Champod model}

Repulsive forces along the lines of smallest stress are an excellent choice to guarantee that the particles form patterns along the lines of smallest stress. Hence we can consider the repulsive coefficient function $0.2 f_A+f_R$ for the force along $s$ with the parameter values in \eqref{eq:parametervaluesRepulsionAttraction} where the coefficient functions $f_A$ and $f_R$ of the attraction and repulsion force are given by \eqref{eq:repulsionforcemodel} and \eqref{eq:attractionforcemodel},  respectively.

Short-range repulsion forces along the lines of largest stress prevent collisions of the particles and medium-range attraction forces are necessary to make the particles form aggregates. However, the long-range forces should not be attractive for modeling complex patterns since strong long-range attraction forces prevent the occurrence of multiple roughly parallel lines as stationary solutions. Motivated by the more complex stationary pattern for $\delta=0.3$ in Figure \ref{fig:numericalsol_dependence_attraction} and its desired structure of the forces along the lines of largest stress (short-range repulsive, medium-range attractive, long-range repulsive as depicted in Figure \ref{fig:forcesmodelvaryattraction}) we consider the coefficient function $0.3f_A+f_R$ along the lines of largest stress for the parameters in \eqref{eq:parametervaluesfingerprintsnew}. Hence, the total force $F$ is given by \eqref{eq:totalforce} where the repulsion force $F_R$ is defined as \eqref{eq:repulsionforcemodel} and the attraction force $F_A$ with coefficient function \eqref{eq:attractionforcemodel} has the new form
\begin{align}\label{eq:attractionforceadapted}
F_A(d=d(x_j,x_k),T(x_j))=f_A(|d|)T(x_j)=f_A(|d|)\bl 0.3(l\cdot d)l+\chi (s\cdot d)s\br
\end{align}
where we set $T(x_j)=0.3(l\cdot d)l+\chi (s\cdot d)s$ and we consider the parameter values in \eqref{eq:parametervaluesfingerprintsnew}.

\begin{figure}[htbp]
	\centering
	\begin{multicols}{4}
		\begin{minipage}{\textwidth}
			\centering
			\subfloat[$s$]{        \includegraphics[width=0.24\textwidth]{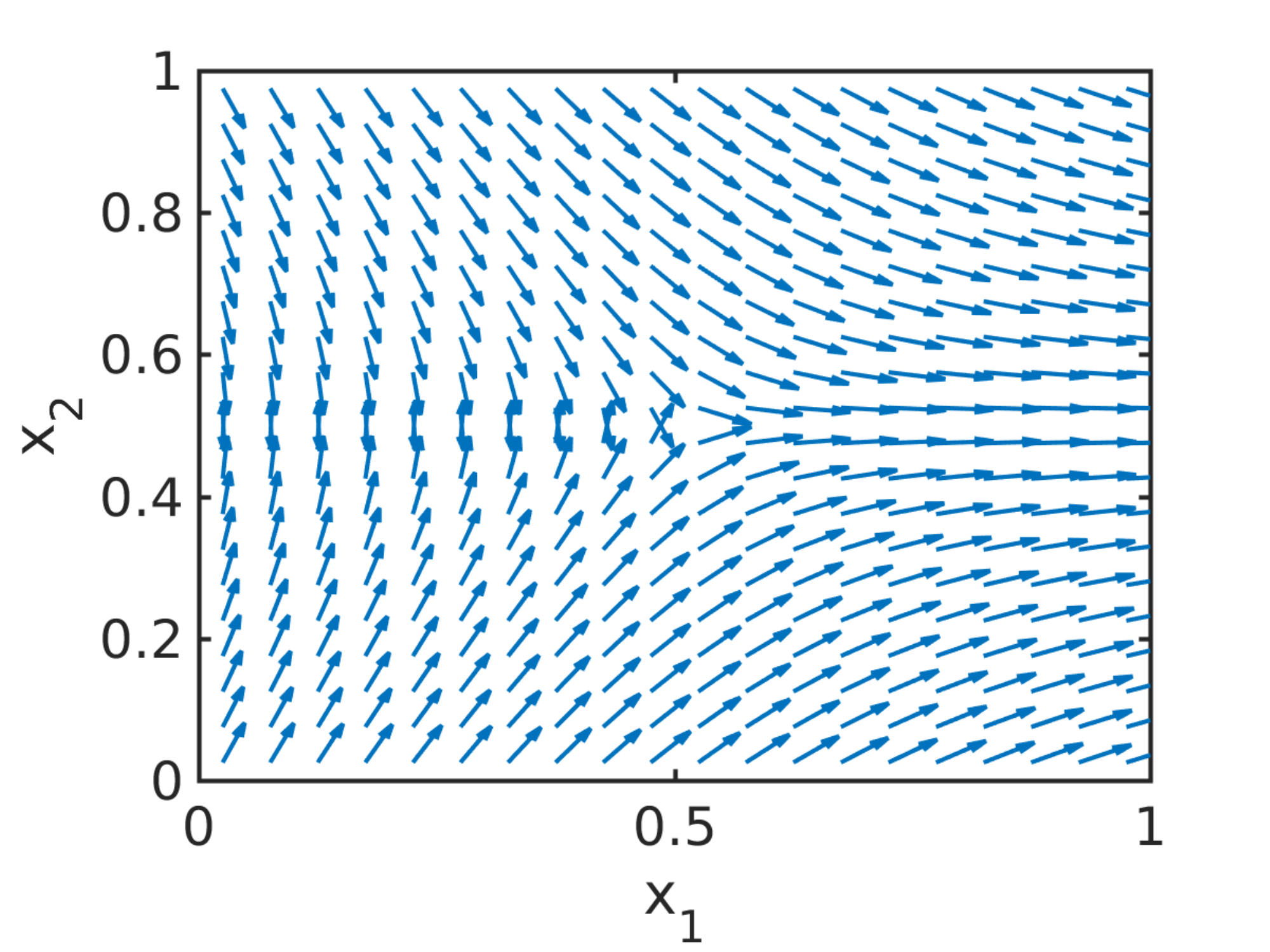}}
			\foreach \x in {40,100,...,400}{%
				\subfloat[$t=\x\cdot 10^4$]{\includegraphics[width=0.24\textwidth]{DeltaN600Time\x.pdf}}\hfill
			}
			\vspace{3mm}
			Delta
		\end{minipage}
		\begin{minipage}{\textwidth}
			\centering
			\subfloat[$s$]{        \includegraphics[width=0.24\textwidth]{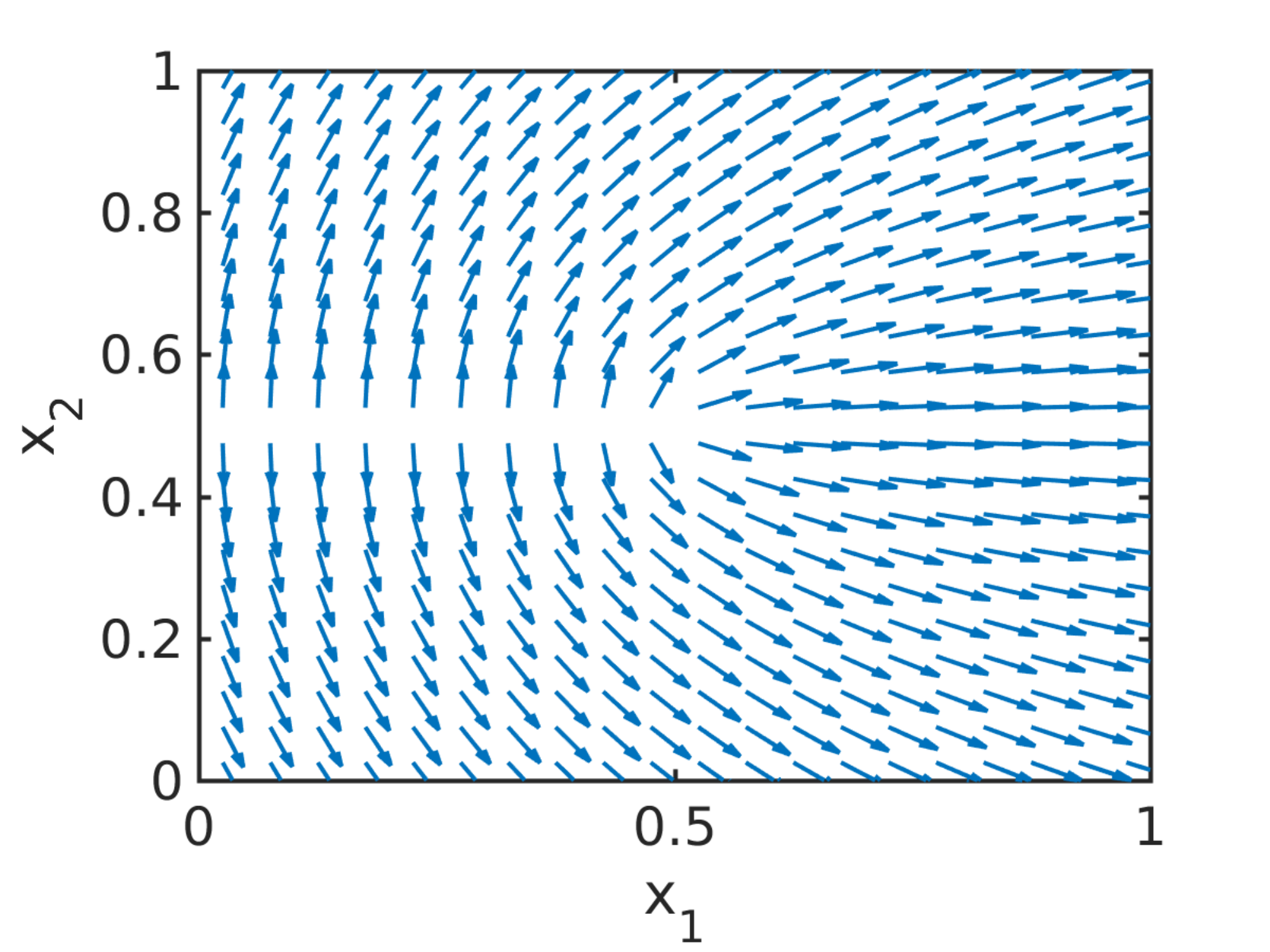}}
			\foreach \x in {40,100,...,400}{%
				\subfloat[$t=\x\cdot 10^4$]{\includegraphics[width=0.24\textwidth]{CoreN600Time\x.pdf}}\hfill
			}
			\vspace{3mm}
			Core
		\end{minipage}
	\end{multicols}
	\caption{Tensor fields $T=T(x)$ for delta (subfigures (A)-(H)) and core (subfigures (I)-(P)) given by $s=s(x)$  and the numerical solution to the extended K\"{u}cken-Champod model \eqref{eq:particlemodel} with attraction force \eqref{eq:attractionforceadapted} at different times $t$ for $\chi=0.2$, $N=600$,   $T=T(x)$ and randomly uniformly distributed initial data}\label{fig:numericalsol_nonhomtensor_deltacore}
\end{figure}

\begin{figure}[htbp]
	\centering
	\begin{multicols}{4}
		\begin{minipage}{\textwidth}
			\centering
			\subfloat[$s$]{        \includegraphics[width=0.24\textwidth]{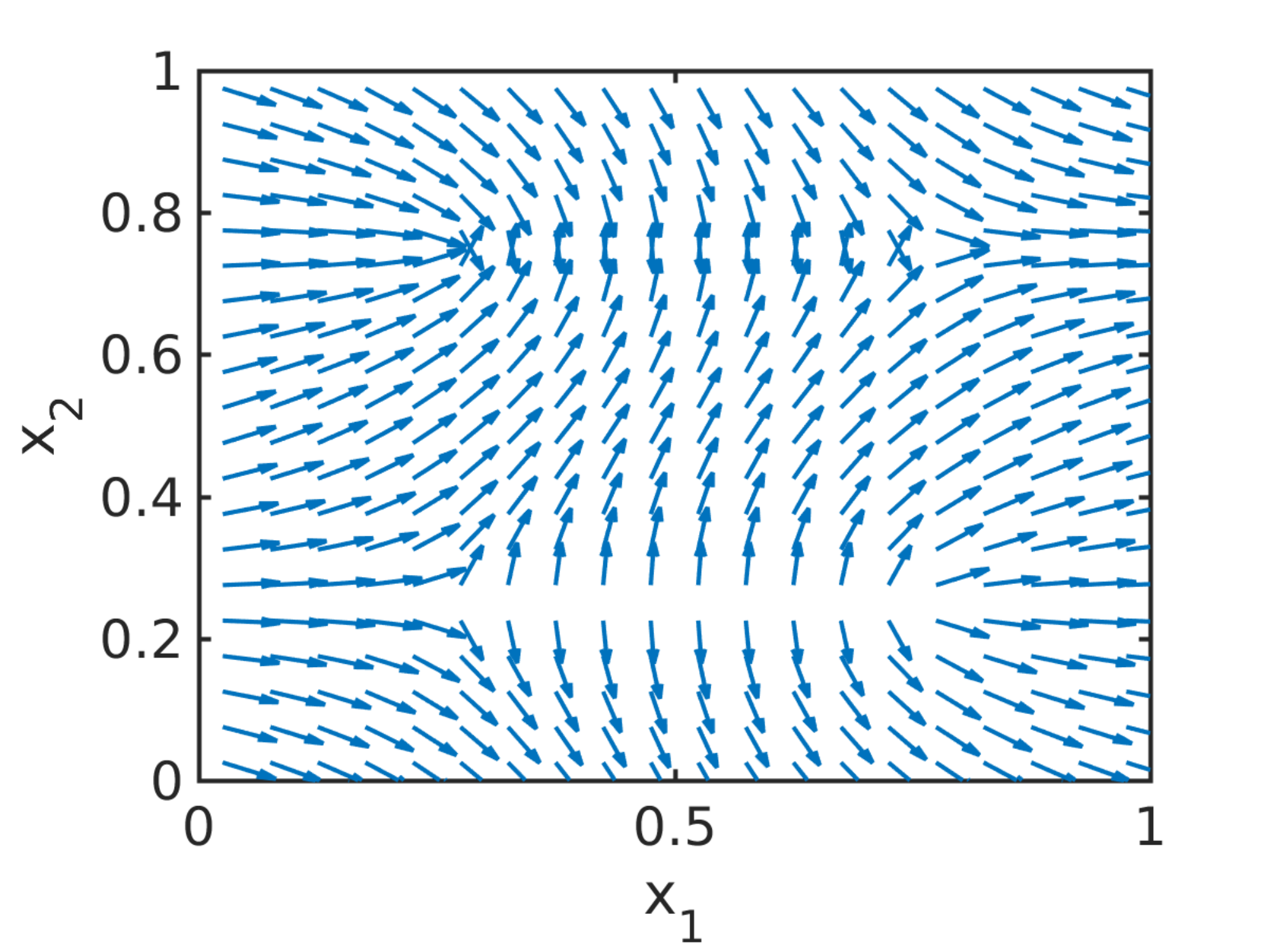}}
			\foreach \x in {40,100,...,400}{%
				\subfloat[$t=\x\cdot 10^4$]{\includegraphics[width=0.24\textwidth]{Tensorfield1N600Time\x.pdf}}\hfill
			}
			\vspace{3mm}
			Example 1
		\end{minipage}
		\begin{minipage}{\textwidth}
			\centering
			\subfloat[$s$]{        \includegraphics[width=0.24\textwidth]{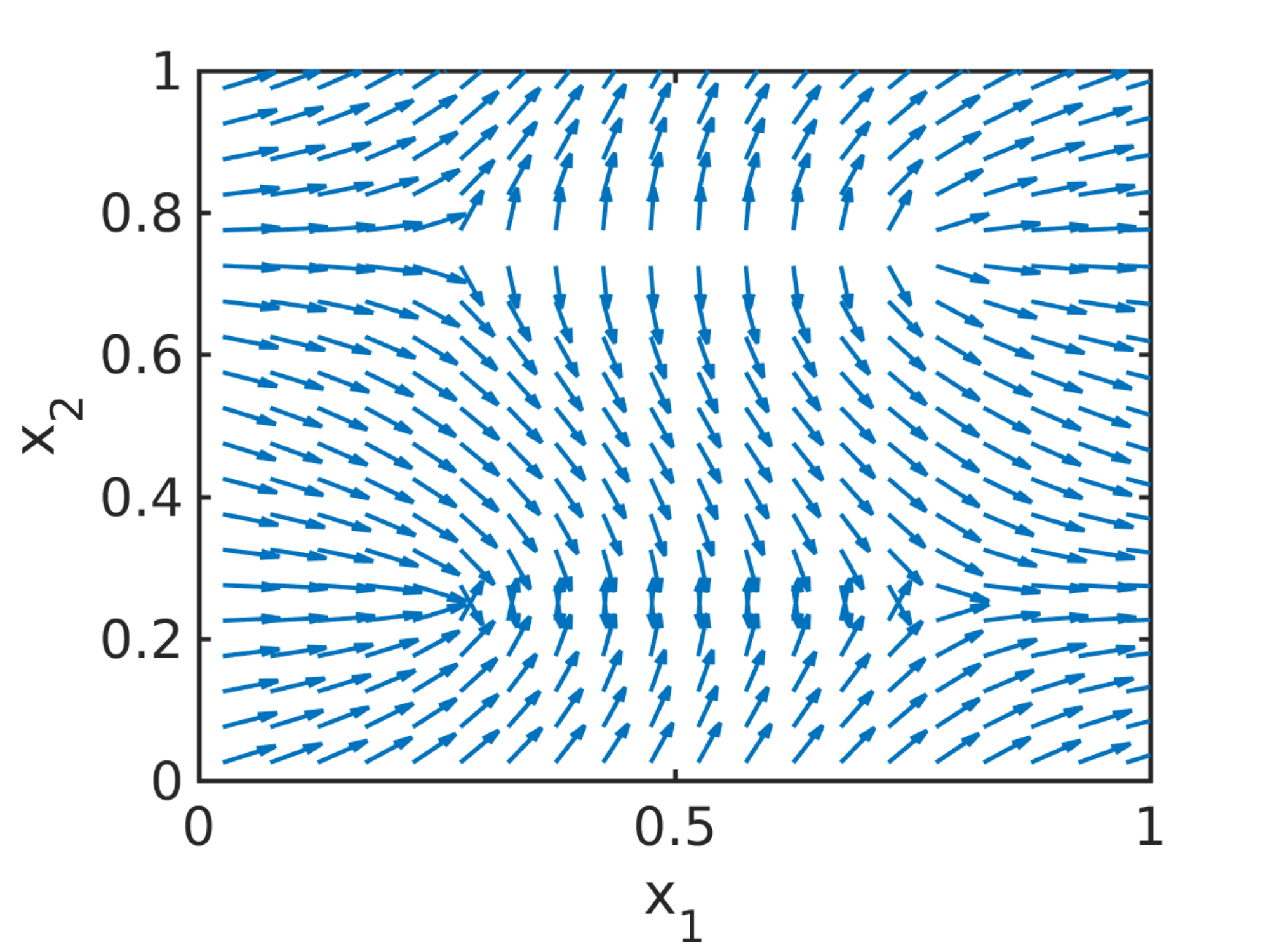}}
			\foreach \x in {40,100,...,400}{%
				\subfloat[$t=\x\cdot 10^4$]{\includegraphics[width=0.24\textwidth]{Tensorfield2N600Time\x.pdf}}\hfill
			}
			\vspace{3mm}
			Example 2
		\end{minipage}
	\end{multicols}
	\caption{Different non-homogeneous tensor fields $T=T(x)$ (Example 1 in subfigures (A)-(H), Example 2 in subfigures (I)-(P)) given by $s=s(x)$  and the numerical solution to the extended K\"{u}cken-Champod model \eqref{eq:particlemodel} with attraction force \eqref{eq:attractionforceadapted} at different times $t$ for $\chi=0.2$,   $N=600$, $T=T(x)$ and randomly uniformly distributed initial data}\label{fig:numericalsol_nonhomtensor_examples}
\end{figure}

\begin{figure}[htbp]
	\centering
	\begin{multicols}{4}
		\begin{minipage}{\textwidth}
			\centering
			\subfloat[$s$]{        \includegraphics[width=0.24\textwidth]{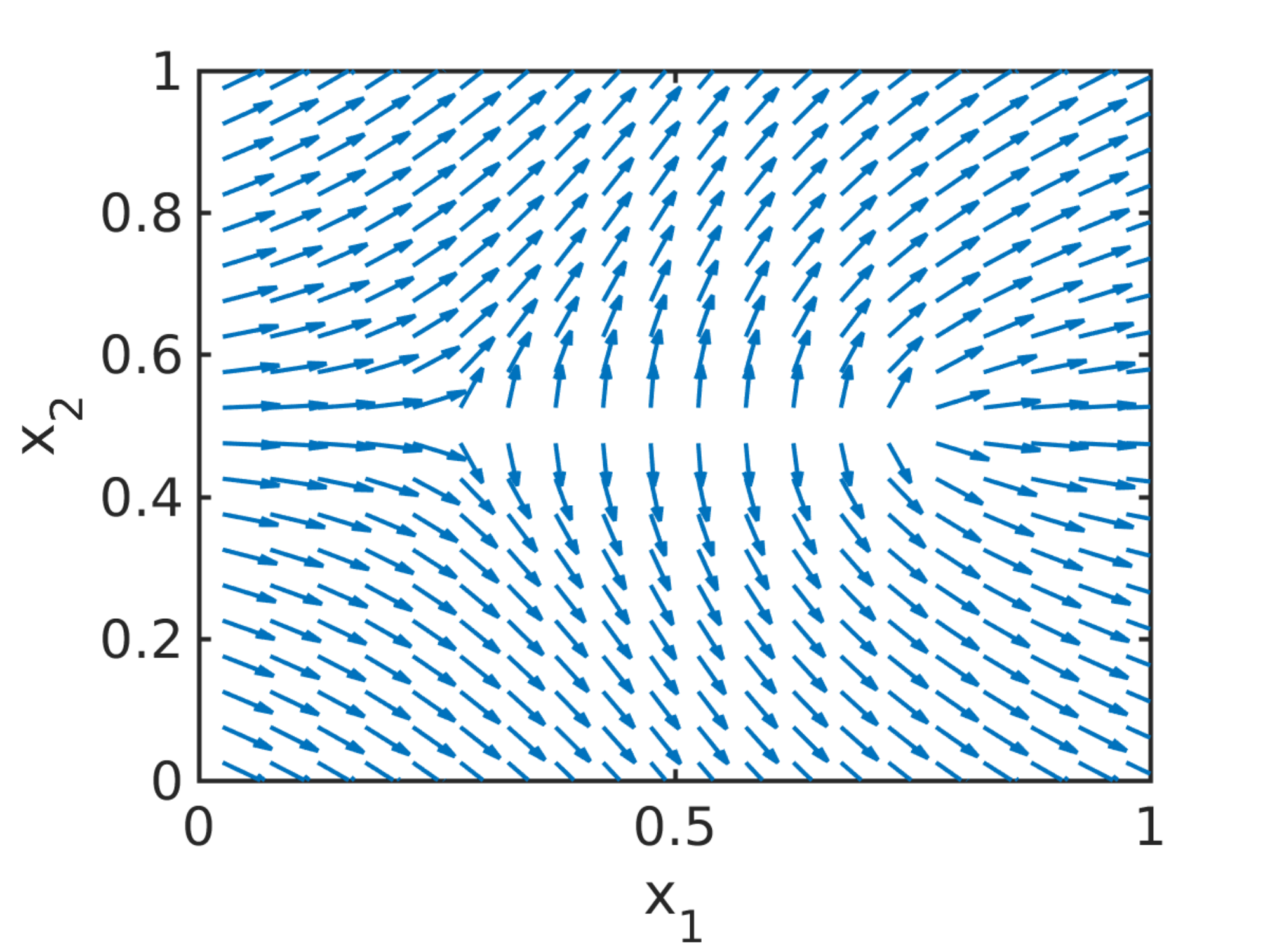}}
			\foreach \x in {40,100,...,400}{%
				\subfloat[$t=\x\cdot 10^4$]{\includegraphics[width=0.24\textwidth]{Tensorfield3N600Time\x.pdf}}\hfill
			}
			\vspace{3mm}
			Example 3
		\end{minipage}
	\end{multicols}
	\caption{Non-homogeneous, non-periodic tensor field $T=T(x)$ given by $s=s(x)$  and the numerical solution to the extended K\"{u}cken-Champod model \eqref{eq:particlemodel} with attraction force \eqref{eq:attractionforceadapted} at different times $t$ for $\chi=0.2$, $N=600$,    $T=T(x)$ and randomly uniformly distributed initial data}\label{fig:numericalsol_nonhomtensor_examples2}
\end{figure}

In Figures \ref{fig:numericalsol_nonhomtensor_deltacore}, \ref{fig:numericalsol_nonhomtensor_examples}  and \ref{fig:numericalsol_nonhomtensor_examples2} the numerical solutions for the repulsive force \eqref{eq:repulsionforce}, the attractive force \eqref{eq:attractionforceadapted} and different realistic tensor fields are illustrated. The tensor fields in Figure \ref{fig:numericalsol_nonhomtensor_deltacore} are given by a delta and a core, respectively,  introduced in Section \ref{sec:fingerprinttensor}, while we consider a combination of deltas and cores for the tensor fields in Figures \ref{fig:numericalsol_nonhomtensor_examples} and \ref{fig:numericalsol_nonhomtensor_examples2}. As desired the particles align in roughly parallel lines along the vector field $s=s(x)$ and because of long-range repulsion forces these nice patterns are not destroyed over time. Further note that the numerical solution in Figures \ref{fig:numericalsol_nonhomtensor_deltacore}, \ref{fig:numericalsol_nonhomtensor_examples} and \ref{fig:numericalsol_nonhomtensor_examples2}  is shown for very large times so that it can be regarded as stationary. In particular, this implies that the adapted forces can be used to simulate fingerprint pattern and more generally any complex patterns is in principal preserved over time.

\begin{figure}[htbp]
	\centering
	\subfloat[$0.3f_A+f_R$ for values in \eqref{eq:parametervaluesfingerprintsnew}]{\includegraphics[width=0.45\textwidth]{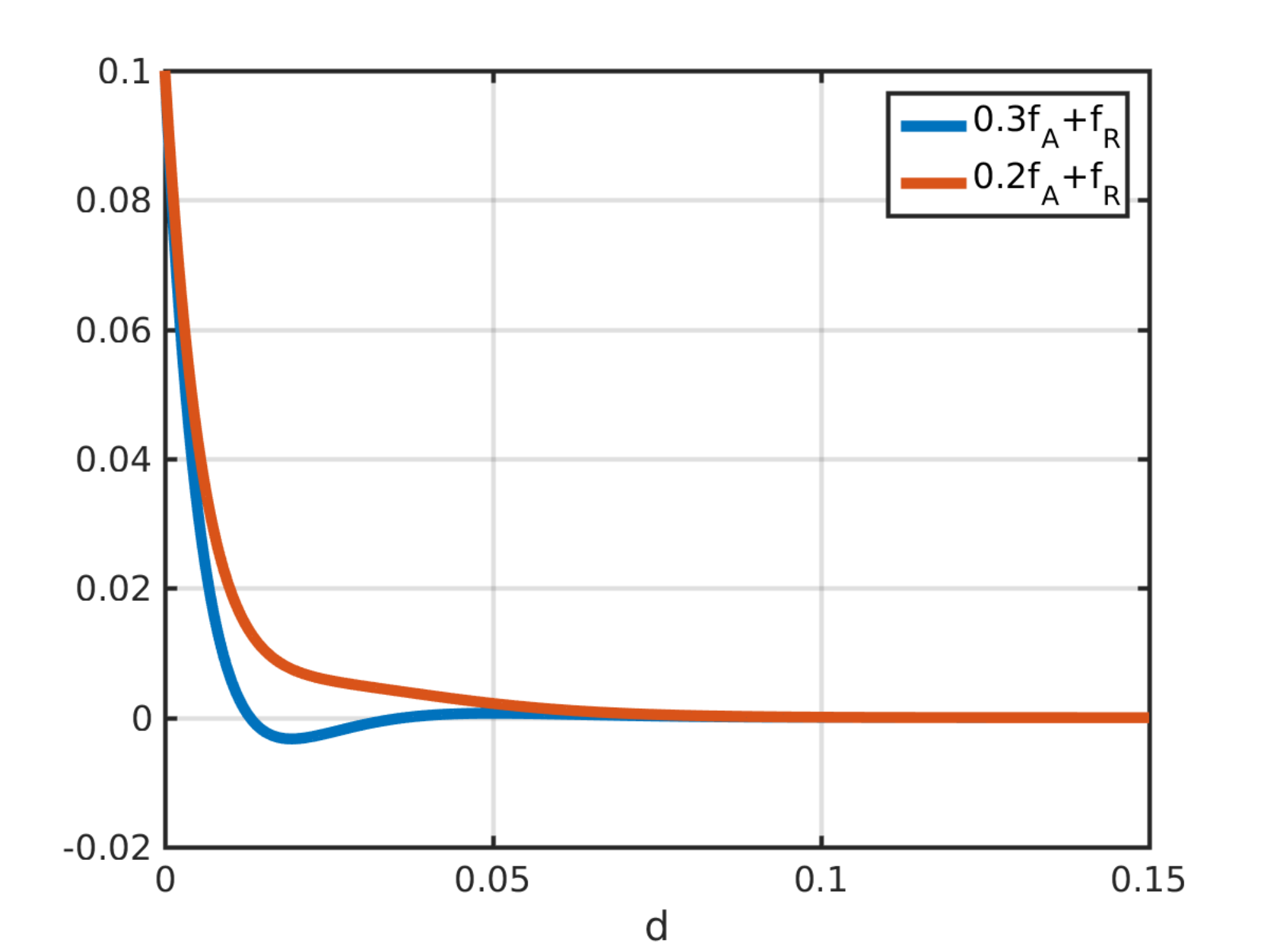}\label{fig:forcesmodeladapt}}
	\subfloat[$f_A+f_R$ for values in \eqref{eq:parametervaluesfingerprintsnew}]{\includegraphics[width=0.45\textwidth]{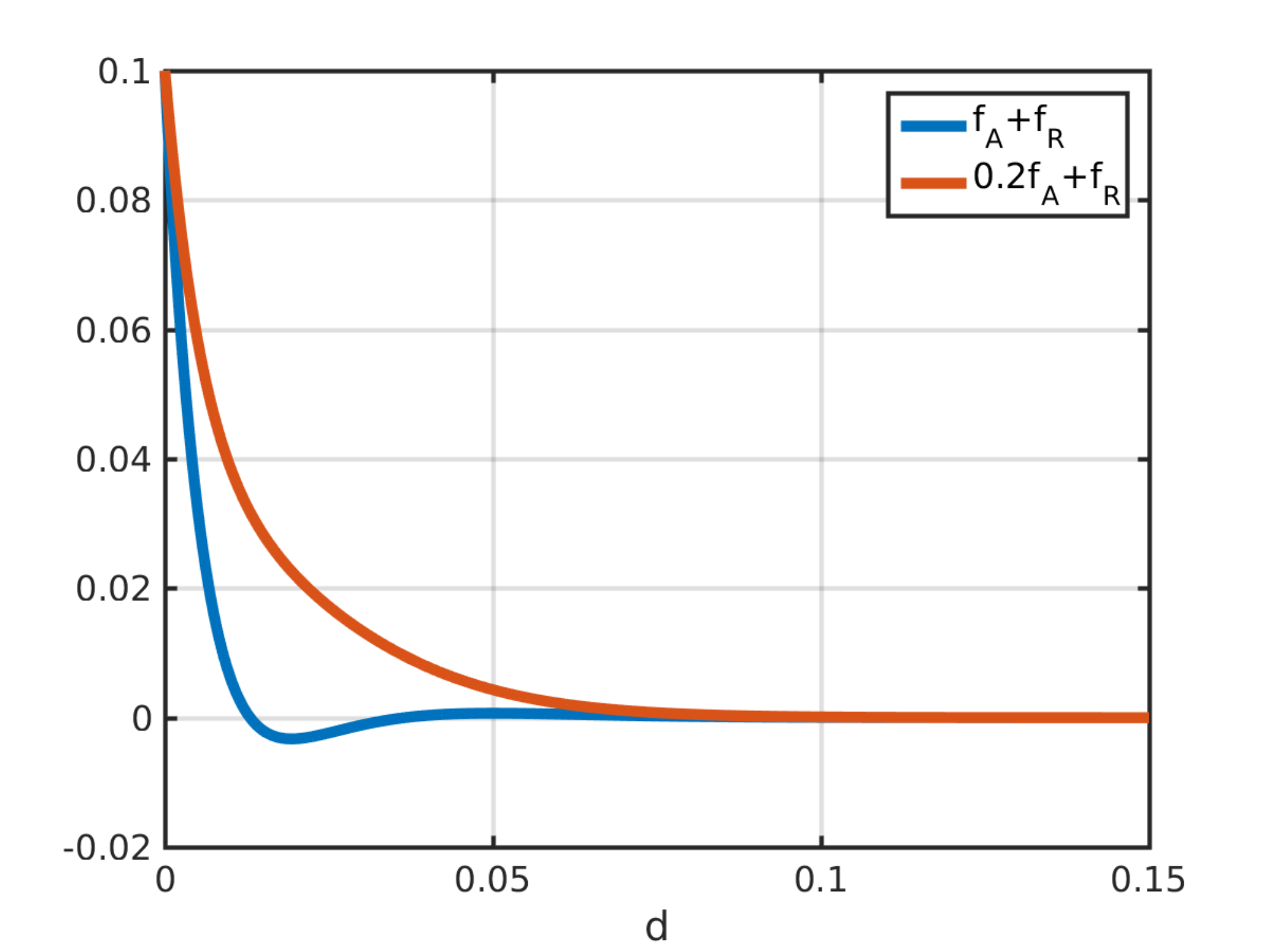}\label{fig:forcesmodelnewparameter}}
	\caption{Total force coefficients $0.2f_A+f_R$ along the lines of  smallest stress, as well as $0.3f_A+f_R$ for parameter values in \eqref{eq:parametervaluesfingerprintsnew} and $f_A+f_R$ for parameter values in \eqref{eq:parametervaluesfingerprintsnew} along the lines of  smallest stress, respectively}
\end{figure}

\begin{figure}[htbp]
	\centering
	\begin{multicols}{4}
		\begin{minipage}{\textwidth}
			\centering
			\subfloat[$s$]{        \includegraphics[width=0.24\textwidth]{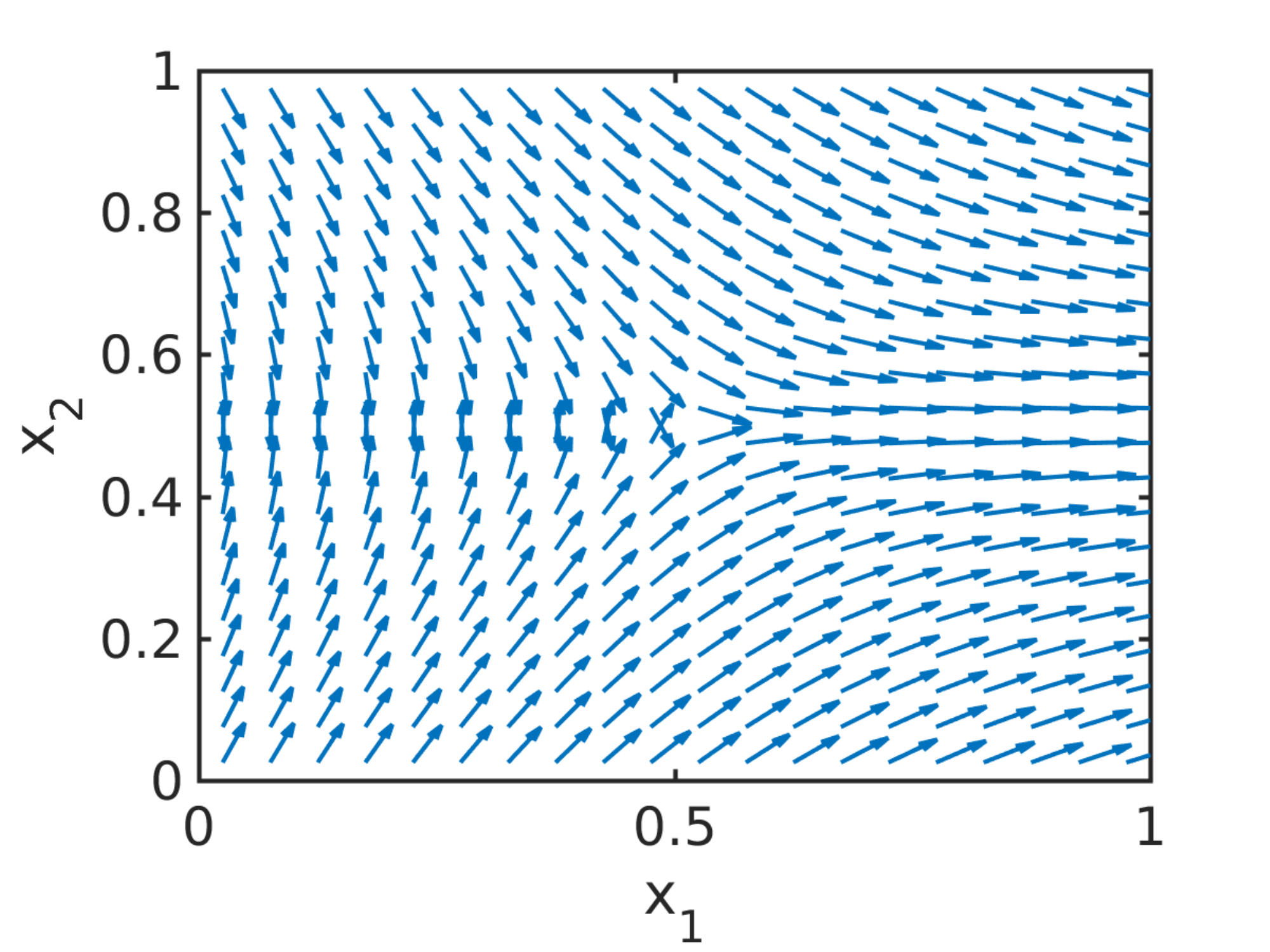}}
			\foreach \x in {40,100,...,400}{%
				\subfloat[$t=\x\cdot 10^4$]{\includegraphics[width=0.24\textwidth]{DeltaN600OriginalModelTime\x.pdf}}\hfill
			}
			\vspace{3mm}
			Example 4
		\end{minipage}
		\begin{minipage}{\textwidth}
			\centering
			\subfloat[$s$]{        \includegraphics[width=0.24\textwidth]{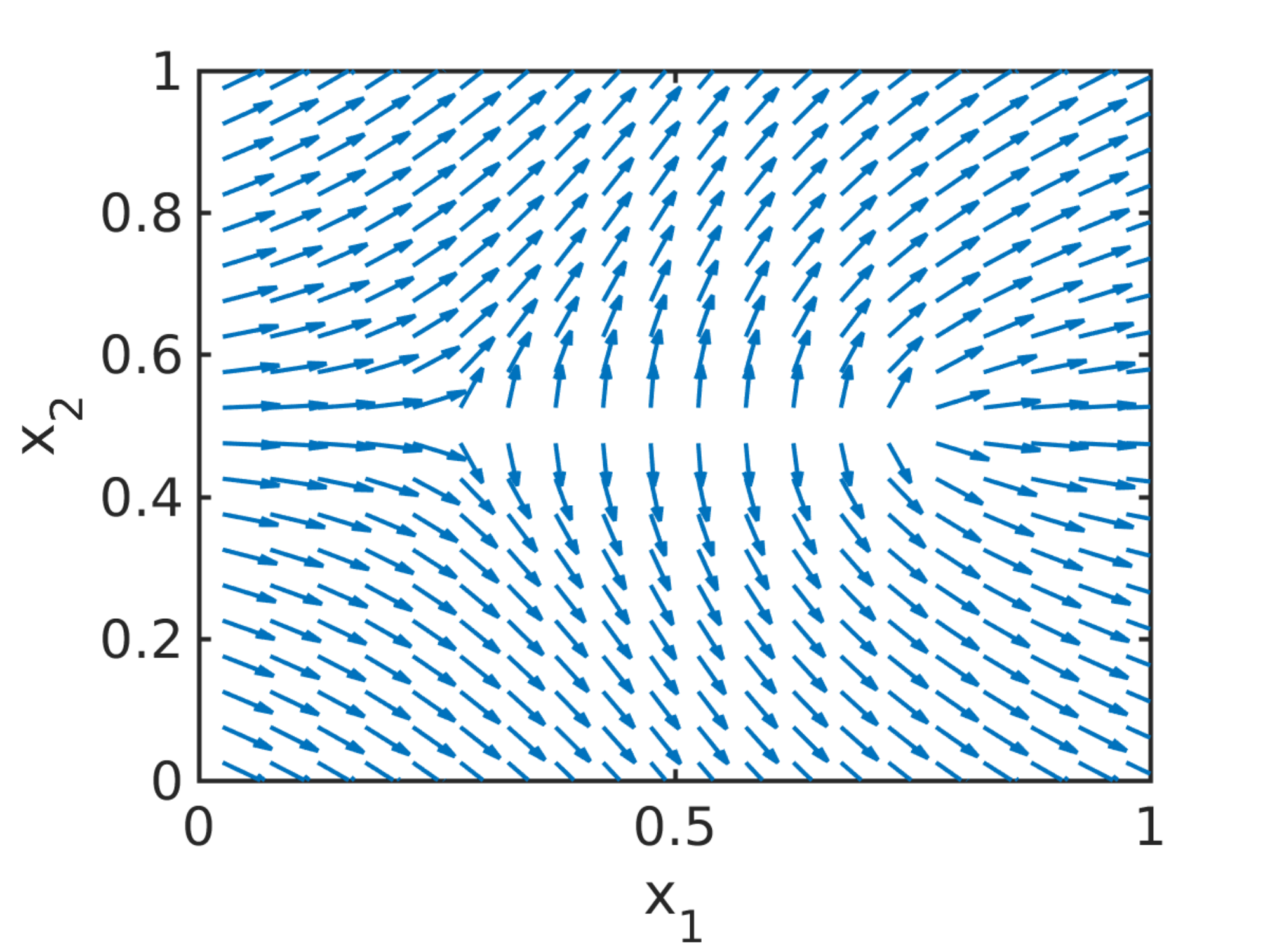}}
			\foreach \x in {40,100,...,400}{%
				\subfloat[$t=\x\cdot 10^4$]{\includegraphics[width=0.24\textwidth]{Tensorfield3N600OriginalModelTime\x.pdf}}\hfill
			}
			\vspace{3mm}
			Example 5
		\end{minipage}
	\end{multicols}
	\caption{Different non-homogeneous tensor fields $T=T(x)$ (Example 4 in subfigures (A)-(H), Example 5 in subfigures (I)-(P))  given by $s=s(x)$  and the numerical solution to the K\"{u}cken-Champod model \eqref{eq:particlemodel} for the parameters in \eqref{eq:parametervaluesfingerprintsnew}  at different times $t$ for $\chi=0.2$, $N=600$,   $T=T(x)$ and randomly uniformly distributed initial data}\label{fig:numericalsol_nonhomtensor_examples_originalmodel}
\end{figure}

After this adaptation of the forces it is desirable to use the original definition of the forces \eqref{eq:totalforce} with repulsion and attraction force given by \eqref{eq:repulsionforce} and \eqref{eq:attractionforce}, respectively, instead of an attraction force of the form \eqref{eq:attractionforceadapted}. Along $l$ the attraction force \eqref{eq:attractionforceadapted} can be regarded as $0.3 f_A$ where $f_A$ is the attraction force along $l$ in the original definition of the attraction force $F_A$ in \eqref{eq:attractionforce}. Note that the parameter $\gamma$ in the definition of the attractive force coefficient $f_A$ in \eqref{eq:attractionforcemodel} is a multiplicative constant.  Hence, we multiply the original value of $\gamma$ in \eqref{eq:parametervaluesRepulsionAttraction} by $0.3$, resulting in
\begin{align}\label{eq:parametervaluesfingerprintsnew}
\begin{split}
\alpha&=270, \quad \beta=0.1, \quad \gamma=10.5, \quad
e_A=95, \quad e_R=100, \quad \chi= 0.2,
\end{split}
\end{align} 
and  consider the original definition of the forces in \eqref{eq:totalforce}, \eqref{eq:repulsionforce} and \eqref{eq:attractionforce}. The forces along the lines of smallest and largest stress are plotted for the parameters in \eqref{eq:parametervaluesfingerprintsnew} in Figure \subref*{fig:forcesmodelnewparameter}. Note that they are of the same form as the adapted forces \eqref{eq:totalforce}, \eqref{eq:repulsionforce} and \eqref{eq:attractionforceadapted} for the original parameter values \eqref{eq:parametervaluesfingerprintsnew}, shown in Figure \subref*{fig:forcesmodeladapt}. Because of the same structure of the forces we can expect similar simulation results. In Figure \ref{fig:numericalsol_nonhomtensor_examples_originalmodel} the numerical solution is shown for two examples, a delta, as well as a combination of a core and a delta. One can clearly see that the particles align along the lines of smallest stress and the resulting patterns are preserved over time. Similarly, one can obtain any complex pattern as stationary solution to the  K\"{u}cken-Champod model \eqref{eq:particlemodel} by adapting the underlying tensor field. In particular, this implies that the K\"{u}cken-Champod model \eqref{eq:particlemodel} with forces defined by \eqref{eq:totalforce}, \eqref{eq:repulsionforce} and \eqref{eq:attractionforce} for the parameters in \eqref{eq:parametervaluesfingerprintsnew} can be used to simulate fingerprint patterns which are in principal preserved over time.

\begin{figure}[htbp]
	\centering
	\foreach \x in {5,10,...,40}{%
		\subfloat[$t=\x\cdot 10^7$]{\includegraphics[width=0.24\textwidth]{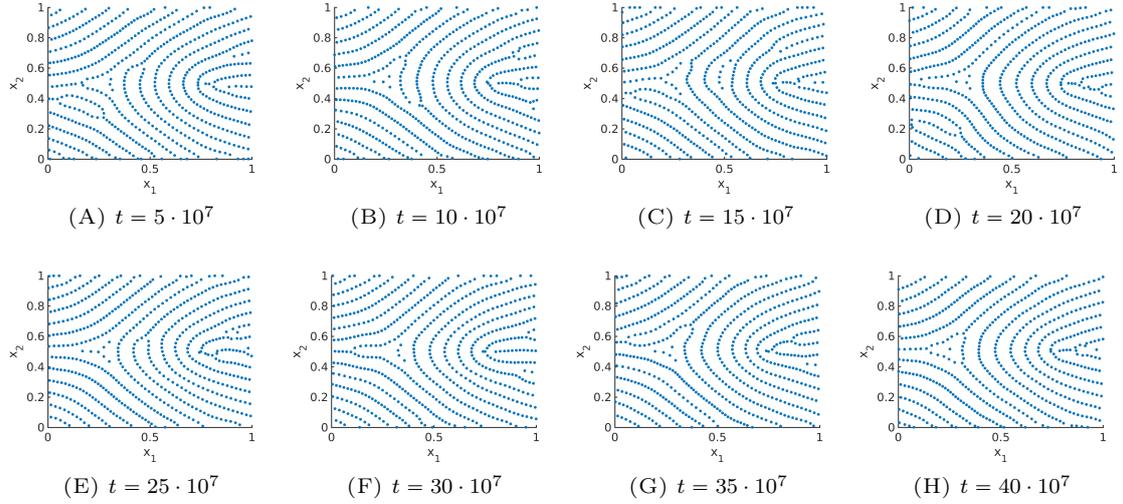}}\hfill
	}
	\caption{Long-time behavior of the numerical solution to the K\"{u}cken-Champod model \eqref{eq:particlemodel} for the parameters in \eqref{eq:parametervaluesfingerprintsnew}  at different times $t$ for $\chi=0.2$, $N=600$, the tensor field  $T=T(x)$ in Example 5 in Figure \ref{fig:numericalsol_nonhomtensor_examples_originalmodel} and randomly uniformly distributed initial data}\label{fig:numericalsol_nonhomtensor_examples_originalmodellongtime}
\end{figure}

The long-time behavior of the numerical solutions to the Kücken-Champod model \eqref{eq:particlemodel} with model parameters \eqref{eq:parametervaluesfingerprintsnew} is investigated in Figure \ref{fig:numericalsol_nonhomtensor_examples_originalmodellongtime} where the numerical solution at large times $t$ is illustrated for the tensor field in Example 5 in Figure \ref{fig:numericalsol_nonhomtensor_examples_originalmodel}. Note that the pattern changes only slightly over large time intervals, demonstrating that these patterns are close to being stationary.

\subsubsection{Pattern formation based on tensor fields from real fingerprints}\label{sec:realfingerprintpatterns}

In this section, we investigate how to simulate fingerprint patterns based on realistic tensor fields. As proposed in \cite{Merkel} the tensor field is constructed based on real fingerprint data.
The tensor field is estimated by a combination of the line sensor method \cite{LineSensor} 
and a gradient based method as described in   \cite[Section~2.1]{OrientedDiffusion}.

\begin{figure}[htbp]
	\centering
	\subfloat[Original]{\includegraphics[width=0.32\textwidth]{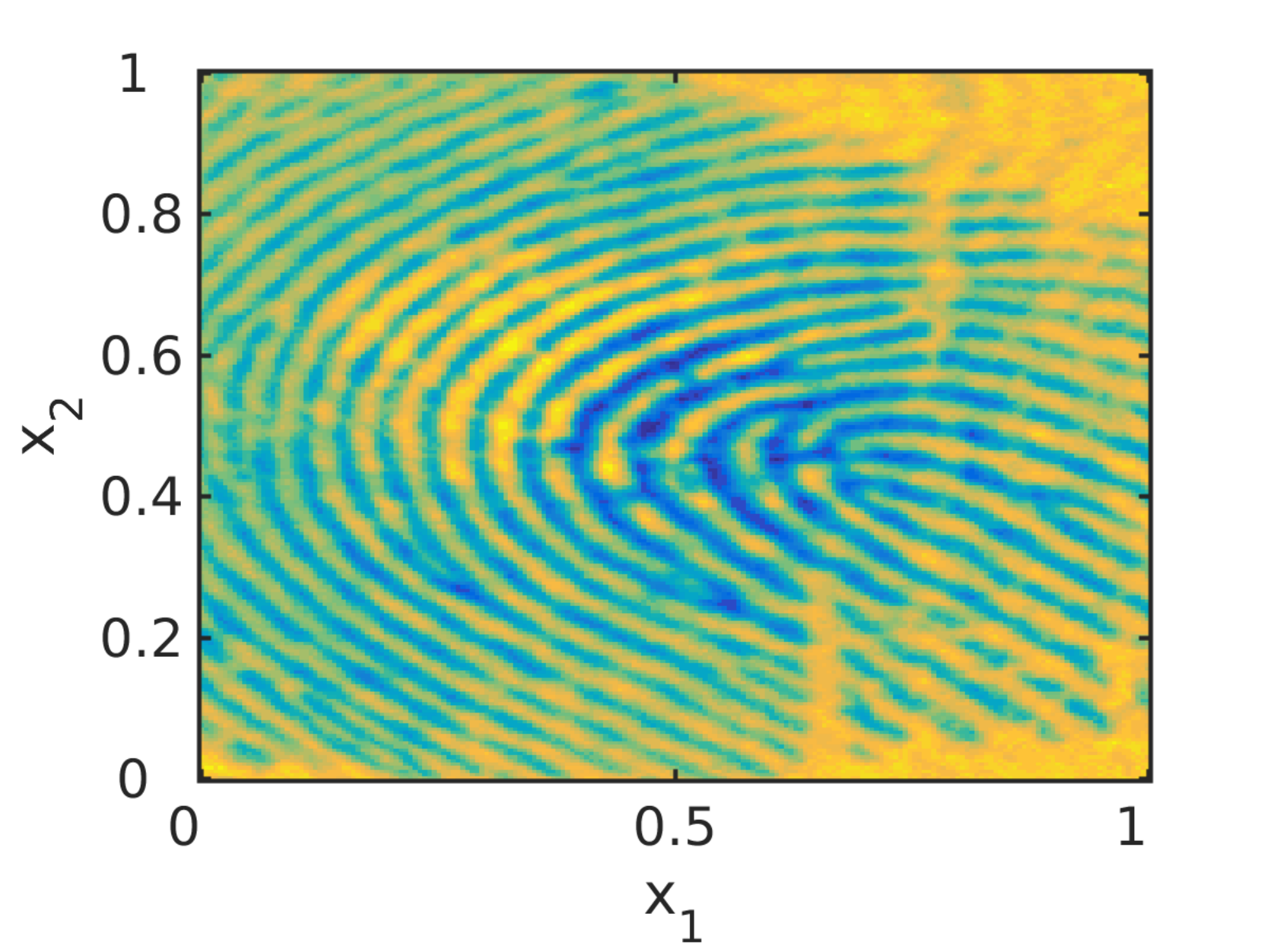}}\hfill
	\subfloat[$\theta$]{\includegraphics[width=0.32\textwidth]{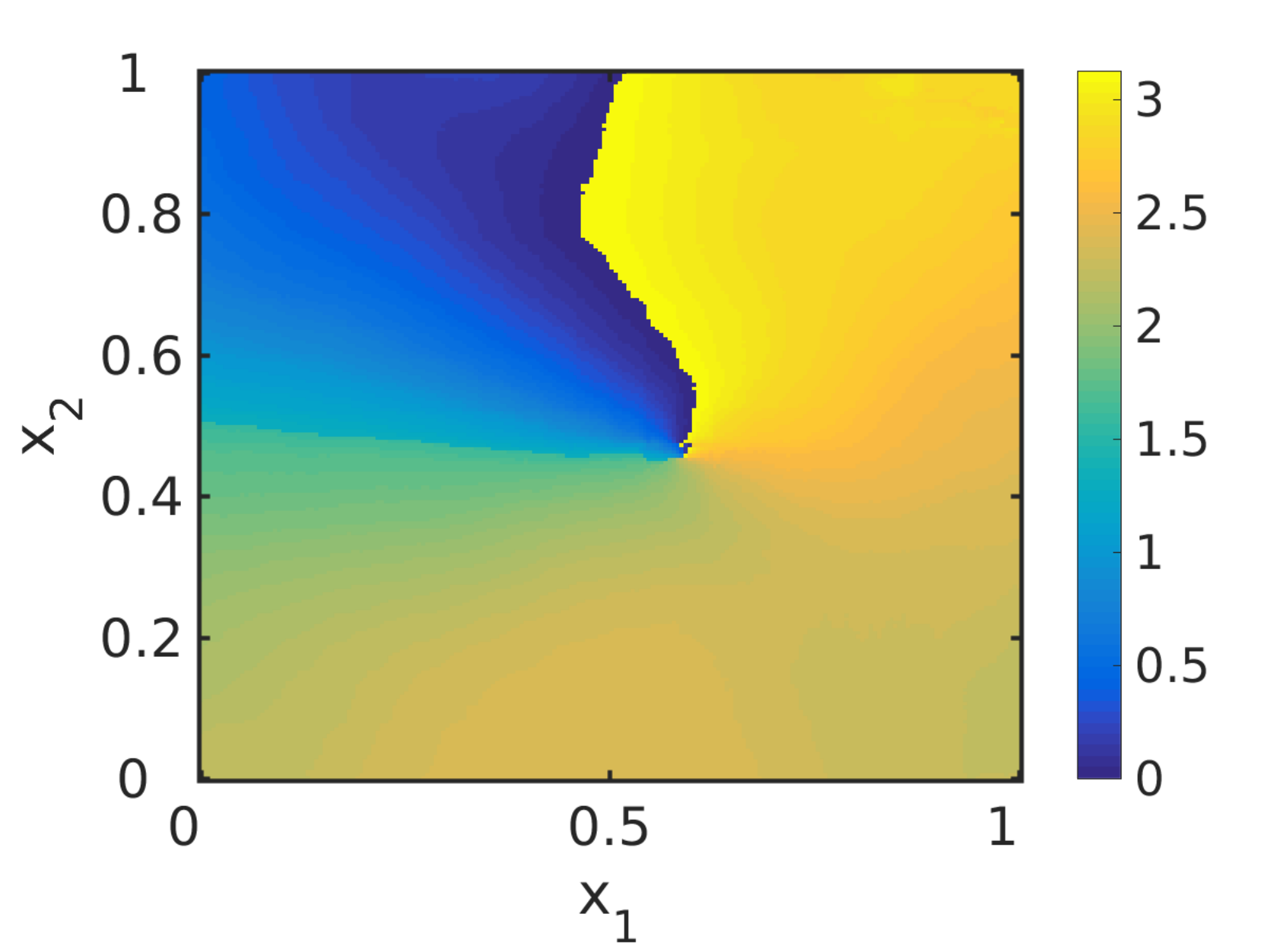}}\hfill
	\subfloat[$s$]{\includegraphics[width=0.32\textwidth]{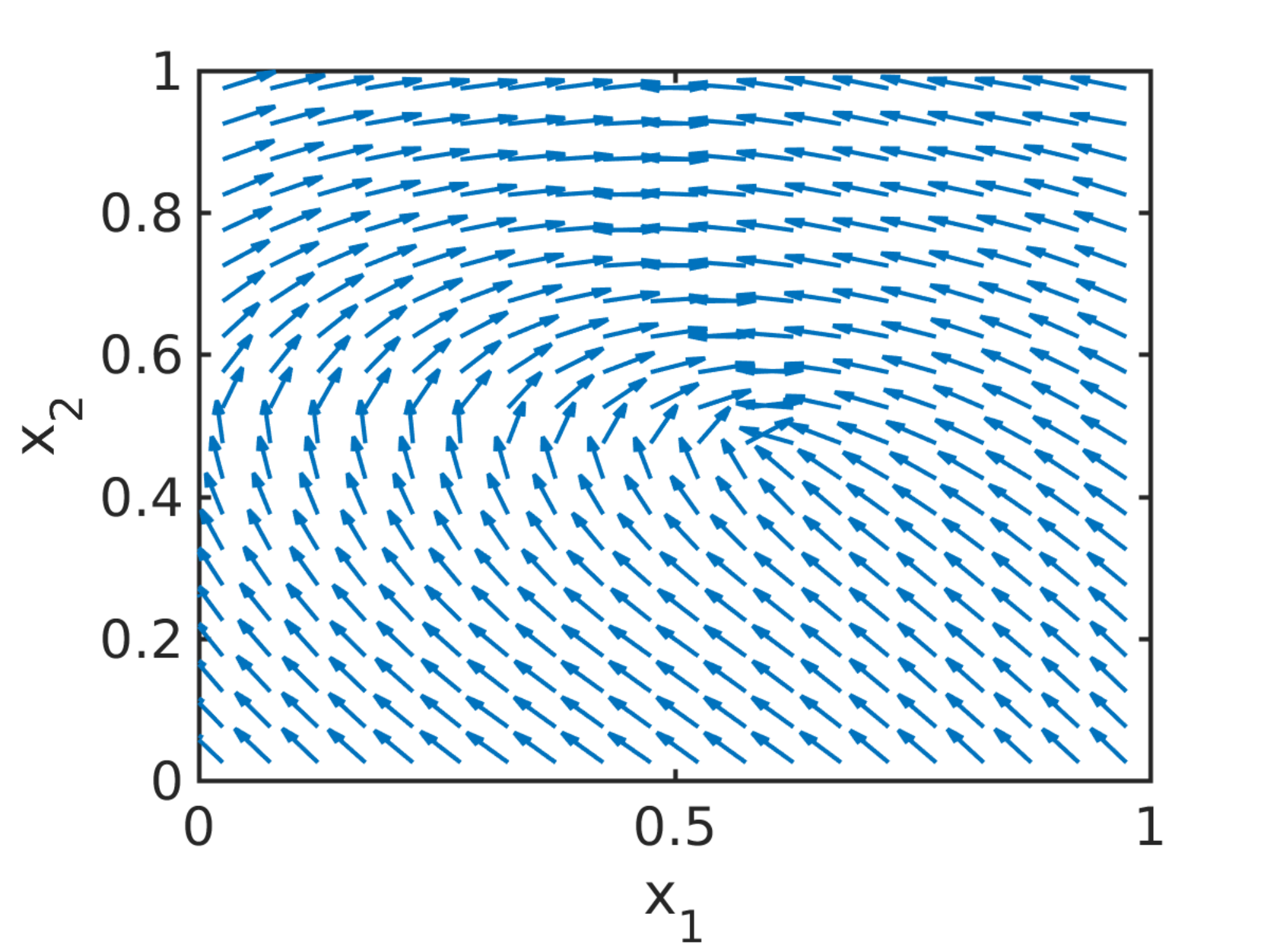}}\hfill
	\caption{Original fingerprint image as well as  angles and lines of smallest stress $s=s(x)$ for the reconstructed tensor field $T=T(x)$}\label{fig:realtensorfield}
\end{figure}

Given some real fingerprint data the aim is to construct the vector field $s=s(x)$ for all $x\in\Omega$ as the tangents to the given fingerprint lines. This is based on the idea that the lines of smallest stress are given by $s$ and the solution to the interaction model \eqref{eq:particlemodel} aligns along $s$. Let  $\theta=\theta(x)$ denote the angle between the vertical axis and the direction of lines of smallest stress $s=s(x)$ at location $x$, then it is sufficient to consider the principal arguments $\theta\in[0,\pi)$ only. Note that for any $x\in \Omega$ and any given  $\theta(x)$ we can reconstruct  $s(x)$ as $(\cos(\theta(x),\sin(\theta(x)))$ since $s(x)$ are defined to be unit vectors. In Figure \ref{fig:realtensorfield}  fingerprint data, the estimated arguments $\theta$ for constructing the tensor field and the lines of smallest stress $s=s(x)$ of the tensor field are shown. Note that the lines of smallest stress $s=s(x)$ of the tensor field and the fingerprint lines in the real fingerprint image coincide. 

\begin{figure}[htbp]
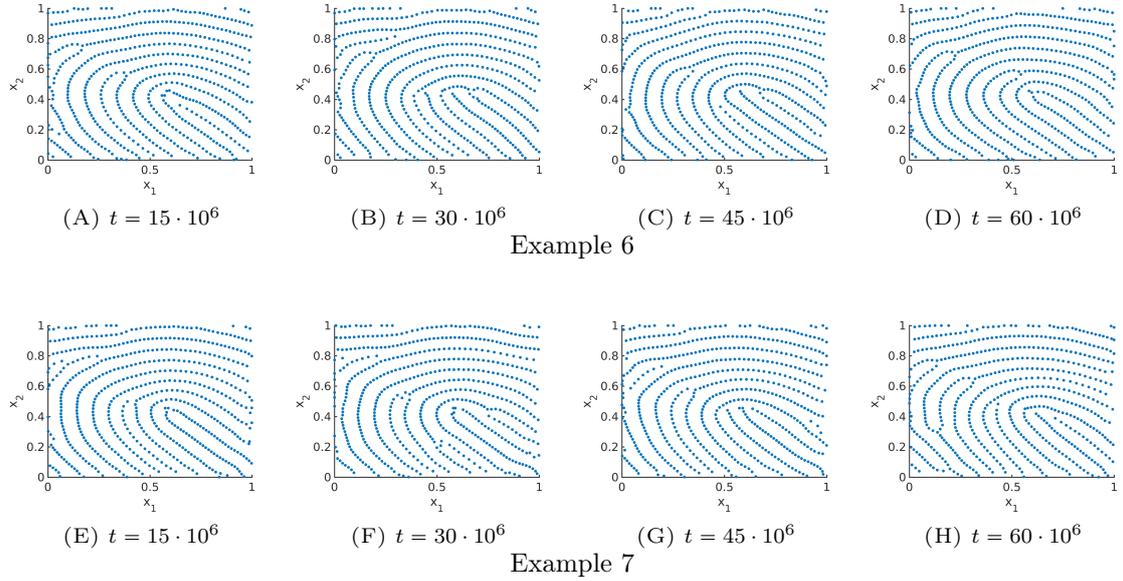

	\centering
	\begin{multicols}{4}
		\begin{minipage}{\textwidth}
			\centering
			\foreach \x in {15,30,...,60}{%120}{%
				\subfloat[$t=\x\cdot 10^6$]{\includegraphics[width=0.24\textwidth]{Example1OriginalModelTime\x00.pdf}}\hfill
			}
			\vspace{3mm}
			Example 6
		\end{minipage}
		\begin{minipage}{\textwidth}
			\centering
			\foreach \x in {15,30,...,60}{%
				\subfloat[$t=\x\cdot 10^6$]{\includegraphics[width=0.24\textwidth]{Example2OriginalModelTime\x00.pdf}}\hfill
			}
			\vspace{3mm}
			Example 7
		\end{minipage}
	\end{multicols}
	\caption{Numerical solution to the K\"{u}cken-Champod model \eqref{eq:particlemodel} for the parameters in \eqref{eq:parametervaluesfingerprintsnew}  at different times $t$ for $\chi=0.2$, the realistic tensor field  $T=T(x)$ in Figure \ref{fig:realtensorfield} and two realisations of randomly uniformly distributed initial data}\label{fig:numericalsol_nonhomtensor_examples_originalmodelrealtensor}
\end{figure}

%Given this tensor field the solution to the K\"ucken-Champod model \eqref{eq:particlemodel} can be determined numerically. 
Considering the tensor field $T=T(x)$ shown in Figure \ref{fig:realtensorfield}   the associated numerical solution is plotted for two realizations of uniformly distributed initial data in Figure \ref{fig:numericalsol_nonhomtensor_examples_originalmodelrealtensor}. One can clearly see that the particles align along the lines of smallest stress $s=s(x)$. Besides, Figure \ref{fig:numericalsol_nonhomtensor_examples_originalmodelrealtensor} illustrates that we obtain similar, but not exactly the same patterns for different realisations of random uniformly distributed initial data. This is consistent with the well-known fact that everyone has unique fingerprints and even the fingerprints of twins can  be distinguished even if the general patterns may seem to be quite similar at first glance \cite{champod2016fingerprints}. 

\begin{figure}[htbp]
	\centering
	\includegraphics[width=0.4\textwidth]{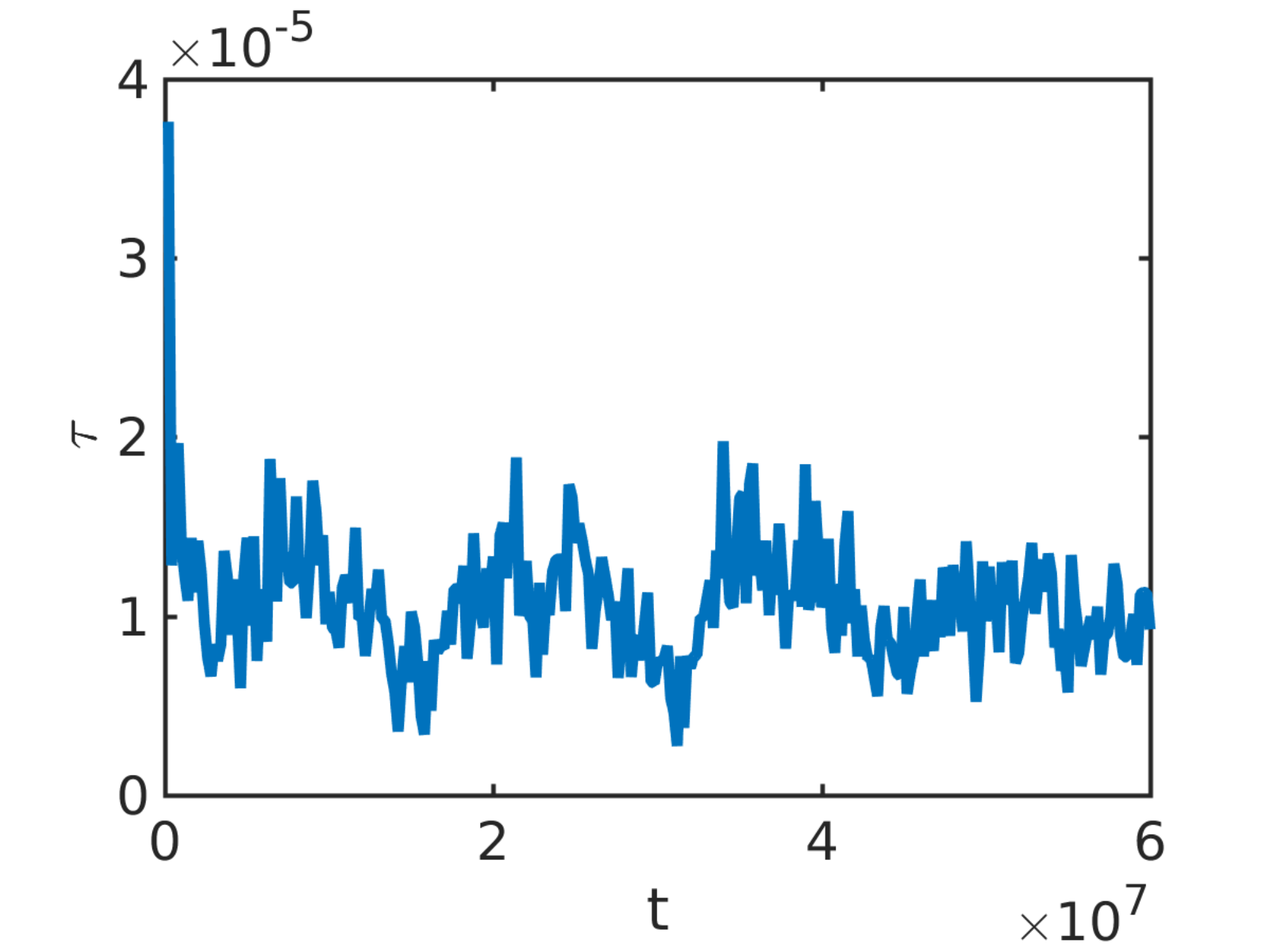}
	\caption{Error $\tau$ in \eqref{eq:changePositionOverTime} between successive time steps for the numerical solution in Example 6 in Figure \ref{fig:numericalsol_nonhomtensor_examples_originalmodelrealtensor} to the K\"{u}cken-Champod model \eqref{eq:particlemodel} for the parameters in \eqref{eq:parametervaluesfingerprintsnew}  at different times $t$ and the realistic tensor field  $T=T(x)$ in Figure \ref{fig:realtensorfield} }\label{fig:imagefingerprinttensorex1error}
\end{figure}

To quantify the distance to the steady state we consider the change of the positions $x_j$ of the particles in successive time steps, given by 
\begin{align}\label{eq:changePositionOverTime}
\tau(t)=\sum_{j=1}^N \|x_j(t+\Delta t)-x_j(t)\|_{L^1}.
\end{align}
In Figure \ref{fig:imagefingerprinttensorex1error} we show the error $\tau$ between successive time steps for the numerical solution in Example 6 in Figure \ref{fig:numericalsol_nonhomtensor_examples_originalmodelrealtensor} to the K\"{u}cken-Champod model \eqref{eq:particlemodel}. After a sharp initial decrease the total change in positions of the particles  is approximately $1.0\cdot 10^{-5}$, i.e.\ the movement of the particles is roughly $1.7\cdot 10^{-8}$ between time steps.

\subsubsection{Interpretation  of the  pattern formation}
\begin{figure}[htbp]
	\centering
	\includegraphics[width=0.45\textwidth]{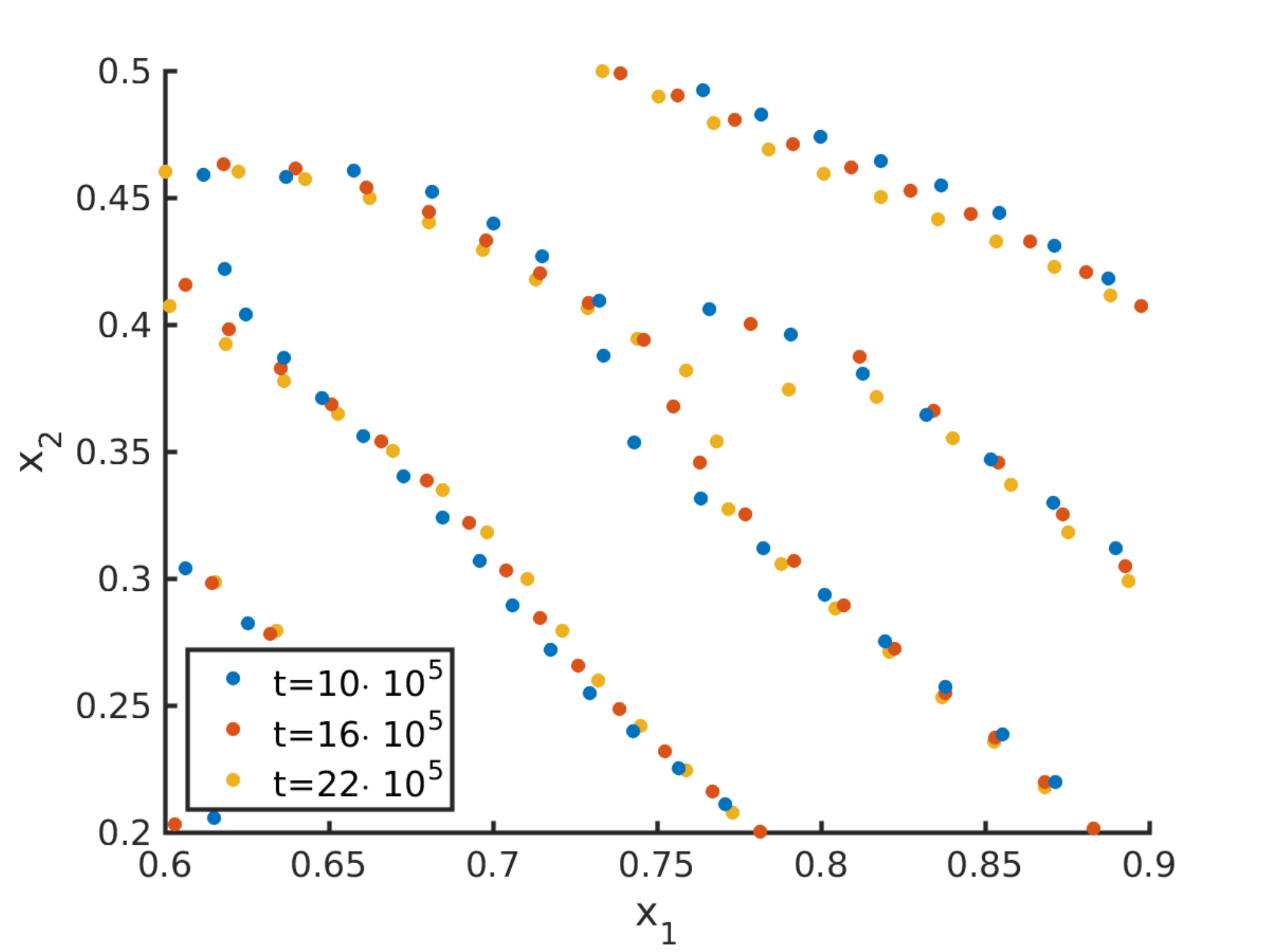}
	\caption{Evolution of the bifurcations in the numerical solution to the K\"{u}cken-Champod model \eqref{eq:particlemodel} for the parameters in \eqref{eq:parametervaluesfingerprintsnew} for the non-homogeneous tensor field $T=T(x)$  in Example 6 in Figure \ref{fig:numericalsol_nonhomtensor_examples_originalmodelrealtensor} at different times $t$ and randomly uniformly distributed initial data}\label{fig:Example1OriginalModelBifurcations}
\end{figure}

In the simulations for  spatially homogeneous tensor fields in \cite{patternformationanisotropicmodel} as well as for realistic tensor fields in Figures \ref{fig:numericalsol_nonhomtensor_deltacore}, \ref{fig:numericalsol_nonhomtensor_examples}, \ref{fig:numericalsol_nonhomtensor_examples2}, \ref{fig:numericalsol_nonhomtensor_examples_originalmodel}, \ref{fig:numericalsol_nonhomtensor_examples_originalmodellongtime}, \ref{fig:numericalsol_nonhomtensor_examples_originalmodelrealtensor} one can see bifurcations in the solution pattern for certain time steps. More precisely, there exist points where two roughly parallel lines  
%in the form of the letter u 
merge with a third roughly parallel line from the other side. These patterns are in the form of the letter `Y'. The evolution of one of these bifurcations is shown in Figure \ref{fig:Example1OriginalModelBifurcations} for the underlying tensor field in Example 6 in Figure \ref{fig:numericalsol_nonhomtensor_examples_originalmodelrealtensor}. Note that all these lines are aligned along the lines of smallest stress $s$ of the tensor field and these bifurcations move towards the two neighboring lines 
%in $u$-form 
over time. This behavior can be explained by attraction forces along the lines of largest stress over medium range distances, i.e.\ as soon as the distance between the particles along the lines of largest stress $l$ is small enough they attract each other. In particular, the particles close to the bifurcation  on the two neighboring lines 
are the first ones to `feel' the attraction force along $l$ and the two roughly parallel lines start merging close to the bifurcation. Hence, the single line on the other side of the bifurcation gets longer over time and the bifurcation moves towards the two parallel lines. While the two roughly parallel lines  get shorter over time until they are finally completely merged, resulting in one single line. 
Since the movement of the particles is mainly along $l$ there is a different particle at the bifurcation at each time step. While the particles on the line in the middle roughly remain at the same position apart from realigning along the lines of smallest stress $s$. 
This realignment along $s$ is due to  the additional number of  particles which are aligned along one single line after the merging, as well as due to the repulsive forces along $s$  spreading the particles to make use of the space along $s$ and to avoid high particles densities after merging.

\subsection{Variable ridge distances}
\subsubsection{Motivation for a new model}

%\subsubsection{Adaptation of the force coefficients}
The results in Section \ref{sec:realfingerprintpatterns} illustrate that it is possible to simulate realistic fingerprints with the Kücken-Champod model \eqref{eq:particlemodel}. 
As seen in the figures, there is some variability in ridge distances and in view of realistic biometric applications, it is of great interest to control them.
%For biometric applications, however, it is also of great interest to be able to vary the distances between the ridges. 
Note that the total force $F$ in \eqref{eq:totalforce}, given by the sum of  repulsion and attraction force $F_R$ and $F_A$ of the form \eqref{eq:repulsionforce} and \eqref{eq:attractionforce}, respectively, can be rewritten as
\begin{align}\label{eq:forceformorig}
F(d(x_j,x_k),T(x_j))=\left[ \chi f_A(|d|)+f_R(|d|)\right] (s\cdot d) s+\left[ f_A(|d|)+f_R(|d|)\right]  (l\cdot d)l
\end{align} 
by using the definition of the tensor field $T$ in \eqref{eq:tensorfield} and the definition of the distance vector $d(x_j,x_k)=x_j-x_k\in\R^2$. The coefficient functions of the  repulsion and attraction forces  \eqref{eq:repulsionforcemodel} and \eqref{eq:attractionforcemodel}, respectively, are plotted along $s$ and $l$ for the parameters in \eqref{eq:parametervaluesfingerprintsnew} in Figure \subref*{fig:forcesmodelnewparameter}. In particular, this motivates us to consider interaction forces of the form \eqref{eq:totalforcenew}.

We are interested in rescaling the forces now to vary the distances between the fingerprint lines, i.e.\ we consider $F(\eta d(x_j,x_k),T(x_j))$ where $\eta>0$ is the rescaling factor. For $\eta=1$ we recover the same solution patterns as in Section \ref{sec:realfingerprintpatterns}, while the distances between the fingerprint lines become larger for $\eta\in(0,1)$ and smaller for $\eta>1$. Note that the force coefficient $f_A+f_R$ along $l$ is repulsive over long distances. For $\eta=1$, the case that has been considered so far, this is fine for the given parameters in \eqref{eq:parametervaluesfingerprintsnew}. For $\eta>1$, however, the scaling results in repulsive interaction forces along $l$ for particles with shorter distances between each other. Besides, short-range forces have a stronger impact on the interactions. Hence, these short-range repulsive interaction forces prevent the accumulation of particles along $l$, resulting in several clusters. Note that the forces along $s$ are purely repulsive so that rescaling by any $\eta$ does not change the nature of the forces. 

\begin{figure}[htbp]
	\centering
	\subfloat[Normal scaling]{\includegraphics[width=0.45\textwidth]{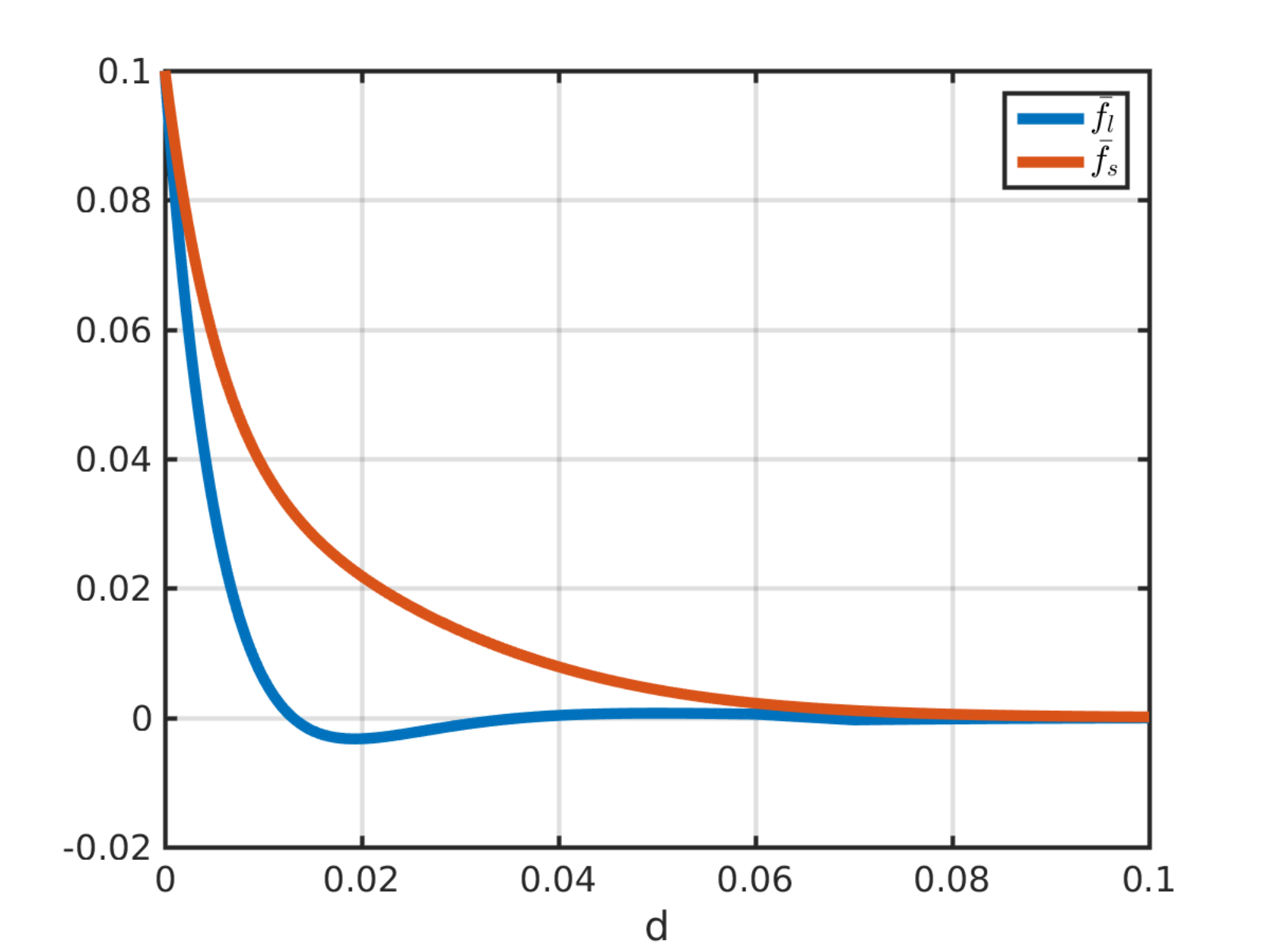}}\hfill
	\subfloat[Zoom]{\includegraphics[width=0.45\textwidth]{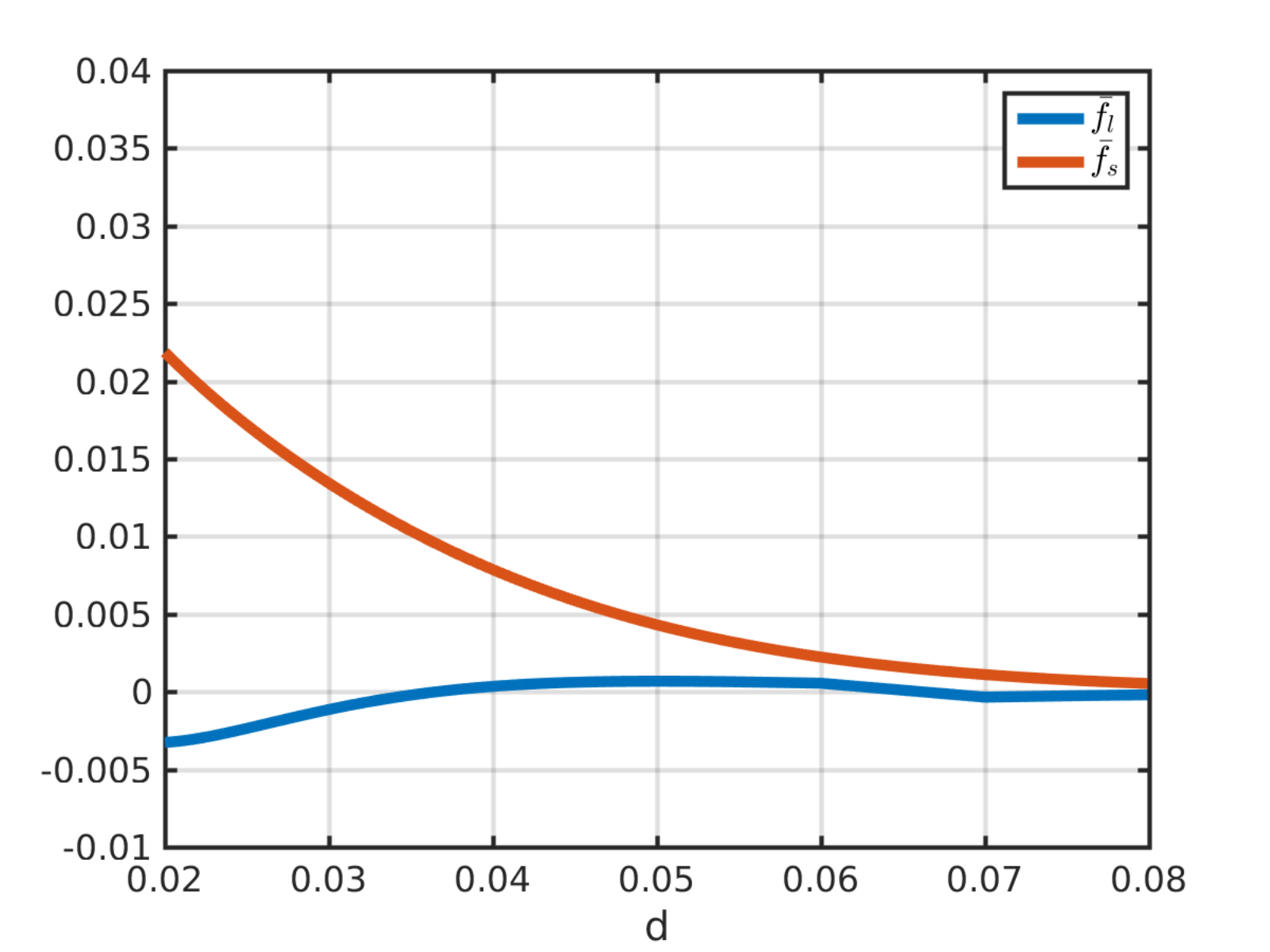}}\hfill
	\caption{Total force coefficients $\bar{f}_l$ and $\bar{f}_s$, defined in \eqref{eq:forcepiecewiseridgedistancel} and \eqref{eq:forcepiecewiseridgedistances} respectively, for interaction forces of the form \eqref{eq:totalforcenew} and parameter values \eqref{eq:parametervaluesfingerprintsnew}}\label{fig:forcesnonsmoothrescaling}
\end{figure}

\begin{figure}[htbp]
	\centering
	\subfloat[$\eta=0.8$]{\includegraphics[width=0.32\textwidth]{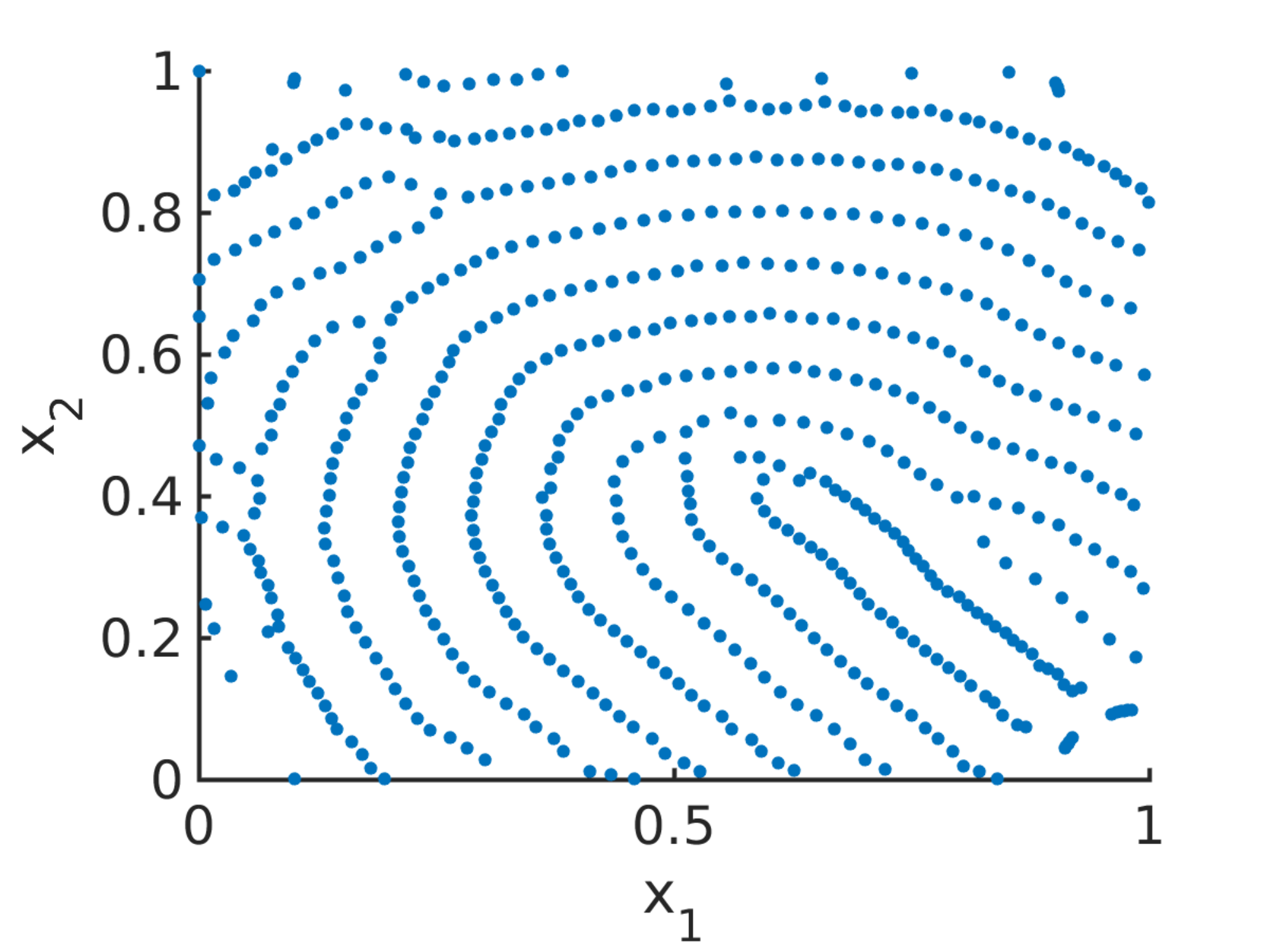}}\hfill
	\subfloat[$\eta=1.0$]{\includegraphics[width=0.32\textwidth]{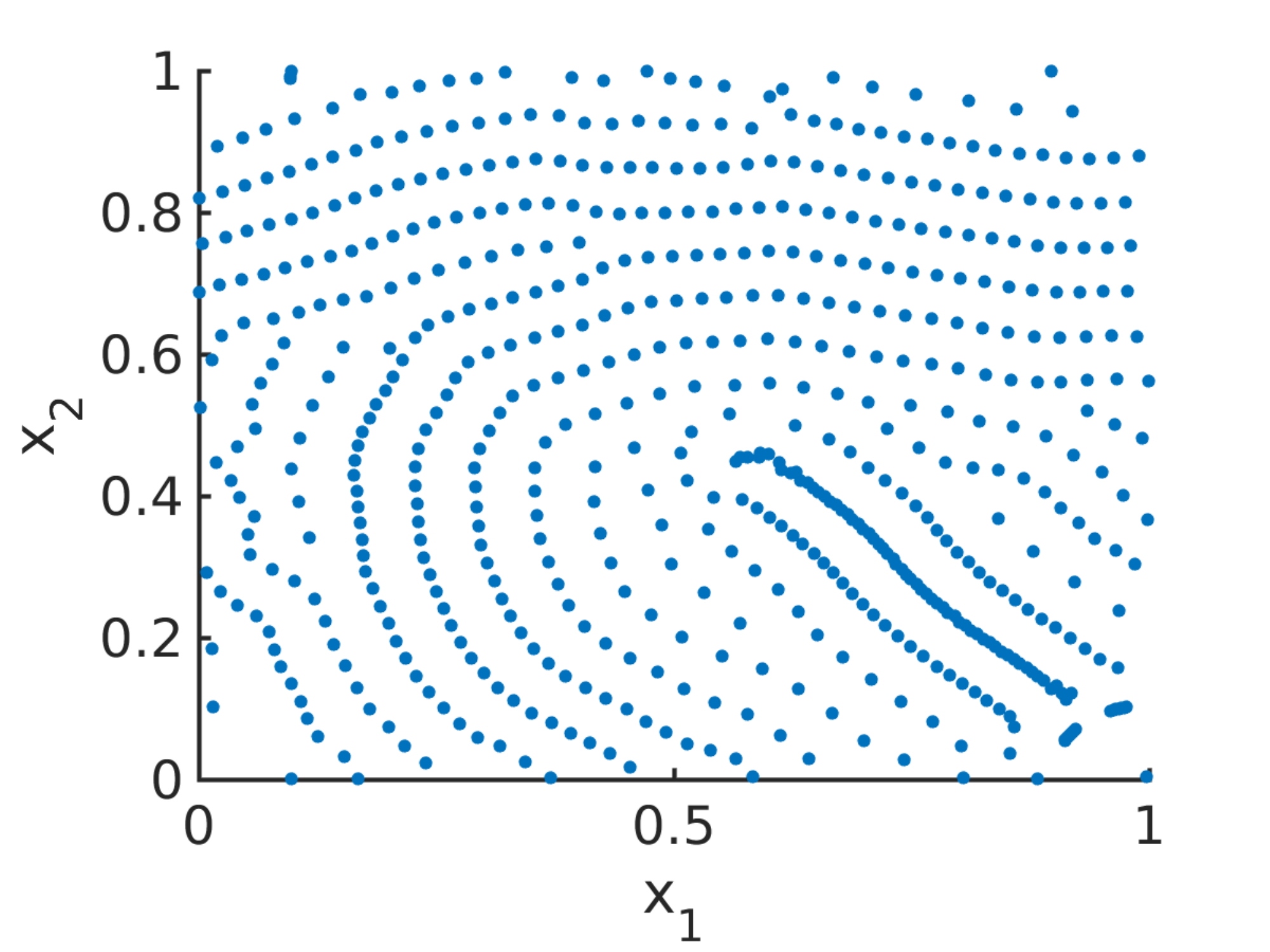}}\hfill
	\subfloat[$\eta=1.2$]{\includegraphics[width=0.32\textwidth]{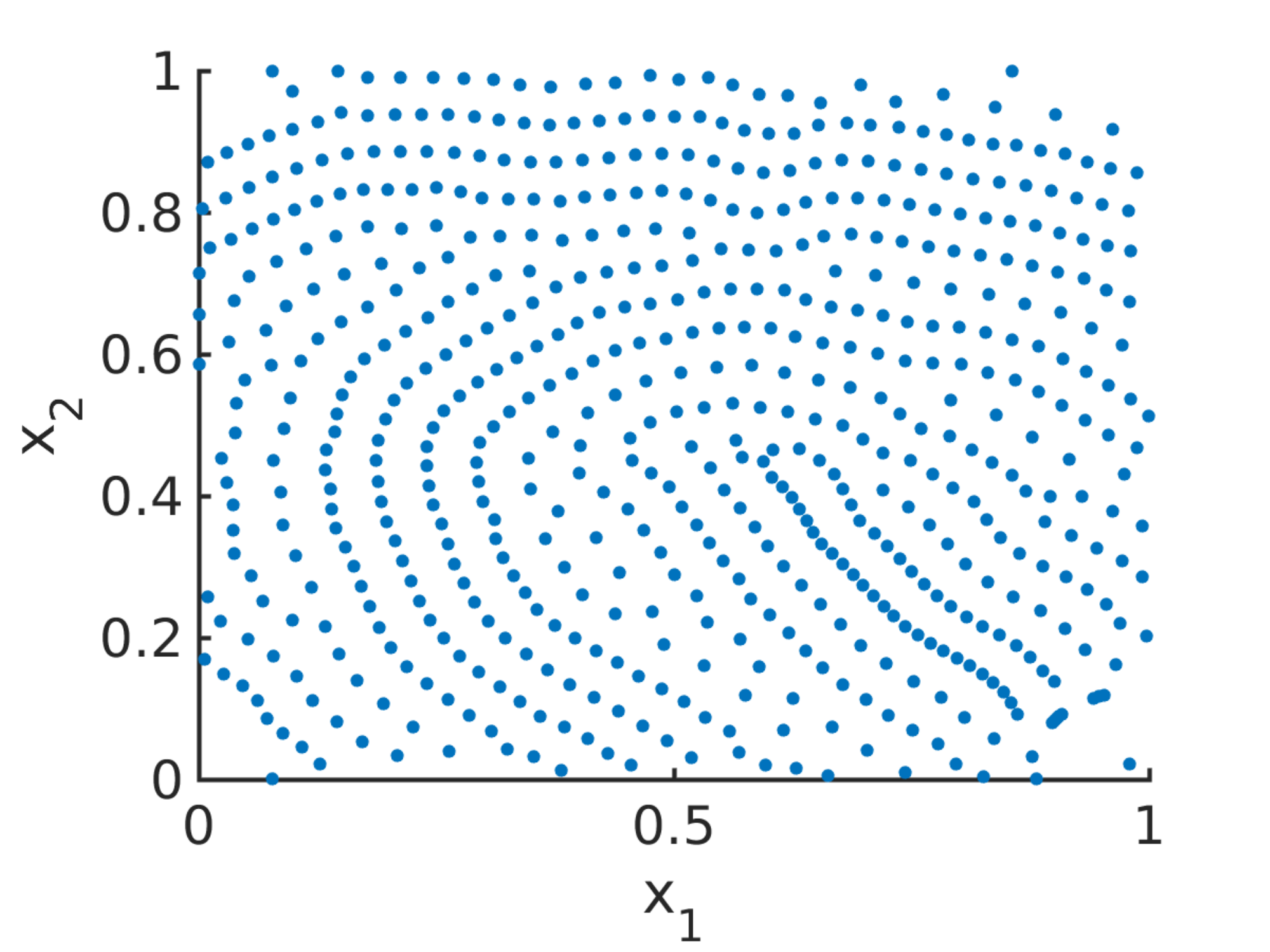}}\hfill
	\caption{Stationary solution to the interaction model \eqref{eq:particlemodel} for interaction forces of the form \eqref{eq:totalforcenew}, force coefficients \eqref{eq:forcepiecewiseridgedistancel}, \eqref{eq:forcepiecewiseridgedistances}, parameter values \eqref{eq:parametervaluesfingerprintsnew},  the realistic tensor field  $T=T(x)$ in Figure \ref{fig:realtensorfield} and $N=2400$ particles initially  distributed uniformly at random}\label{fig:stationarypiecewisesmooth}
\end{figure}

In order to prevent this behavior and to obtain an interaction model that can be used for different rescalings, the forces need to be changed slightly so that we have very small attractive forces along $l$ for $\eta=1$.  This does not influence the pattern formation for $\eta=1$, but for rescaling by $\eta>1$ we can obtain the desired line patterns with smaller distances between each other. In order to achieve this, we consider a straight-forward approach first. We consider two cutoffs $c_1$ and $c_2$ and define the adapted force $F$ piece-wise such that for $|d|<c_1$ the force $F$ is of the form \eqref{eq:forceformorig} as before while for $|d|>c_2$ we consider an attraction force  tending to zero as $d\to \infty$. To obtain a continuous force we consider a linear interpolation of the force on $[c_1,c_2]$. 
%Note that the force in the third case has to be attractive and tend to zero as $d\to \infty$.
Setting $$f(|d|,\chi):=\chi f_A(|d|)+f_R(|d|)$$ we consider the force coefficients $\bar{f}_s$ and $\bar{f}_l$ for interaction forces of the form \eqref{eq:totalforcenew} where the force coefficients are defined as
\begin{align}\label{eq:forcepiecewiseridgedistancel}
\begin{split}
\bar{f}_l(d)=\begin{cases}
f(|d|,1)   & |d|<c_1\\
f(c_1,1)+\frac{|d|-c_1}{c_2-c_1}\bl -f(c_2,1)-f(c_1,1)\br & |d|\in[c_1,c_2]\\
-f(|d|,1) & |d|>c_2
\end{cases}
\end{split}
\end{align}
and
\begin{align}\label{eq:forcepiecewiseridgedistances}
\bar{f}_s(d)=f(|d|,\chi).
\end{align}
% following force $F$ which is of the above form:
%\begin{align}\label{eq:forcepiecewiseridgedistance}
%\begin{split}
%&F(d(x_j,x_k),T(x_j))\\&=\begin{cases}
%f(|d|,\chi)(s\cdot d) s+f(|d|,1)  (l\cdot d)l & |d|<c_1\\
%f(|d|,\chi) (s\cdot d) s+\bl f(c_1,1)+100(|d|-c_1)\bl -f(c_2,1)-f(c_1,1)\br\br (l\cdot d)l& |d|\in[c_1,c_2]\\
%f(|d|,\chi) (s\cdot d) s-f(|d|,1)  (l\cdot d)l & |d|>c_2.
%\end{cases}
%\end{split}
%\end{align}
%\begin{align*}
%F(d(x_j,x_k),T(x_j))=\begin{cases}
%\left[ \chi f_A(|d|)+f_R(|d|)\right] (s\cdot d) s+\left[ f_A(|d|)+f_R(|d|)\right]  (l\cdot d)l & |d|<c_1\\
%\text{linear interpolation to have continuous force}& |d|\in[c_1,c_2]\\
%\left[ \chi f_A(|d|)+f_R(|d|)\right] (s\cdot d) s-\left[ f_A(|d|)+f_R(|d|)\right]  (l\cdot d)l & |d|>c_2
%\end{cases}
%\end{align*}
Here, we consider the parameter values $c_1=0.06,c_2=0.07$ and the parameters in the force coefficients  \eqref{eq:repulsionforcemodel}, \eqref{eq:attractionforcemodel} are given by \eqref{eq:parametervaluesfingerprintsnew}. The force coefficient $f_l$ along $l$ for $|d|>c_2$ is obtained by multiplying the original force along $l$ by $-1$. This is based on the fact that the  force coefficient $f(d,1)$ is repulsive for large distances along $l$ for the parameters in \eqref{eq:parametervaluesfingerprintsnew}. In Figure \ref{fig:forcesnonsmoothrescaling} the force coefficients $\bar{f}_l$ and $\bar{f}_s$ in \eqref{eq:forcepiecewiseridgedistancel} and \eqref{eq:forcepiecewiseridgedistances}, respectively, are shown. In particular, the piecewise definition of $\bar{f}_l$ only has a small influence of the form.  In Figure \ref{fig:stationarypiecewisesmooth}, the stationary solution to the particle model \eqref{eq:particlemodel} for  interaction forces of the form \eqref{eq:totalforcenew}, force coefficients \eqref{eq:forcepiecewiseridgedistancel}, \eqref{eq:forcepiecewiseridgedistances}, parameter values \eqref{eq:parametervaluesfingerprintsnew}, the underlying tensor field $T=T(x)$ in Figure \ref{fig:realtensorfield}  and different rescaling factors $\eta$ is shown and one can clearly see that $\eta>1$ leads to smaller ridge distances whereas $\eta<1$ results in larger ridge distances. In particular, the interaction model \eqref{eq:particlemodel} with interaction forces of the form \eqref{eq:totalforcenew} and force coefficients in \eqref{eq:forcepiecewiseridgedistancel} and \eqref{eq:forcepiecewiseridgedistances} can be used to simulate fingerprints with variable ridge distances. Due to the smaller distances between the fingerprint lines  for $\eta=1.2$ this leads to a larger number of fingerprint lines on the given domain. Due to this increased number of lines it is desirable to run simulations with larger numbers of particles. However, particle simulations can only be applied efficiently as long as the total particle number is not too large. In order to solve this remedy one can introduce the density $\rho=\rho(t,x)$ associated with the particle positions and consider the associated macroscopic model 
\begin{align}\label{eq:macroscopiceq}
\begin{split}
\partial_t \rho(t,x)+\nabla_x\cdot \left[ \rho(t,x)\bl F\bl\cdot,T(x)\br \ast \rho(t,\cdot)\br\bl x\br\right]=0\qquad \text{in }\R^2\times \R_+.
\end{split}
\end{align}
In future work,  advanced numerical methods for solving the macroscopic model \eqref{eq:macroscopiceq} with anisotropic interaction forces could be developed for simulating fingerprint patterns.

%To vary the distances between the ridges we consider rescaled forces $F(d(x_j,x_k)*\eta,T(x_j))$ where $\eta$ is the rescaling factor. In particular, the larger $\eta$ the smaller the distance between the 

\subsubsection{A bio-inspired model for simulating stationary fingerprints with variable ridge distances}
\begin{figure}[htbp]
	\centering
	\includegraphics[width=0.45\textwidth]{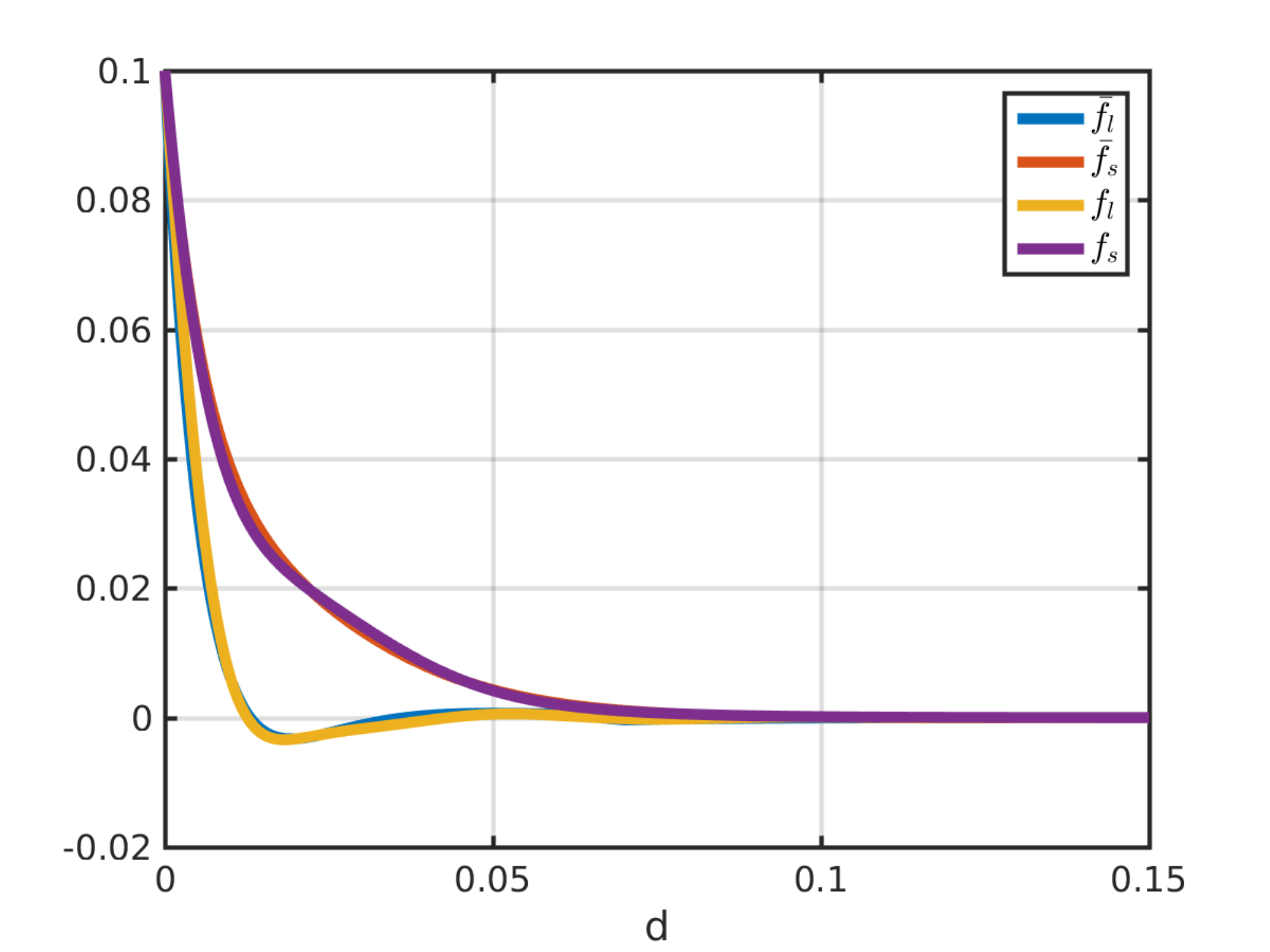}
	\caption{Coefficients $f_l$ and $f_s$ in \eqref{eq:forcecoeffharmonic}  for parameter values in  \eqref{eq:parameterharmonic} as well as piecewise defined coefficients $\bar{f}_l$ and $\bar{f}_s$ in \eqref{eq:forcepiecewiseridgedistancel},\eqref{eq:forcepiecewiseridgedistances}}\label{fig:forcesharmonicoscillatorcomparison}
\end{figure}

\begin{figure}[htbp]
	\centering
	\subfloat[$\eta=0.6$]{\includegraphics[width=0.24\textwidth]{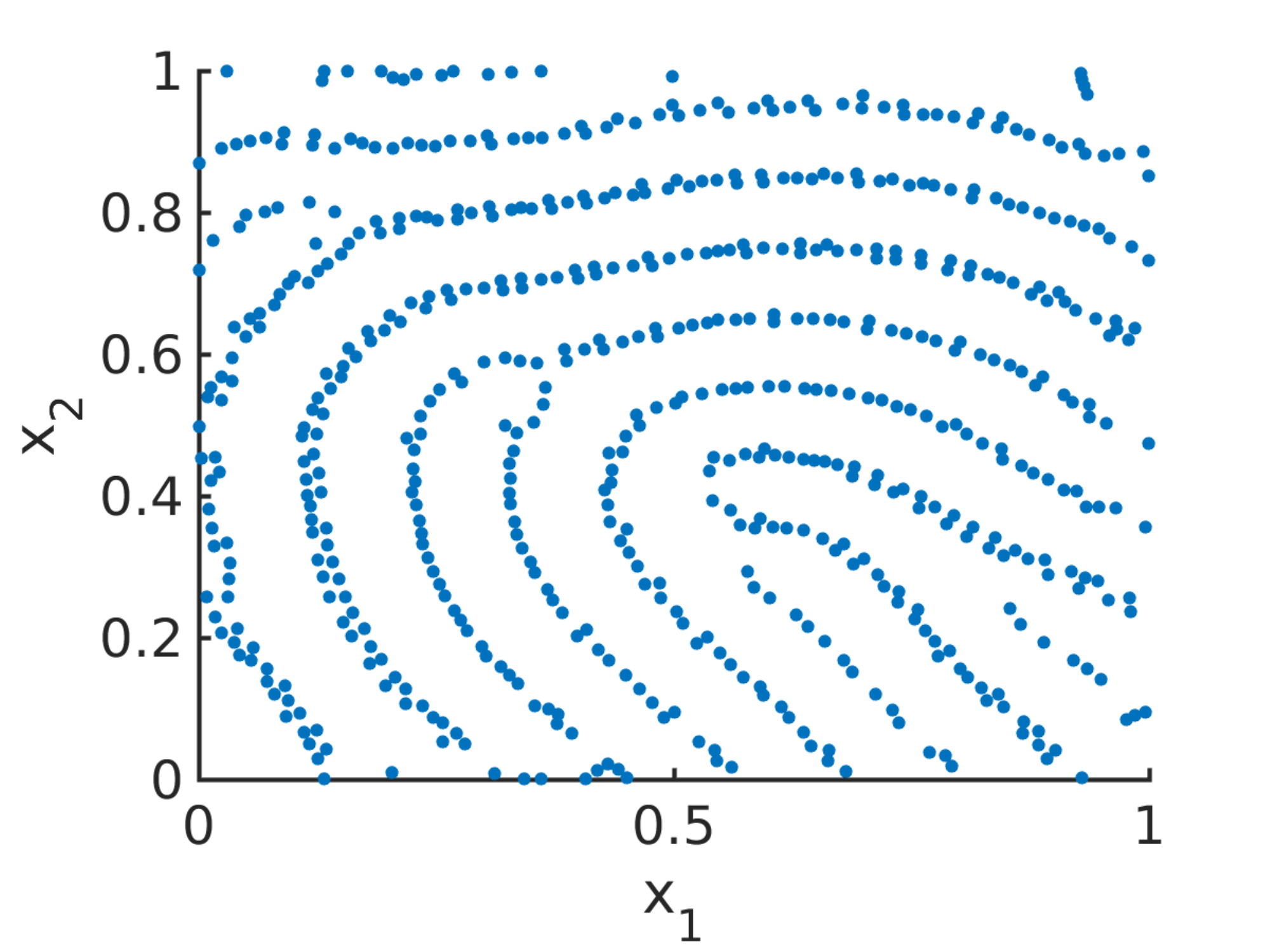}}\hfill
	\subfloat[$\eta=0.8$]{\includegraphics[width=0.24\textwidth]{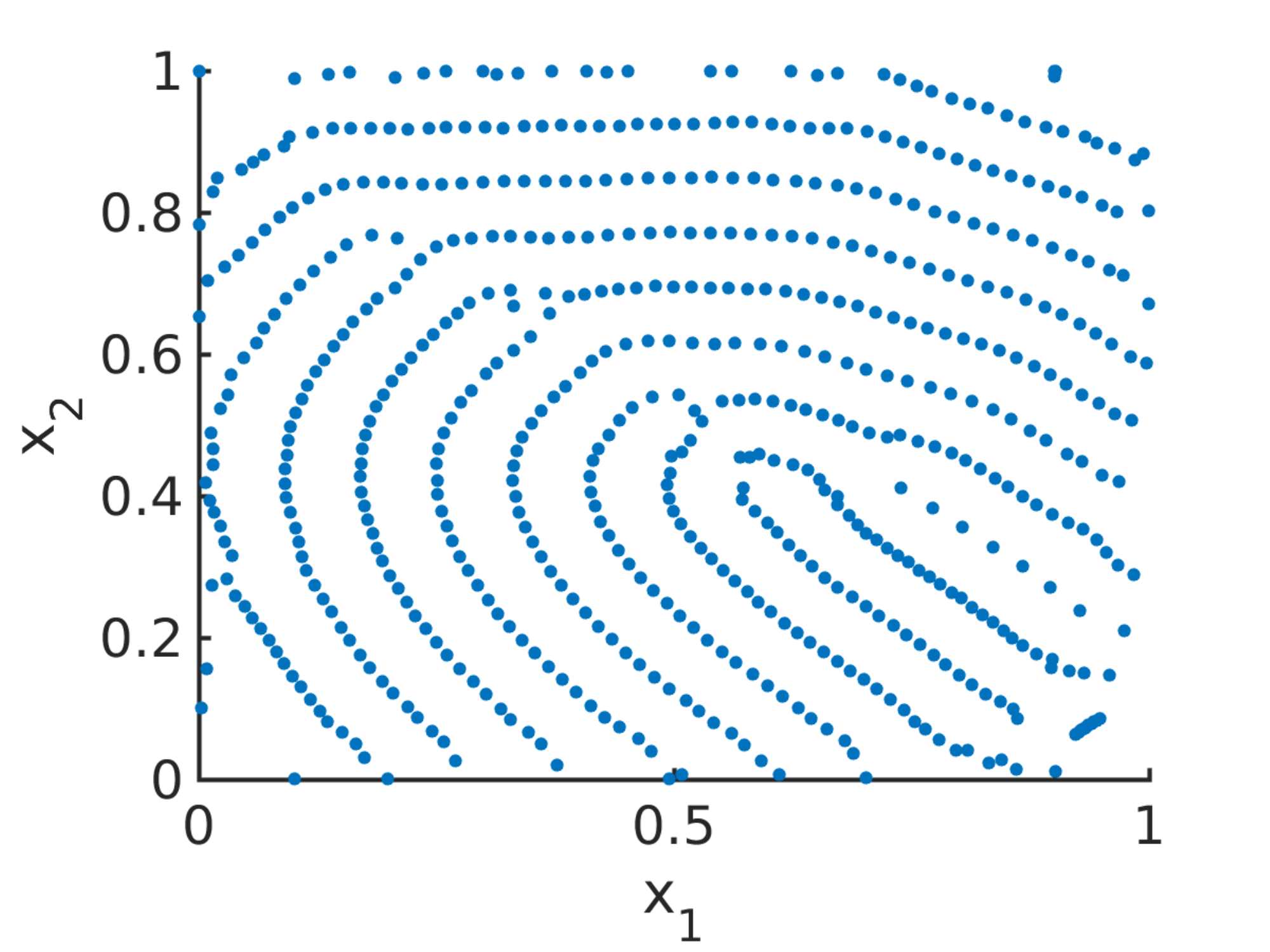}}\hfill
	\subfloat[$\eta=1.0$]{\includegraphics[width=0.24\textwidth]{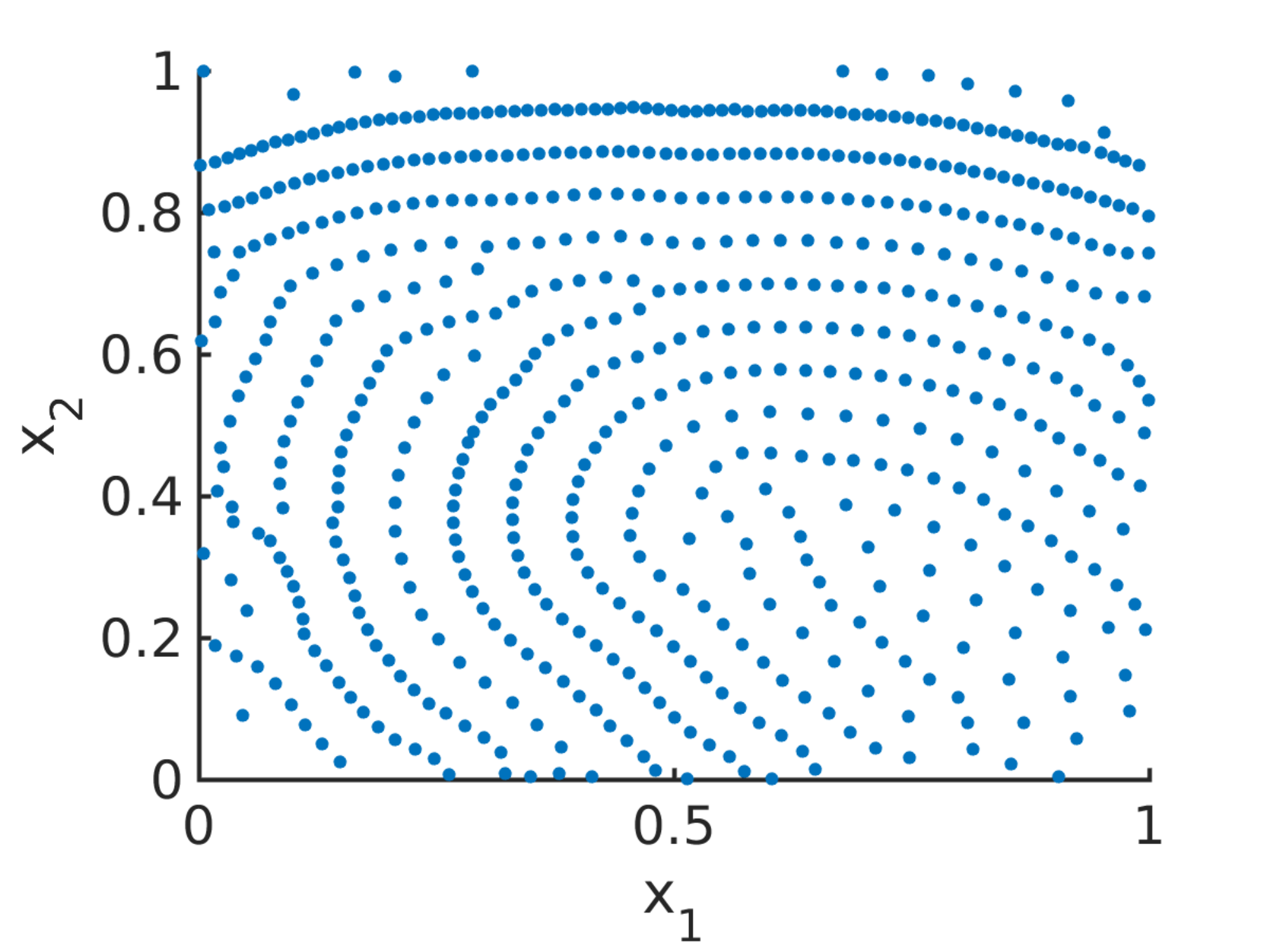}}\hfill
	\subfloat[$\eta=1.2$]{\includegraphics[width=0.24\textwidth]{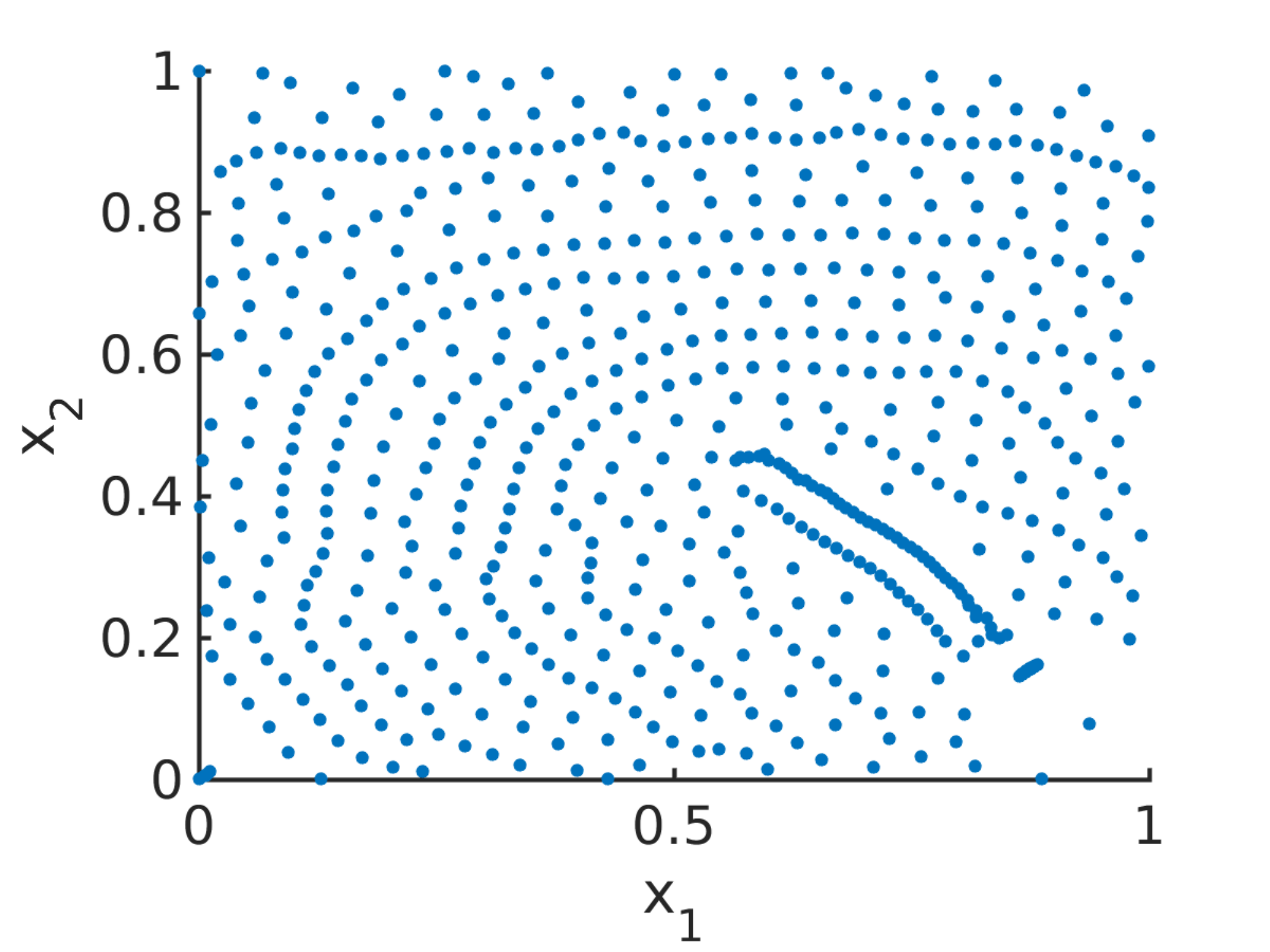}}\hfill
	\caption{Stationary solution to the interaction model \eqref{eq:particlemodel} for interaction forces of the form \eqref{eq:totalforcenew}, force coefficients \eqref{eq:forcecoeffharmonic}, parameter values \eqref{eq:parameterharmonic},  the realistic tensor field  $T=T(x)$ in Figure \ref{fig:realtensorfield} and $N=2400$ particles initially distributed uniformly at random}\label{fig:stationaryharmonic}
\end{figure}

In this section, we consider interaction forces of the form \eqref{eq:totalforcenew} as before with the aim of simulating fingerprints with variable ridge distances based on a bio-inspired approach. The coefficient functions $f_l$ and $f_s$ in \eqref{eq:forcepiecewiseridgedistancel} and \eqref{eq:forcepiecewiseridgedistances}, respectively, are defined piecewise and it is desirable to obtain a closed form for the coefficient functions. As before  we consider  exponentially decaying forces  describing that short-range interactions between the particles  are much stronger than long-range interactions. Since the forces are repulsive and attractive on different regimes, this interplay between repulsion and attraction forces can be regarded  as oscillations. Motivated by this, we  model the force coefficients $f_s$ and $f_l$ in \eqref{eq:totalforcenew} as solutions to a damped harmonic oscillator. Note that harmonic oscillators are a common modeling approach in cell biology and the force coefficients $f_l, f_s$ are given by \eqref{eq:forcecoeffharmonic} and are shown in Figure \ref{fig:forcesharmonicoscillatorcomparison} for the parameters in \eqref{eq:parameterharmonic} in comparison with the piecewise defined force coefficients $\bar{f}_l,\bar{f}_s$ for the parameters in \eqref{eq:parametervaluesfingerprintsnew}. Note that the parameters \eqref{eq:parameterharmonic} are chosen in such a way that the coefficient functions $f_l,f_s$ of the harmonic oscillator approximate the piecewise defined coefficient functions $\bar{f}_l,\bar{f}_s$ in \eqref{eq:forcepiecewiseridgedistancel},\eqref{eq:forcepiecewiseridgedistances}, respectively. In Figure \ref{fig:stationaryharmonic} the stationary patterns to \eqref{eq:particlemodel} for different rescaling factors $\eta$ are shown. As expected the larger the value of $\eta$ the smaller the distances between the fingerprint lines and the more lines occur. 

\subsubsection{Whole fingerprint simulations}
\begin{figure}[htbp]
	\centering
	\subfloat[Original]{\includegraphics[width=0.24\textwidth]{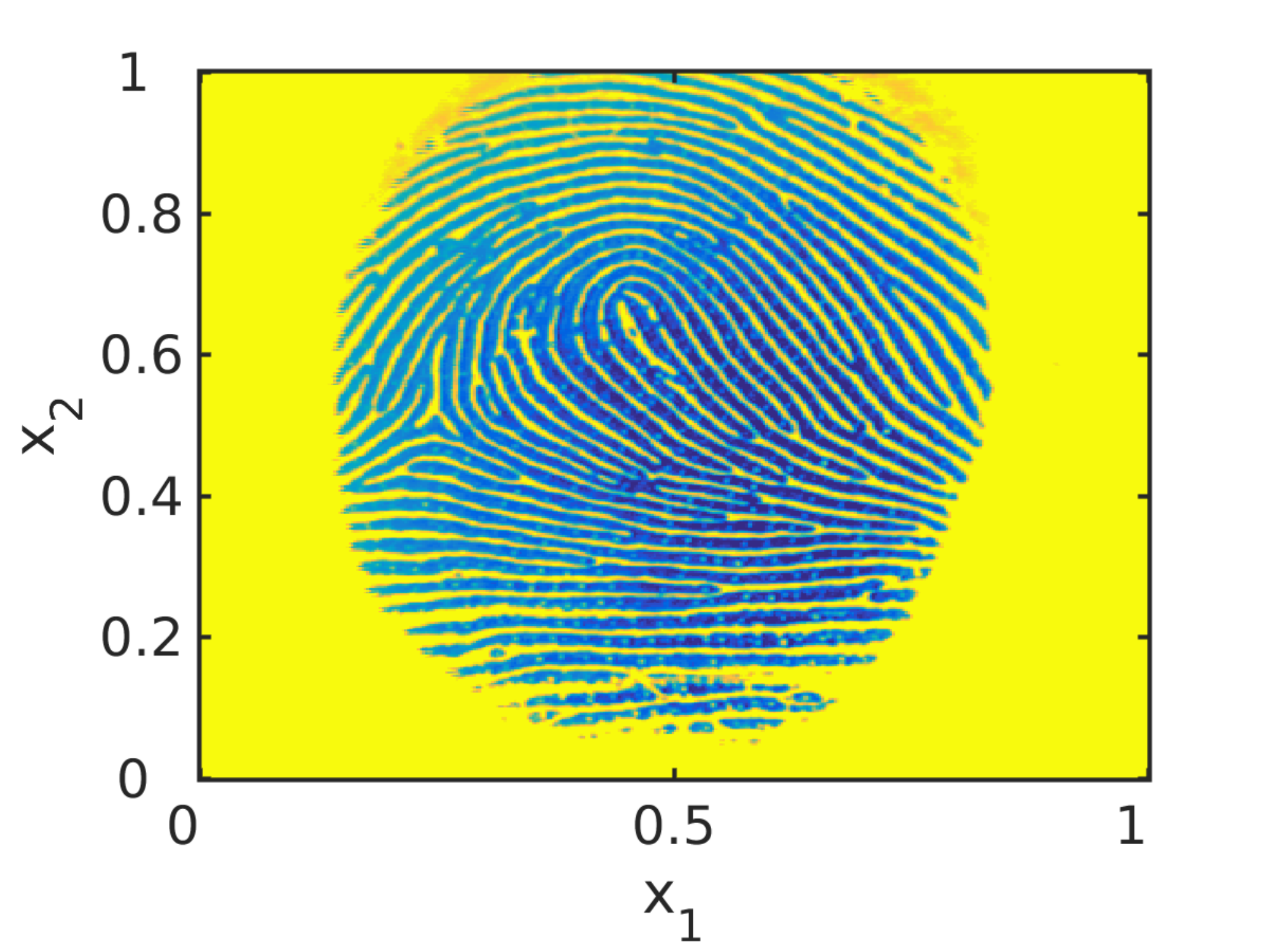}\label{fig:realfingerimage}}\hfill
	\subfloat[$\theta$ and original ]{\includegraphics[width=0.24\textwidth]{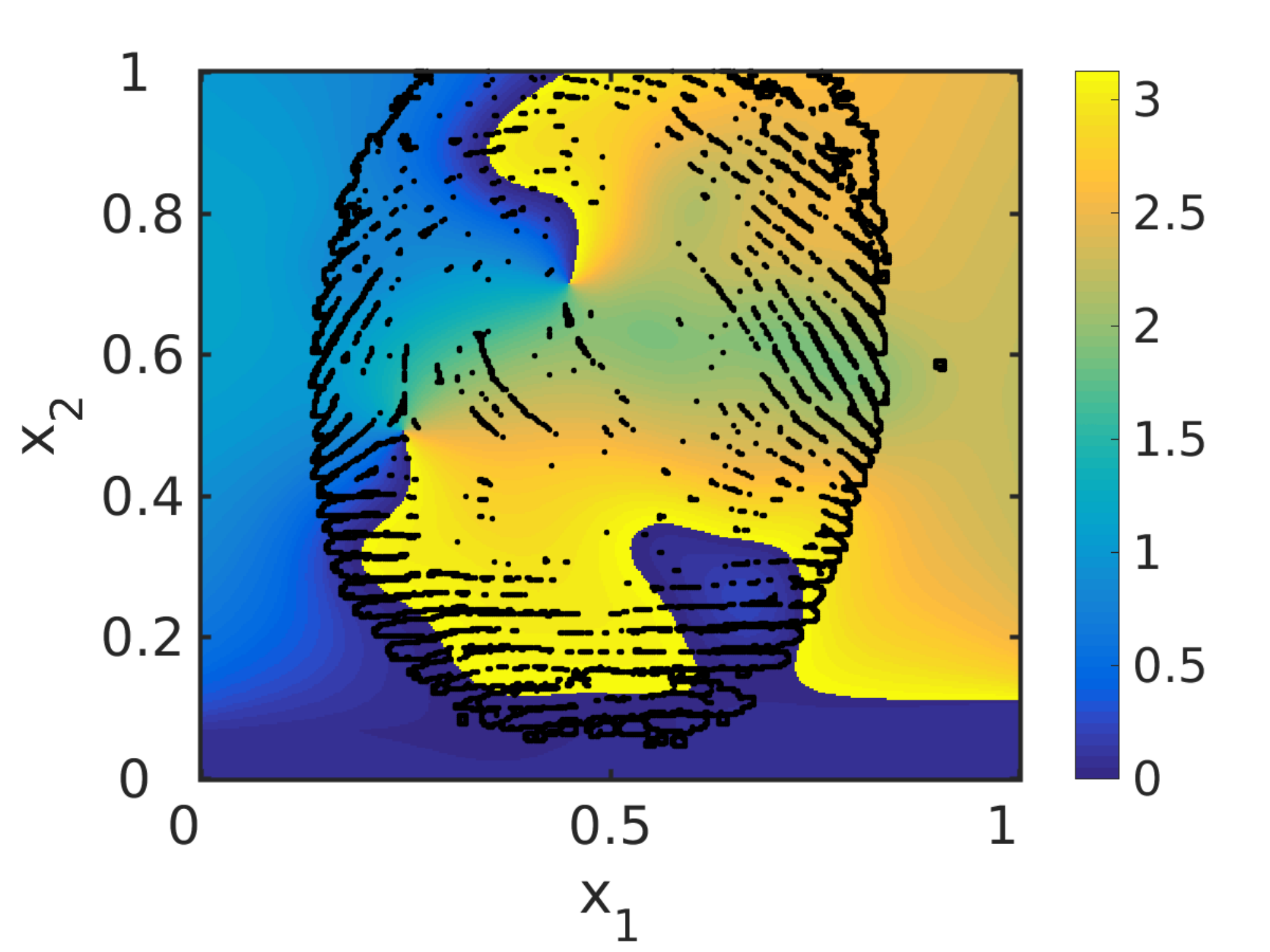}\label{fig:realfingerargmask}}\hfill
	\subfloat[$\theta$]{\includegraphics[width=0.24\textwidth]{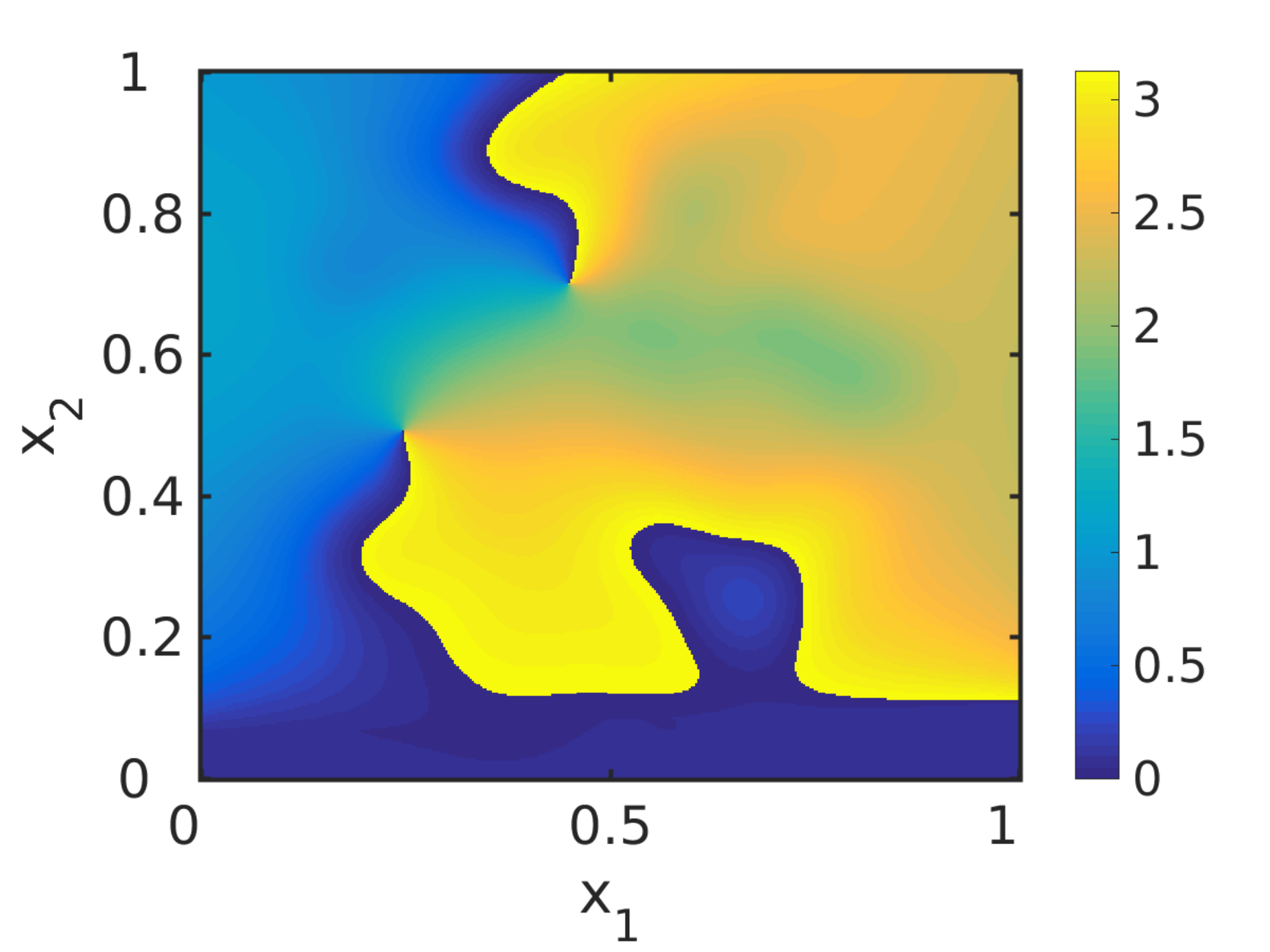}\label{fig:realfingerarg}}\hfill
	\subfloat[$s$ and original]{\includegraphics[width=0.24\textwidth]
		{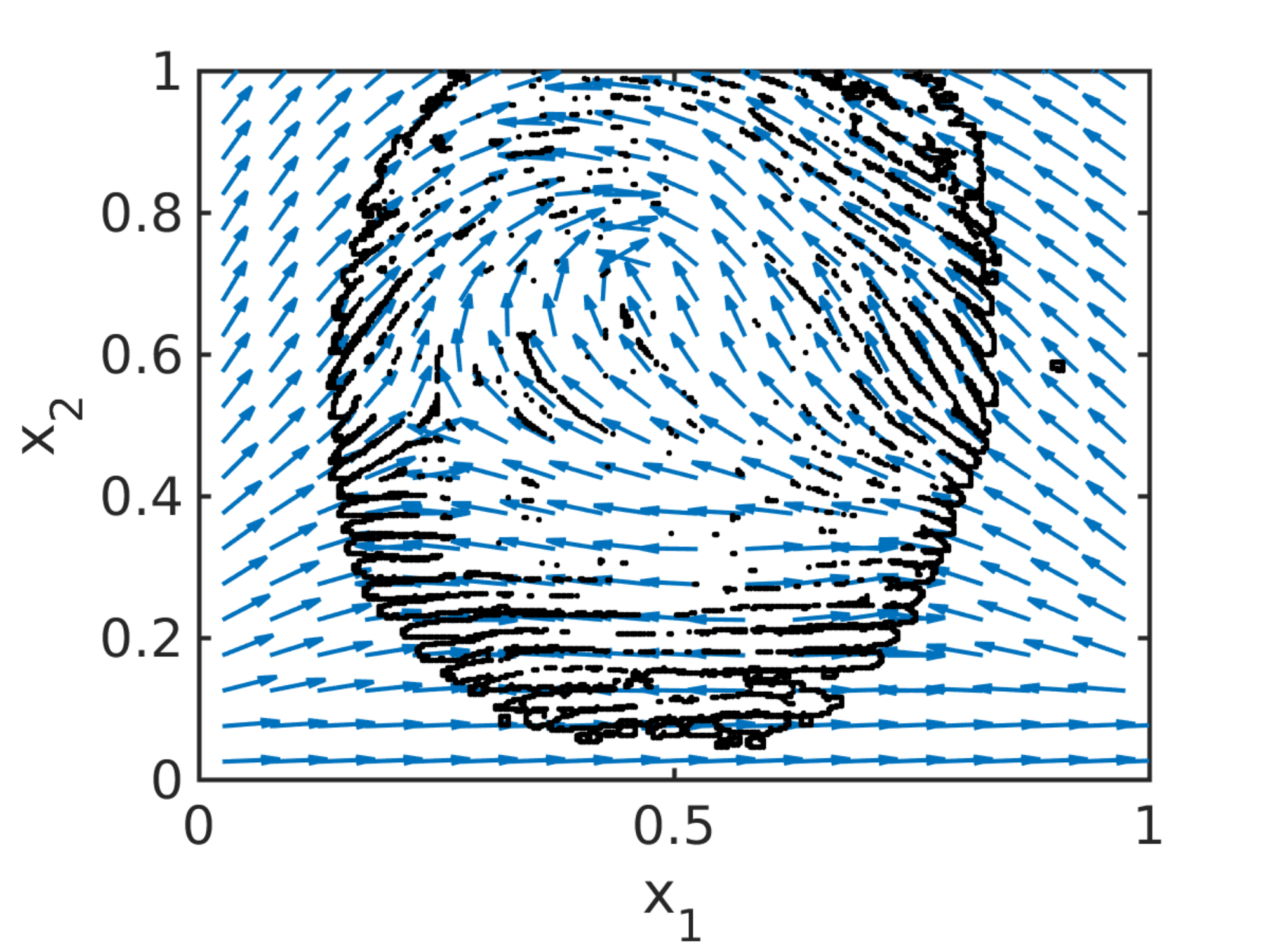}\label{fig:realfingertensormask}}\hfill
	\subfloat[$s$]{\includegraphics[width=0.24\textwidth]
		{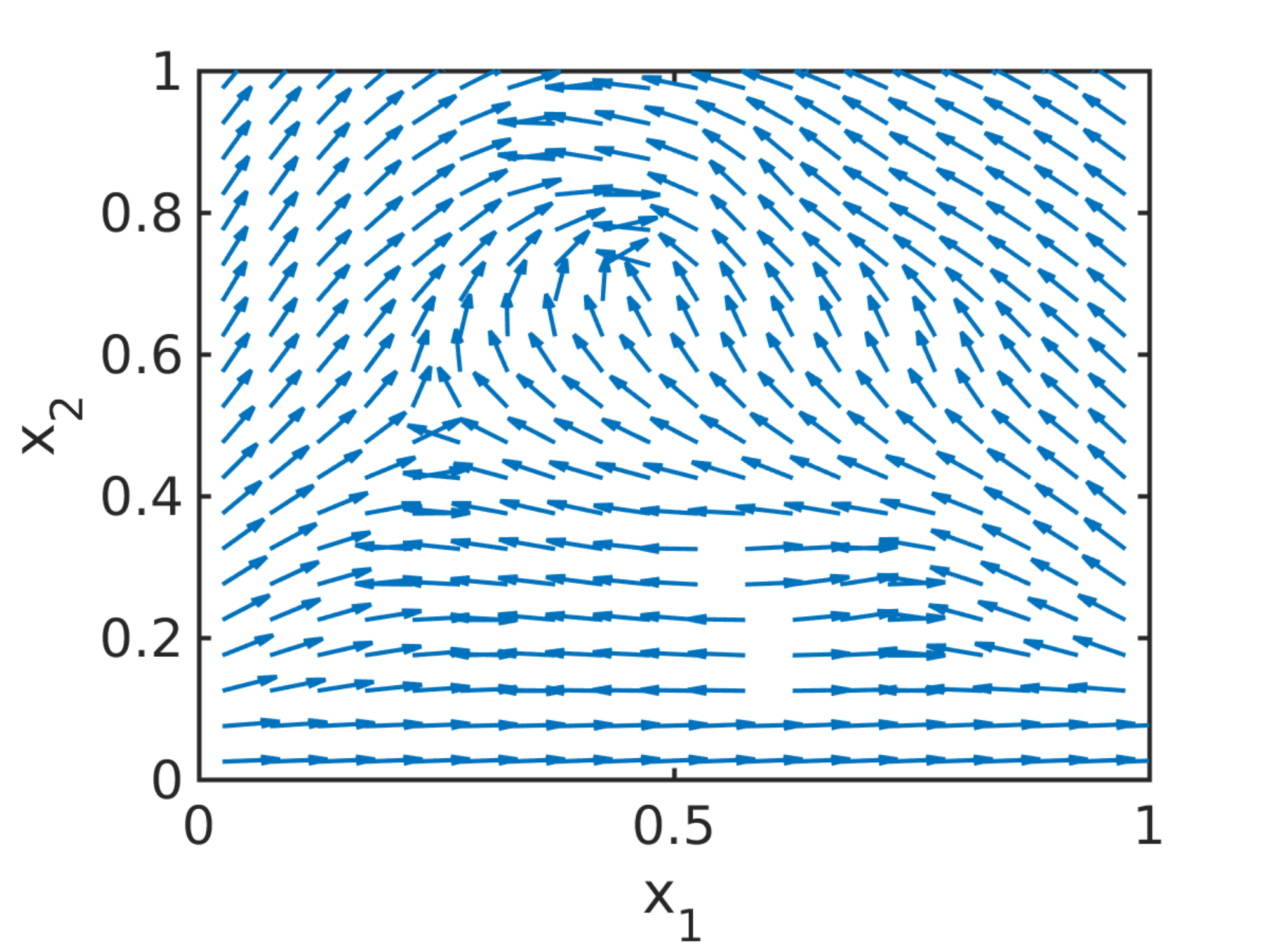}\label{fig:realfingertensor}}\hfill
	\subfloat[$\eta=0.6$]{\includegraphics[width=0.24\textwidth]{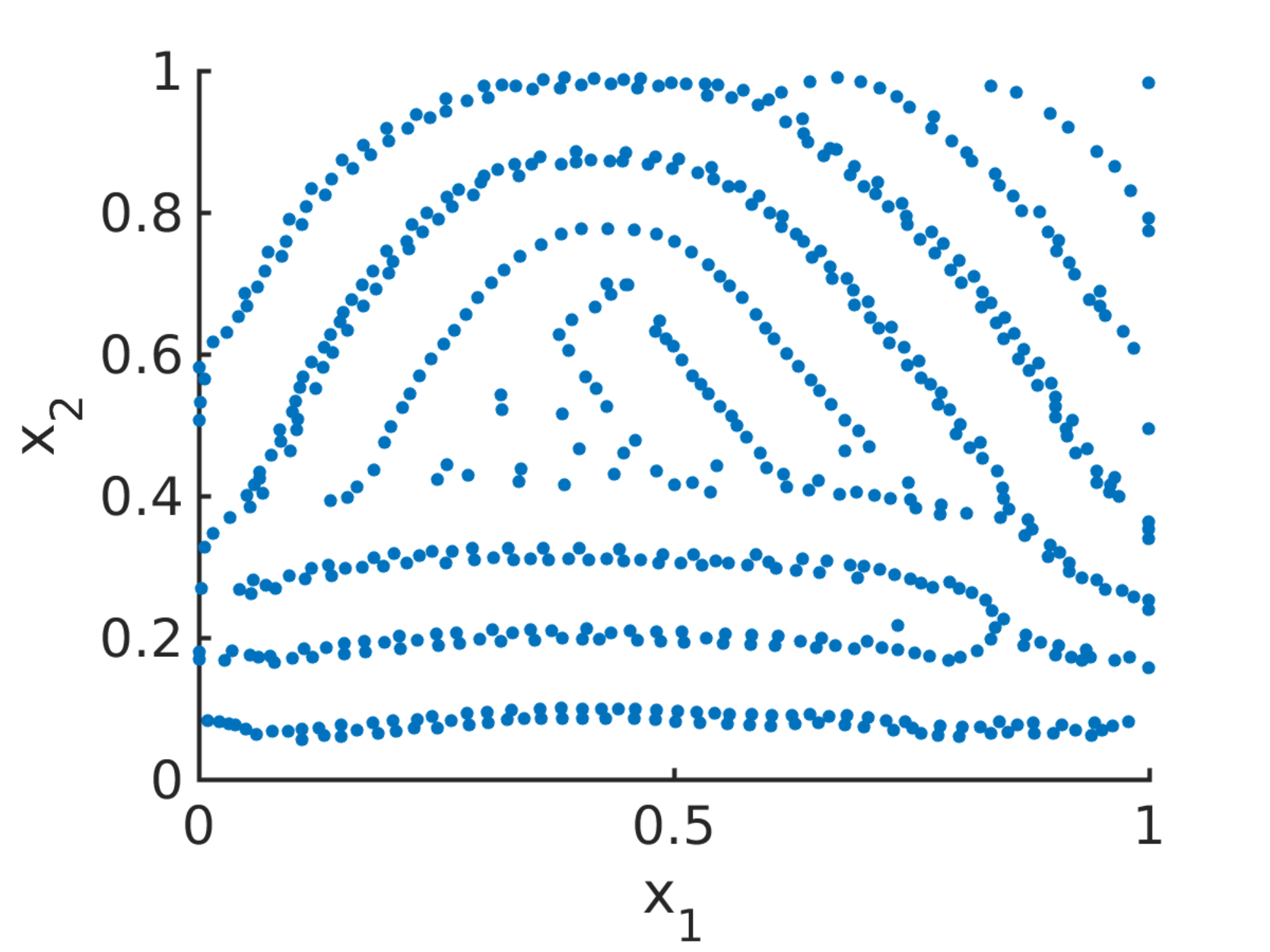}}\hfill
	\subfloat[$\eta=0.8$]{\includegraphics[width=0.24\textwidth]{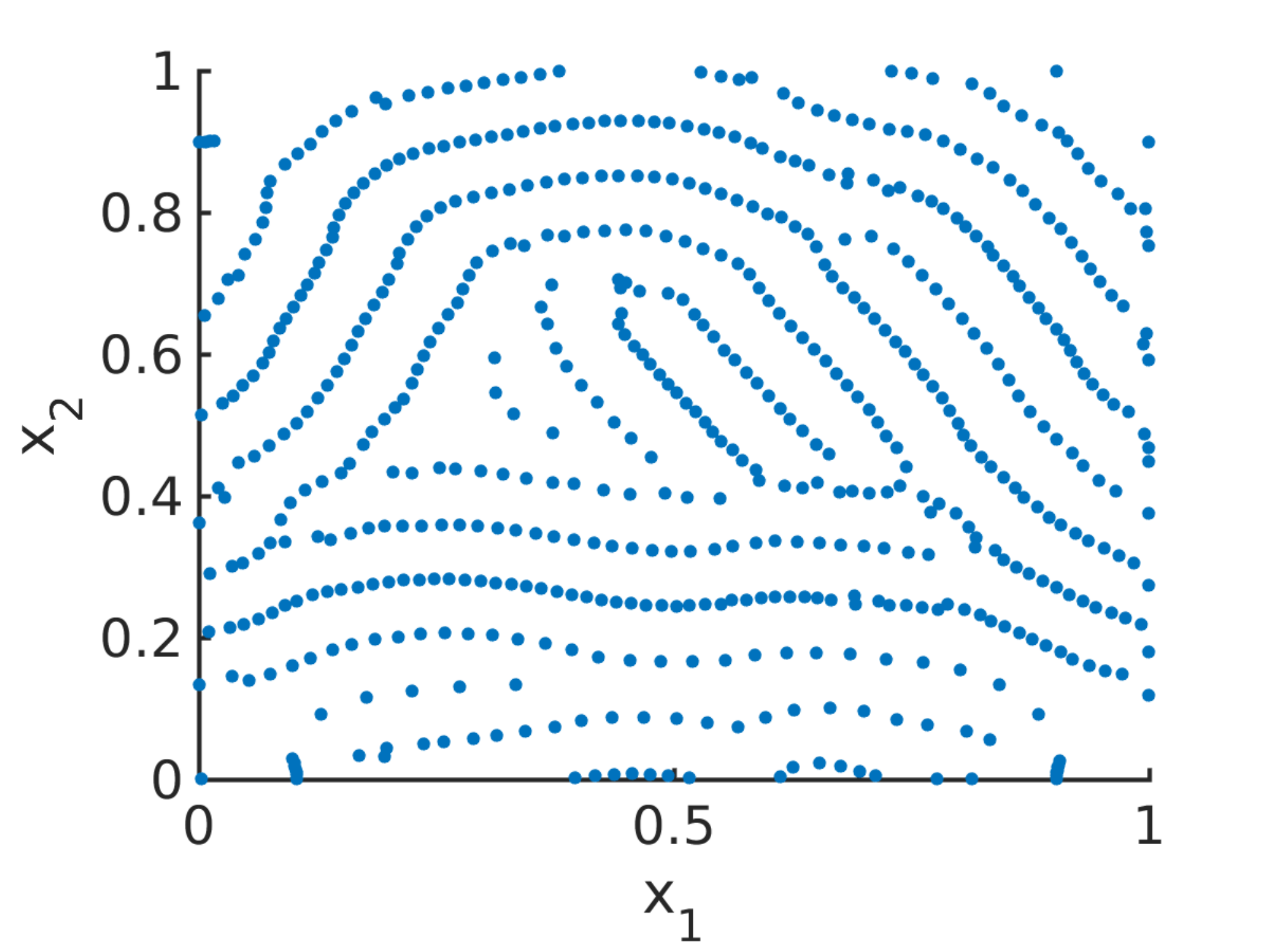}}\hfill
	\subfloat[$\eta=1.0$]{\includegraphics[width=0.24\textwidth]
		{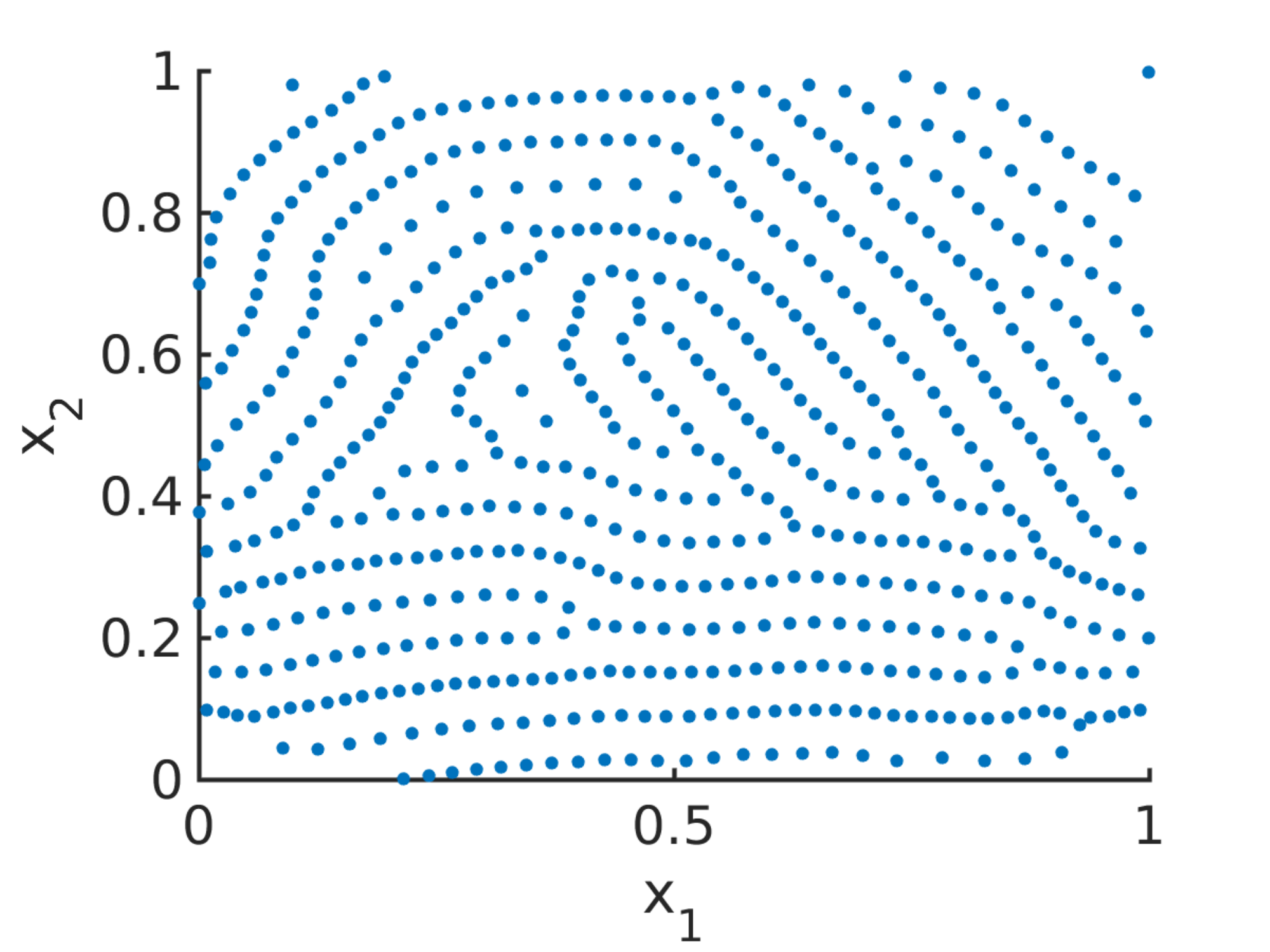}}\hfill
	\caption{Original fingerprint image, arguments and lines of smallest stress $s=s(x)$ for the reconstructed tensor field $T=T(x)$ with an overlying mask of the original fingerprint image in black, as well as stationary solution to the interaction model \eqref{eq:particlemodel} for interaction forces of the form \eqref{eq:totalforcenew}, force coefficients \eqref{eq:forcecoeffharmonic}, parameter values \eqref{eq:parameterharmonic} and $N=2400$ particles initially distributed  uniformly at random}\label{fig:stationaryharmonicrealfinger}
\end{figure}

In Figure \ref{fig:stationaryharmonicrealfinger} we construct tensor fields from real fingerprint data based on the methods discussed in Section \ref{sec:realfingerprintpatterns}. We consider a whole fingerprint image shown in Figure \subref*{fig:realfingerimage} and determine the underlying tensor field by estimating the arguments $\theta=\theta(x)$ for every $x\in\Omega$. Since we consider the domain $\Omega=\mathbb{T}^2$ we extend the tensor field via extrapolation from the original fingerprint image in Figure \subref*{fig:realfingerimage}, based on \cite{LineSensor}. In Figures \subref*{fig:realfingerargmask} and \subref*{fig:realfingerarg} the arguments $\theta=\theta(x)$ are shown and the arguments $\theta$ are overlayed by the mask of the original fingerprint in black in Figure \subref*{fig:realfingerargmask}. Since $s(x)$ is a unit vector and hence uniquely determined by its argument $\theta(x)$ we reconstruct the lines of smallest stress $s(x)$ as  $(\cos(\theta(x),\sin(\theta(x)))$ in Figures \subref*{fig:realfingertensormask} and \subref*{fig:realfingertensor}, and overlay the direction field $s$ by the original fingerprint image in black in Figure \subref*{fig:realfingertensormask}. We run simulations for these realistic tensor fields using our new bio-inspired model \eqref{eq:particlemodel} with interaction forces of the form \eqref{eq:totalforcenew}, force coefficients \eqref{eq:forcecoeffharmonic} inspired from harmonic oscillators and parameter values in \eqref{eq:parameterharmonic} for randomly uniformly distributed initial data and $N=2400$ particles. Note that the patterns are preserved over time.

\section*{Acknowledgments}
BD has been supported by the Leverhulme Trust research project grant `Novel discretisations for higher-order nonlinear PDE' (RPG-2015-69). SH acknowledges support from the Niedersachsen Vorab of the Volkswagen Foundation and the DFG Graduate Research School 2088. LMK has been supported by the UK Engineering and Physical Sciences Research Council (EPSRC) grant
EP/L016516/1. CBS acknowledges support from Leverhulme Trust project on `Breaking the non-convexity barrier', EPSRC grant Nr. EP/M00483X/1, the EPSRC Centre Nr. EP/N014588/1 and the Cantab Capital Institute for the Mathematics of Information.

\bibliographystyle{plain}
\bibliography{referencesfingerprints}
\end{document}